%% file: fibrmp.tex
\documentclass[onefignum,onetabnum]{siamart171218} 

\pdfoutput=1

\usepackage{amssymb,amsmath}

\usepackage{calc}
\usepackage[margin=1.in]{geometry}

\input texFiles/defs

\usepackage{tikz}

\input texFiles/trimFig.tex

\begin{document}


\title{A stable added-mass partitioned (AMP) algorithm 
       for elastic solids and incompressible flow: model problem analysis.\thanks{Submitted to the editors DATE.\funding{This research was supported by the National Science Foundation under grants DGE-1744655, DMS-1519934, and DMS-1818926, 
as well as by DOE contracts from the ASCR Applied Math Program.}}}

\author{D.~A.~Serino\thanks{Department of Mathematical Sciences, Rensselaer Polytechnic Institute, Troy, NY 12180, USA
(\email{serind@rpi.edu}, \email{banksj3@rpi.edu}, \email{henshw@rpi.edu}, \email{schwed@rpi.edu}).
}
\and J.~W.~Banks\footnotemark[2]
\and W.~D.~Henshaw\footnotemark[2]
\and D.~W.~Schwendeman\footnotemark[2]
}
\maketitle










\begin{abstract}
A stable added-mass partitioned (AMP) algorithm is developed for fluid-structure
interaction (FSI) problems involving viscous incompressible flow and compressible elastic-solids.
The AMP scheme remains stable and second-order accurate even when added-mass and 
added-damping effects are large.
The fluid is updated with an implicit-explicit (IMEX) fractional-step scheme whereby the velocity
is advanced in one step, treating the viscous terms implicitly, and the pressure is computed
in a second step. The AMP interface conditions for the fluid arise from the outgoing 
characteristic variables in the solid and are partitioned into a Robin (mixed) interface condition
for the pressure, and interface conditions for the velocity. The latter conditions include an 
impedance-weighted average between fluid and solid velocities using a fluid impedance of a 
special form. A similar impedance-weighted average is used to define interface values for 
the solid. 
The fluid impedance is defined using material and discretization parameters and follows 
from a careful analysis of the discretization of the governing equations and coupling 
conditions near the interface. 
A normal mode analysis is performed to show that the AMP scheme is stable for 
a few carefully-selected model problems. 
Two extensions of the analysis in~\cite{fib2014} are considered, 
including a first-order accurate discretization of a viscous model problem 
and a second-order accurate discretization of an inviscid model problem.
The AMP algorithm is shown to be stable for any ratio of solid and fluid
densities, including when added-mass effects are large.  
On the other hand, the traditional algorithm involving a Dirichlet-Neumann coupling 
is shown to be unconditionally unstable as added-mass effects become large with grid 
refinement. The algorithm is verified for accuracy and stability for 
set of new exact benchmark solutions. These new solutions are \textit{elastic piston}
problems, where finite interface deformations are permitted.
The AMP scheme is found to be stable and second-order accurate even for very difficult 
cases of very light solids.
\end{abstract}







\input texFiles/intro
\input texFiles/governing-new
\input texFiles/algorithm-new

\input texFiles/monolithic-new2
\input texFiles/viscousAnalysis-new

\input texFiles/inviscidAnalysis-new

\input texFiles/numericalResults

\input texFiles/conclusions

\appendix
\input texFiles/analysisAppendix-new



\clearpage
\bibliographystyle{elsart-num}
\bibliography{journal-ISI,jwb,henshaw,henshawPapers,fsi,fibrmp}

\end{document}

%% file: texFiles/defs.tex
\newcommand{\citeCount}[1]{}
\newcommand{\f}[2]{\frac{#1}{#2}}
\newcommand{\defeq}{\overset{{\rm def}}{=}}
\newcommand{\eqdef}{\overset{{\rm def}}{=}}

\newcommand{\dx}{\Delta x}
\newcommand{\dxs}{{\Delta \bar{x}}}
\newcommand{\dy}{\Delta y}
\newcommand{\dt}{\Delta t}

\def\ba#1\ea{\begin{align}#1\end{align}}

\def\bas#1\eas{\begin{align*}#1\end{align*}}

\def\bat#1\eat{\begin{alignat}{3}#1\end{alignat}}

\def\bats#1\eats{\begin{alignat*}{3}#1\end{alignat*}}

\newcommand{\bse}{\begin{subequations}}
\newcommand{\ese}{\end{subequations}}

\newcommand{\Lvh}{\Lv_h}
\newcommand{\Nvh}{\Nv_h}
    
\newcommand{\bogus}[1]{{}}

\newcommand{\av}{\mathbf{ a}}

\newcommand{\iv}{\mathbf{ i}}

\newcommand{\kv}{\mathbf{ k}}

\newcommand{\nv}{\mathbf{ n}}

\newcommand{\qv}{\mathbf{ q}}
\newcommand{\rv}{\mathbf{ r}}

\newcommand{\tv}{\mathbf{ t}}
\newcommand{\uv}{\mathbf{ u}}
\newcommand{\vv}{\mathbf{ v}}

\newcommand{\xv}{\mathbf{ x}}

\newcommand{\Fv}{\mathbf{ F}}

\newcommand{\Iv}{\mathbf{ I}}

\newcommand{\Lv}{\mathbf{ L}}

\newcommand{\Nv}{\mathbf{ N}}

\newcommand{\Vv}{\mathbf{ V}}

\newcommand{\tauv}{\boldsymbol{\tau}}

\newcommand{\sigmav}{\boldsymbol{\sigma}}

\newcommand{\grad}{\nabla}
\newcommand{\curl}{\times}
\newcommand{\grads}{\nabla_{\xsv}}

\clearpage

\newcommand{\tnv}{\tv}

\newcommand{\rhos}{\bar{\rho}}
\newcommand{\cp}{\bar{c}_p}
\newcommand{\cs}{\bar{c}_s}

\newcommand{\zs}{{\bar {z}}}  
\newcommand{\zp}{z_p}  
\newcommand{\us}{\bar{u}}
\newcommand{\rs}{\bar{r}}

\newcommand{\vs}{\bar{v}}
\newcommand{\sigmas}{\bar{\sigma}}

\newcommand{\xs}{\bar{x}}
\newcommand{\xsv}{\bar{\xv}}
\newcommand{\ys}{{\bar{y}}}

\newcommand{\zf}{z_f}

\newcommand{\ssy}{\bar{\sigma}_{22}}
\newcommand{\ssx}{\bar{\sigma}_{21}}
\newcommand{\vsy}{\vs_{2}}
\newcommand{\vfy}{v_{2}}

\newcommand{\ssyh}{{\bar{\sigma}}_{22}}
\newcommand{\ssxh}{{\bar{\sigma}}_{21}}
\newcommand{\vsyh}{{\bar{v}}_2}
\newcommand{\vfyh}{{v}_2}
\newcommand{\ph}{{p}}

\newcommand{\ac}{b} 
\newcommand{\bc}{a} 
\newcommand{\dc}{d}
\newcommand{\acv}{\mathbf{\ac}}
\newcommand{\am}{a_m}
\newcommand{\bm}{b_m}
\newcommand{\aone}{a_1}
\newcommand{\bone}{b_1}
\newcommand{\atwo}{a_2}
\newcommand{\btwo}{b_2}

\newcommand{\Heff}{H_{\text{eff}}}

\newcommand{\p}{\phi} 
\newcommand{\rev}{\mathbf{r}} 
\newcommand{\qn}[1]{\tilde{q}_{#1}}
\newcommand{\rn}[1]{\tilde{r}_{#1}}
\newcommand{\sn}[1]{\tilde{s}_{#1}}

\newcommand{\tac}{\tilde{\ac}}
\newcommand{\tbc}{\tilde{\bc}}
\newcommand{\tdc}{\tilde{\dc}}
\newcommand{\tacv}{\tilde{\acv}}

\newcommand{\dComp}{\mathcal{H}}
\newcommand{\dMat}{\mathcal{G}}
\newcommand{\dMatI}{\mathcal{H}}

\newcommand{\Usv}{\bar{\uv}}
\newcommand{\Vsv}{\bar{\vv}}
\newcommand{\Fsv}{\bar{\Fv}}

\newcommand{\Qsv}{\bar{\qv}}

\newcommand{\Sigmasv}{\bar{\sigmav}}

\newcommand{\Grad}{\grad_h}
\newcommand{\Grads}{\bar{\grad}_h}

\newcommand{\isv}{{\bar{\iv}}}
\newcommand{\Sigmav}{{\sigmav}}

\newcommand{\usv}{\bar{\uv}}
\newcommand{\vsv}{\bar{\vv}}

\newcommand{\xvs}{\bar{\xv}}

\newcommand{\sigmasv}{\bar{\sigmav}}

\newcommand{\amp}{A}

\newcommand{\normalss}{\sffamily}

\newcommand{\Hs}{\bar{H}}

\newcommand{\lambdas}{\bar{\lambda}}
\newcommand{\mus}{\bar{\mu}}

\newcommand{\OmegaF}{\Omega}
\newcommand{\OmegaS}{{\bar\Omega}_0}

\newcommand{\OmegaFh}{\Omega_h}
\newcommand{\OmegaSh}{{\bar\Omega}_h}

\newcommand{\GammaIs}{{\bar\Gamma}_0}
\newcommand{\GammaIf}{\Gamma}
\newcommand{\GammaIfh}{\Gamma_h}
\newcommand{\GammaIsh}{{\bar\Gamma}_h}
\newcommand{\GammaIFh}{\Gamma_h}

\renewcommand{\zp}{\bar{z}_p}
\renewcommand{\zs}{\bar{z}_s}
\newcommand{\zb}{\bar{z}}

\newcommand{\Mr}{\mathcal{M}}
\newcommand{\ACFL}{\mathcal{A}_{\text{CFL}}}

\newcommand{\cb}{\bar{c}}
\newcommand{\kx}{k_x} 

\newcommand{\ly}{\lambda_y}
\newcommand{\lx}{\lambda_x}



%% file: texFiles/trimFig.tex
\newlength{\tfwidth}
\newlength{\tfheight}
\newlength{\tfxa}
\newlength{\tfxb}
\newlength{\tfya}
\newlength{\tfyb}
\newcommand{\trimFigWithBox}[6]{%
\setlength\fboxsep{0pt}%
\setlength\fboxrule{1.0pt}
\fbox{\includegraphics[width=#2, clip, trim=#3 #4 #5 #6]{#1}}%
}
\newcommand{\trimFigNoBox}[6]{%
\setlength\fboxsep{1pt}
\setlength\fboxrule{0.0pt}
\fbox{\includegraphics[width=#2, clip, trim=#3 #4 #5 #6]{#1}}%
}
\newcommand{\trimFigHeightWithBox}[6]{%
\setlength\fboxsep{0pt}%
\setlength\fboxrule{1.0pt}
\fbox{\includegraphics[height=#2, clip, trim=#3 #4 #5 #6]{#1}}%
}
\newcommand{\trimFigHeightNoBox}[6]{%
\setlength\fboxsep{1pt}
\setlength\fboxrule{0.0pt}
\fbox{\includegraphics[height=#2, clip, trim=#3 #4 #5 #6]{#1}}%
}

\newsavebox\figBox

\newcommand{\trimw}[6]{%
\sbox\figBox{\includegraphics{#1}}
\setlength{\tfwidth}{\the\wd\figBox}
\setlength{\tfheight}{\the\ht\figBox}
\setlength{\tfxa}{\tfwidth*\real{#3}}%
\setlength{\tfxb}{\tfwidth*\real{#4}}%
\setlength{\tfya}{\tfheight*\real{#5}}%
\setlength{\tfyb}{\tfheight*\real{#6}}%
\trimFigNoBox{#1}{#2}{\tfxa}{\tfya}{\tfxb}{\tfyb}%
}

\newcommand{\trimwb}[6]{%

\sbox\figBox{\includegraphics{#1}}
\setlength{\tfwidth}{\the\wd\figBox}
\setlength{\tfheight}{\the\ht\figBox}
\setlength{\tfxa}{\tfwidth*\real{#3}}%
\setlength{\tfxb}{\tfwidth*\real{#4}}%
\setlength{\tfya}{\tfheight*\real{#5}}%
\setlength{\tfyb}{\tfheight*\real{#6}}%
\trimFigWithBox{#1}{#2}{\tfxa}{\tfya}{\tfxb}{\tfyb}%
}

\newcommand{\trimh}[6]{%
\sbox\figBox{\includegraphics{#1}}
\setlength{\tfwidth}{\the\wd\figBox}
\setlength{\tfheight}{\the\ht\figBox}
\setlength{\tfxa}{\tfwidth*\real{#3}}%
\setlength{\tfxb}{\tfwidth*\real{#4}}%
\setlength{\tfya}{\tfheight*\real{#5}}%
\setlength{\tfyb}{\tfheight*\real{#6}}%
\trimFigHeightNoBox{#1}{#2}{\tfxa}{\tfya}{\tfxb}{\tfyb}%
}

\newcommand{\trimhb}[6]{%

\sbox\figBox{\includegraphics{#1}}
\setlength{\tfwidth}{\the\wd\figBox}
\setlength{\tfheight}{\the\ht\figBox}
\setlength{\tfxa}{\tfwidth*\real{#3}}%
\setlength{\tfxb}{\tfwidth*\real{#4}}%
\setlength{\tfya}{\tfheight*\real{#5}}%
\setlength{\tfyb}{\tfheight*\real{#6}}%
\trimFigHeightWithBox{#1}{#2}{\tfxa}{\tfya}{\tfxb}{\tfyb}%
}

%% file: texFiles/intro.tex


\section{Introduction} \label{sec:intro}

In this paper, we consider fluid-structure interaction (FSI) of incompressible fluids 
and bulk elastic solids. FSI arises in many scientific and engineering applications
including flow-induced vibrations of structures (i.e.,~aircraft wings, undersea cables, 
wind turbines, and bridges) and simulating blood flow in arteries and veins.
FSI is typically modeled by solving suitable partial differential equations in the
fluid and solid domains together with coupling conditions involving
velocity and stress at the fluid-solid interface. 
FSI algorithms can either be categorized as monolithic schemes, where the solutions 
for the fluid and solid are advanced implicitly as one large system, or 
partitioned schemes, where the evolution of the fluid and solid solutions 
are decoupled from the interface conditions.
While partitioned schemes are generally more computationally efficient than 
monolithic schemes, 
they may suffer from instabilities when added-mass effects are present.
The goal of current work is to develop a robust partitioned scheme.

In the \textit{traditional partitioned} (TP) scheme, the solid provides a 
Dirichlet (no-slip) boundary condition for the fluid, and then the fluid supplies a 
Neumann (traction) boundary condition for the solid.
Without sub-iterations, it is well known that the TP scheme is unstable for light solids. 
On the other hand, we can consider the \textit{anti-traditional partitioned} (ATP) 
scheme, which reverses the role of the solid and fluid. In this scheme, the solid
provides a Neumann (traction) boundary condition for the fluid and the fluid 
supplies a Dirichlet (no-slip) boundary condition for the solid.
This scheme is only stable when the solid is sufficiently light, see~\cite{fib2014}.
Stability of the aforementioned schemes is typically addressed by using an iteration 
of the interface conditions.  
Research on acceleration techniques such as Aitken acceleration or quasi-Newton methods
has proven that the number of sub-iterations can be greatly 
reduced, see~\cite{Küttler2008,MEHL2016869}.
Reduction of sub-iterations can also be achieved by considering 
Robin-Neumann or Robin-Robin coupling instead of the traditional Dirichlet-Neumann coupling 
\cite{Wang2018,BASTING2017312,BadiaNobileVergara2008,MokWallRamm2001,FernandezMullaertVidrascu2014,FernandezMullaertVidrascu2013,FernandezLandajuela2014,Gerardo_GiordaNobileVergara2010,BadiaNobileVergara2009,NobilePozzoliVergara2014}.


In recent work~\cite{fib2014}, we developed a new class of Added-Mass Partitioned (AMP)
algorithms for FSI problems coupling incompressible flow and elastic solids.
The fluid AMP interface conditions are Robin conditions motivated from the propagation of 
characteristics out of the solid domain.
One important property of the AMP scheme is that it approaches the TP scheme in the 
heavy solid limit and the ATP scheme in the light solid limit.
The algorithms use a fractional-step approach for the fluid in 
which the velocity is advanced in one stage followed by the solution of a 
Poisson problem for the pressure. 
During the pressure step, a Robin condition which involves both the interface traction and 
acceleration is applied.
This condition is derived by manipulating the characteristic conditions using 
the governing equations of the fluid.
The central aim of this paper is to extend the AMP algorithm in~\cite{fib2014} for
IMEX fluid solvers, where the viscous terms are treated implicitly and the convection 
and pressure gradient terms are treated explicitly.
The adaptation to IMEX solvers necessitated important modifications for stability. 
An impedance-based projection of the interface velocity and stress
(motivated from~\cite{fsi2012} for the case of compressible fluids and linearly-elastic solids)
is applied after the fluid and solid updates to ensure the heavy and 
light solid limits are achieved. 
An implicit version of this projection is also used during the stage for the fluid velocity.
These projections rely on defining the fluid impedance, which is a 
key ingredient which leads to a stable scheme.


A fluid impedance was defined for incompressible flows and bulk solids in~\cite{fib2014}.
This definition leads to a stable scheme for explicit time step restrictions
which are typically limited by viscous terms, 
but led to instabilities for IMEX schemes where the time step is less restrictive 
and chosen by the convective terms.
In this paper, a suitable definition for the fluid impedance is derived through an analysis 
of a discretization of the fluid equations near the interface.
In this analysis, it was revealed that the impedance has an inertial component for
treating added-mass instabilities and a viscous component for treating 
added-damping instabilities.
For this choice of fluid impedance, an analysis of the AMP scheme is performed
to show stability for any ratio of the mass of the solid to that of the fluid.
The stability analysis extends~\cite{fib2014} in two different directions.
Previously in~\cite{fib2014}, stability was analyzed for a \textit{first-order accurate}
scheme applied to an \textit{inviscid} incompressible fluid and acoustic solid 
in a Cartesian geometry.
First we consider an extension of the model problem to include \textit{viscous} fluids.
Due to the introduction of viscosity, there are both pressure 
and shear forces at the interface and 
instabilities can arise from both added-mass and added-damping effects.
This new analysis is complicated by the introduction of more dimensionless parameters, namely the 
viscous CFL number in the fluid.
Despite the added complication, the AMP scheme is shown to be stable for all possible mass ratios
and viscous CFL numbers when applied to this new model problem.
Secondly, we consider the same inviscid model problem in~\cite{fib2014} and analyze a 
\textit{second-order accurate} scheme. 
In both extensions, we consider a detailed study of the traditional and anti-traditional schemes. 
The stability and accuracy of the AMP algorithm is tested on new exact solutions 
for a wide range of fluid and solid densities.
The stability and accuracy of the AMP algorithm is verified numerically for new exact solutions.
In our companion paper~\cite{fibr2019}, the AMP algorithm is implemented using 
deforming composite grids for curvilinear geometries.

The remaining sections of the paper are organized as follows.  The equations governing the FSI problem are described in Section~\ref{sec:governingEquations}.  
The AMP algorithm is summarized in Section~\ref{sec:algorithm}.  
In Section~\ref{sec:monolithic}, the fluid impedance is derived following an analysis of 
an FSI discretization near the interface. 
The stability of the AMP algorithm is analyzed in Section~\ref{sec:viscousAnalysis}
for a viscous model problem and in Section~\ref{sec:inviscidAnalysis} for an inviscid model problem.
Section~\ref{sec:numericalResults} provides numerical results confirming the stability and accuracy of the scheme.  Some of the results use new exact solutions of benchmark FSI problems.  
Conclusions are given in Section~\ref{sec:conclusions}.

%% file: texFiles/governing-new.tex
\section{Governing equations} \label{sec:governingEquations}

We consider the coupled evolution of an incompressible fluid and a linear elastic solid. 
The fluid occupies the domain $\xv \in \OmegaF(t)$, where
$\xv=(x_1,x_2,x_3)$ is a vector of physical coordinates and $t$ is time. 
The velocity-pressure form of the incompressible Navier-Stokes equations is given by
\bse
\begin{alignat}{2}
\rho \vv_t + \rho (\vv \cdot \grad) \vv + \grad p &=
\mu \Delta \vv, \qquad &&\xv \in \OmegaF(t),
\label{eq:fluidMomentum} \\
\Delta p &= -\rho \grad \vv : \left(\grad \vv\right)^T, \qquad &&\xv \in \OmegaF(t),
\label{eq:pressurePoissonEquation}
\end{alignat}
\ese
where
\begin{align}
\grad \vv : \left(\grad \vv\right)^T \equiv 
\sum_{i=1}^3 \sum_{j=1}^3 \frac{\partial v_i}{\partial x_j} \frac{\partial v_j}{\partial x_i}.
\end{align}
Here, $\vv(\xv,t)$ is the velocity, $p(\xv,t)$ is the pressure,
$\rho$ is the (constant) density, and $\mu$ is the (constant) dynamic viscosity. 
The fluid stress tensor is given by
\begin{align}
\sigmav = - p \Iv + \tauv, 
\qquad 
\tauv &= \mu \left(\grad \vv + \left(\grad \vv\right)^T \right),
\label{eq:fluidStressTensor} 
\end{align}
where $\Iv$ is the identity matrix and $\tauv(\xv,t)$ is
the viscous stress tensor.  
In the velocity-pressure form of the equations, an extra boundary condition is required and a suitable choice is to impose $\grad \cdot \vv = 0$ for $\xv\in\partial\OmegaF(t)$, see~\cite{splitStep2003}.

The equations for the solid are written in terms of the Lagrangian coordinate 
$\xsv=(\bar x_1,\bar x_2,\bar x_3)$ for
a reference configuration $\xsv \in \OmegaS$ at $t=0$.  
(An overbar is used here and elsewhere to denote
quantities belonging to the solid.)
The position of the solid in physical space is determined by the mapping
\begin{align}
\xv = \xsv + \usv(\xsv,t),
\label{eq:solidPhysicalCoordinates}
\end{align}
where $\usv(\xsv,t)$ is the displacement of the solid.
The Cauchy stress tensor $\sigmasv(\xsv,t)$ for a linearly-elastic solid is defined by
\begin{align}
\sigmasv &= \lambdas (\grads \cdot \usv) \Iv 
+ \mus \left(\grads \usv + \left(\grads \usv\right)^T\right),
\label{eq:solidStressTensor}
\end{align}
where $\lambdas$ and $\mus$ are Lam\'e parameters (taken to be constants).
The solid equations are considered as a 
first-order system of PDEs in time and space, following~\cite{smog2012}, and are given by
\bse
\label{eq:solidEquations}
\begin{alignat}{2}
\usv_t &= \vsv, \qquad &&\xvs \in \OmegaS ,
\label{eq:solidDisplacement} \\
\rhos \vsv_t &= \grads \cdot \sigmasv , \qquad &&\xvs \in \OmegaS,
\label{eq:solidMomentum} \\
\sigmasv_t &= \lambdas (\grads \cdot \vsv) \Iv 
+ \mus \left(\grads \vsv + \left(\grads \vsv \right)^T\right),\qquad &&\xvs \in \OmegaS,
\label{eq:solidStressTensorDeriv}
\end{alignat}
\ese
where $\vsv(\xsv,t)$ is the velocity of the solid, and $\rhos$ is its density (assumed constant).  
In this form, upwind solvers can be used to advance displacement, velocity and stress of the 
solid. We note that~\eqref{eq:solidStressTensor} is enforced at $t=0$.

The fluid and solid are coupled at an interface described by
$\xv\in\GammaIf(t)$ in physical space and $\xsv \in \GammaIs$
in the corresponding reference space. At the interface, the conversion between reference 
and physical coordinates is determined by the mapping in~\eqref{eq:solidPhysicalCoordinates}. 
The interface is assumed to be smooth so that a well-defined normal to the interface exists.  
Along the interface, the following matching conditions hold:
\bse
\label{eq:matching}
\begin{alignat}{2}
\vv &= \vsv, \qquad &&\xv \in \GammaIf(t),
\label{eq:velocityMatching} \\
\sigmav \nv &= 
\sigmasv \nv, \qquad &&\xv \in \GammaIf(t),
\label{eq:tractionMatching}
\end{alignat}
\ese
%
where $\nv(\xv,t)$ is the outward unit normal to the fluid domain, i.e.~$\nv$ points from the fluid domain to the solid domain.  Suitable boundary
conditions are applied on the boundaries of the fluid and solid domians not included in $\GammaIf(t)$, and initial
conditions on $\vv$, $\usv$ and $\vsv$ are set to close the problem.

%% file: texFiles/algorithm-new.tex
\section{AMP interface conditions and algorithm} 
\label{sec:interface}
\label{sec:algorithm}

In this section, we derive the AMP interface conditions at a continuous level and discuss their implementation in the AMP algorithm.  The derivation follows the work in~\cite{fib2014}, but there are important modifications discussed to accommodate the IMEX fractional-step scheme used in the AMP algorithm to solve the equations in the fluid domain.  These modifications are guided by a consideration of the behavior of the AMP interface conditions in the limits of very light and very heavy solids.

\subsection{AMP interface conditions}

The starting point for the derivation is the matching conditions involving velocity and stress in~\eqref{eq:matching}.  Following~\cite{fib2014}, a linear combination of these conditions are expressed in terms of the outgoing characteristic variables of the solid, i.e.
\bse
\label{eq:fluidInterfaceConditions}
\begin{align}
-p + \nv^T \tauv \nv + \zp \nv^T \vv = \nv^T \sigmasv \nv + \zp \nv^T \vsv, \quad \phantom{m = 1,2,}\qquad \xv \in \GammaIf(t),
\label{eq:AMPNormalCharacteristic} \\
\tnv_m^T \tauv \nv + \zs \tnv_m^T \vv = \tnv_m^T \sigmasv \nv + \zs \tnv_m^T \vsv,  \quad m = 1,2, \qquad \xv \in \GammaIf(t),\label{eq:AMPTangentialCharacteristic}
\end{align}
\ese
where $\nv$ is the unit normal, $\tnv_m$, $m=1,2$, are mutually orthogonal unit vectors tangent to the interface, and $\zp$ and $\zs$ are impedances involving the characteristic velocities of the solid given by
\[
\zp=\rhos\cp,\qquad \zs=\rhos\cs,\qquad \cp = \sqrt{\frac{\lambdas+2 \mus}{\rhos}},\qquad \cs = \sqrt{\frac{ \mus}{\rhos}}.
\]
In the AMP algorithm, the conditions in~\eqref{eq:fluidInterfaceConditions} are interpreted as providing boundary conditions for the fluid at the interface with the outgoing characteristic quantities of the solid on the left-hand side considered to be known from a previous stage of the algorithm.  While these conditions, along with $\nabla\cdot\vv=0$ for $\xv \in \GammaIf(t)$, are sufficient conditions for the fluid equations in velocity-pressure form, a further manipulation is required to obtain suitable conditions to be used for the fractional-step solver.  The objective is to separate the conditions in~\eqref{eq:fluidInterfaceConditions} to obtain a condition to be used in the IMEX time-stepping scheme for the fluid velocity and a condition for the subsequent update for the pressure.

For the Poisson problem for the fluid pressure, the interface condition in~\eqref{eq:AMPNormalCharacteristic} is used with the momentum equation in~\eqref{eq:fluidMomentum} to derive a Robin condition for the pressure.  The momentum equation involves the acceleration of the fluid, and this quantity may be obtained on the moving fluid-solid interface using the Taylor approximation
\newcommand{\fluidPoint}{{\mathcal P}}
\newcommand{\solidPoint}{{\bar{\mathcal P}}_0}
\begin{align}
\Bigl.\vv(\xv,t-\dt)\Bigr\vert_{\xv=\fluidPoint(t-\dt)}\approx\Bigl.\bigl(\vv(\xv,t)-\dt D_t\vv(\xv,t)\bigr)\Bigr\vert_{\xv=\fluidPoint(t)},
\label{eq:fluidTaylorApprox}
\end{align}
where $D_t=\partial_t+\vv\cdot\grad$ is the material derivative, $\fluidPoint(t)$ is a point on the moving interface and $\dt$ is a time-step.  The corresponding approximation for the solid is
\begin{align}
\Bigl.\vsv(\xsv,t-\dt)\Bigr\vert_{\xsv=\solidPoint}\approx\Bigl.\bigl(\vsv(\xsv,t)-\dt\vsv_t(\xsv,t)\bigr)\Bigr\vert_{\xsv=\solidPoint},
\label{eq:solidTaylorApprox}
\end{align}
where $\solidPoint$ is the Lagrangian position associated with $\fluidPoint(t)$.  Using~\eqref{eq:fluidTaylorApprox} and~\eqref{eq:solidTaylorApprox} in~\eqref{eq:AMPNormalCharacteristic}, and assuming the fluid and solid velocities match on the interface at times $t-\dt$ and $t$, we obtain
\begin{align}
-p + \nv^T \tauv \nv + \zp \dt \nv^T D_t\vv = \nv^T \sigmasv \nv + \zp \dt \nv^T \vsv_t, \qquad \xv \in \GammaIf(t).
\label{eq:pressureBCtemp}
\end{align}
We may now eliminate the fluid acceleration using~\eqref{eq:fluidMomentum} to obtain the following Robin condition for the fluid pressure:
\begin{align}
-p - {\zp \dt\over\rho} \partial_n p =  \nv^T(\sigmasv \nv-\tauv \nv) + \zp \dt\nv^T\bigl( \vsv_t + \nu\grad\times\grad\times\vv \bigr), \qquad \xv \in \GammaIf(t),
\label{eq:AMPpressureBC}
\end{align}
where $\partial_n=\nv\cdot\grad$ is the normal derivative and $\nu=\mu/\rho$ is the kinematic viscosity of the fluid.  Following~\cite{splitStep2003}, we have used the identity, $\Delta \vv = - \grad \curl \grad \curl \vv$, noting that $\grad \cdot \vv = 0$, to replace $\Delta \vv$ on the right-hand side of~\eqref{eq:AMPpressureBC} in favor of the curl-curl operator.  This is done for improved stability of the fractional-step scheme.  The condition in~\eqref{eq:AMPpressureBC}, along with suitable conditions for $\xv\in\partial\OmegaF(t)\backslash\GammaIf(t)$, is used for the Poisson equation in~\eqref{eq:pressurePoissonEquation} for the pressure.

\newcommand{\Cam}{{\mathcal C}_{\rm AM}}
\newcommand{\Cad}{{\mathcal C}_{\rm AD}}

As was noted in~\cite{fib2014}, the remaining interface conditions
in~\eqref{eq:AMPTangentialCharacteristic}, together with the continuity equation, can be used as
boundary conditions to advance the fluid velocity.  This was found to be an effective approach for
an explicit integration of the momentum equations.
To ensure that the fluid velocity and tractions match at the end of the time step, an interface projection
is performed to give a common interface velocity $\vv^I$ and interface traction $\sigmav^I\nv$.
In analogy to the interface projection used for compressible fluids in~\cite{fsi2012,flunsi2016},
which is based on a characteristic analysis, the projection for incompressible fluids is also
proposed to be of the form of an impedance-weighted average.  For the velocity, the projection has the form 
%
\bse
\label{eq:AMPinterfaceVelocity}
\begin{alignat}{2}
\nv^T\vv^I=\;&{\zf\over\zf+\zp}\nv^T\vv+{\zp\over\zf+\zp}\nv^T\vsv+{1\over\zf+\zp}\nv^T\bigl(\sigmasv\nv-\sigmav\nv\bigr), \quad &&\phantom{m = 1,2,}
\label{eq:AMPinterfaceNormalVelocity} \\
\tnv_m^T\vv^I=\;&{\zf\over\zf+\zs}\tnv_m^T\vv+{\zs\over\zf+\zs}\tnv_m^T\vsv+{1\over\zf+\zs}\tnv_m^T\bigl(\sigmasv\nv-\sigmav\nv\bigr) ,  \quad &&m = 1,2,
\label{eq:AMPinterfaceTangentialVelocity}
\end{alignat}
\ese
while for the traction, the projection is based on an inverse impedance-weighted average of the form
\bse
\label{eq:AMPinterfaceTraction}
\begin{alignat}{2}
\nv^T\sigmav^I\nv=\;&{\zf^{-1}\over\zf^{-1}+\zp^{-1}}\nv^T\sigmav\nv+{\zp^{-1}\over\zf^{-1}+\zp^{-1}}\nv^T\sigmasv\nv
+{1\over\zf^{-1}+\zp^{-1}}\nv^T\bigl(\vsv-\vv\bigr), \quad &&\phantom{m = 1,2,}
\label{eq:AMPinterfaceNormalTraction}\\
\tnv_m^T\sigmav^I \nv=\;&{\zf^{-1}\over\zf^{-1}+\zs^{-1}}\tnv_m^T\sigmav\nv+{\zs^{-1}\over\zf^{-1}+\zs^{-1}}\tnv_m^T\sigmasv\nv
+{1\over\zf^{-1}+\zs^{-1}}\tnv_m^T\bigl(\vsv-\vv\bigr) ,  \quad &&m = 1,2.
\label{eq:AMPinterfaceTangentialTraction}
\end{alignat}
\ese
These projections introduce a {\em fluid impedance}, $\zf$, which is well defined for compressible fluids, but
has no obvious definition for incompressible fluids. However, an analysis of a discrete approximation
to the governing equations given in Section~\ref{sec:monolithic} suggests a form for $\zf$ given by
\begin{equation}
\zf \eqdef \Cam \Bigl(\frac{\rho h}{\dt}\Bigr) + \Cad \Bigl(\frac{\mu}{h}\Bigr),
\label{eq:fluidImpedanceAlg}
\end{equation}
where $h$ is an  appropriate mesh spacing and $(\Cam,\Cad)$ are constants whose approximate values are provided by the analysis.
The projections in~\eqref{eq:AMPinterfaceVelocity} and~\eqref{eq:AMPinterfaceTraction} can be used to set values of the fluid and solid velocity and traction at the interface.

For the IMEX scheme considered here, a modification of the previous approach in \cite{fib2014} is required in the implementation of the interface conditions for the fluid velocity.  The issue is informed by considering the limits of very light and heavy solids.  In the limit of a very light solid ($\zp,\zs\rightarrow0$), for example, the Robin condition in~\eqref{eq:AMPpressureBC} becomes a Dirichlet condition for the pressure, while the interface conditions in~\eqref{eq:AMPTangentialCharacteristic} reduce to matching conditions involving the shear stress of the fluid.  The latter conditions, along with the continuity constraint, provide Neumann conditions on the fluid velocity.  These conditions for the fluid pressure and velocity correspond to those for a free surface, and the latter are suitable for the implicit solution of the fluid velocity in the IMEX fractional-step scheme.

The difficulty is revealed in the limit of a very heavy solid ($\zp,\zs\rightarrow\infty$).  In this limit, the Robin condition in~\eqref{eq:AMPpressureBC} becomes a Neumann condition for the fluid pressure balancing the acceleration of the interface as determined by the solid.  This condition is analogous to the usual Neumann boundary condition for the pressure at a rigid boundary obtained from the fluid momentum equations as a compatibility condition (see~\cite{splitStep2003} for example).  The interface conditions in~\eqref{eq:AMPTangentialCharacteristic} reduce to matching conditions involving the tangential components of velocity.  However, the matching condition on the normal component of velocity,
\begin{equation}
\nv^T \vv =  \nv^T \vsv,\qquad \xv \in \GammaIf(t),
\label{eq:missingCondition}
\end{equation}
implied by~\eqref{eq:AMPNormalCharacteristic} in the limit of a heavy solid has been lost in the derivation of~\eqref{eq:AMPpressureBC}.
A remedy can be obtained by using the interface projection for the normal component of the velocity in~\eqref{eq:AMPinterfaceNormalVelocity} as a boundary condition for the implicit solution of the fluid velocity in the IMEX fractional-step scheme. 
The implementation of this approach is described next in the discussion of the AMP algorithm.

\subsection{AMP algorithm}

Algorithm~\ref{alg:ampAlgorithm} provides a consise description of the AMP time-stepping scheme (see~\cite{fibr2019} for
additional details of the implementation of the algorithm).
The algorithm advances the solution from a time $t^n$ to $t^{n+1}=t^n+\dt$.  It is assumed that the
fluid domain is represented by a grid consisting of interior points $\iv \in \OmegaFh$, boundary
points $\iv \in \partial \OmegaFh$ and interface points $\iv \in \GammaIfh$, where
$\iv=(i_1,i_2,i_3)$ is a multi-index.  Similarly, the solid reference domain is covered by a grid
with interior points $\isv \in \OmegaSh$, boundary points $\isv \in \partial \OmegaSh$ and interface
points $\isv \in \GammaIsh$.  Discrete operators, such as $\grad_{\! h}$ and $\Delta_h$, denote
approximations of the corresponding differential operators on the grid.

\input texFiles/ampAlg

The time-stepping scheme uses a predictor-corrector approach.  Steps 1--5 of Algorithm~\ref{alg:ampAlgorithm} describe the preditor steps.  Predicted values for the solid displacement $\Usv_\isv$ are obtained in Step~1 using a Lax-Wendroff-type scheme, while the solid velocity and stress $\Qsv_\isv=(\Vsv_\isv,\Sigmasv_{\isv})$ are advanced using a Godunov type scheme with numerical fluxes $\Fsv_{m, \; \isv}^\pm$ corresponding to the $\xs_m$ coordinate direction.  In Step~2, the solid displacement is used to compute the deformed fluid grid at time $t^{n+1}$.

The fluid velocity is predicted in Step~3.  Here, $\Nvh$ and $\Lvh$ represent grid operators associated with the explicit and implicit terms in the velocity update, respectively, given by
\ba
\Nvh(\vv_\iv,p_\iv) = - \bigl( (\vv_\iv - \dot{\xv}_\iv) \cdot \Grad \bigr) \vv_\iv 
- \frac{1}{\rho} \Grad \,p_\iv ,\qquad
\Lvh(\vv_\iv) = \nu \Delta_h \vv_\iv,
\ea
where $\dot{\xv}_\iv$ is the velocity of the grid.  The explicit terms are advanced using an Adams-Bashforth scheme, while the implicit terms use Crank-Nicholson.
The boundary conditions on the interface makes use of a predicted velocity, coming from the interior equation
applied on the boundary, and defined by
\ba
\Vv\sp{p}_h(\vv_\iv\sp{(p)}) \eqdef \vv_\iv^n + \frac{\dt}{2}\left(3\Nvh(\vv_\iv^ n,p_\iv^ n) - \Nvh(\vv_\iv^{n-1},p_\iv^{n-1})\right)
+ \frac{\dt}{2}\Big(\Lvh(\vv_\iv^{(p)}) + \Lvh(\vv_\iv^ n)\Big).
\label{eq:veDef}
\ea
In particular, this velocity is used in the impedance-weighted average condition
\ba
\nv^T\vv_\iv^{(p)}=\f{\zf}{\zf+\zp}\nv^T\Vv^p_h(\vv_\iv\sp{(p)})+ \f{\zp}{\zf+\zp}\nv^T\vsv_\isv^{(p)},\qquad \iv \in \GammaIfh,\quad \isv \in \GammaIsh,
\label{eq:normalVelociityFix}
\ea
which is obtained from the projection in~\eqref{eq:AMPinterfaceNormalVelocity}. Here the 
  the term involving the jump in the stress is dropped (because a suitable approximation for the fluid stress is unavailable at the predictor stage). Notice that~\eqref{eq:normalVelociityFix} is an implicit condition on $\vv_\iv^{(p)}$ which
  appears on the left and right-hand sides.
In the light-solid limit ($\zp\rightarrow0$), the boundary condition in~\eqref{eq:normalVelociityFix} reduces to $\nv^T\vv_\iv^{(p)}=\nv^T\Vv^p_h(\vv_\iv\sp{(p)})$, which simply sets the normal component of the fluid velocity to be equal to that given by the interior time-stepping scheme applied on the boundary.  In the heavy-solid limit ($\zp\rightarrow\infty$), \eqref{eq:normalVelociityFix} becomes
\[
\nv^T\vv_\iv^{(p)}=\nv^T\vsv_\isv^{(p)},\qquad \iv \in \GammaIfh,\quad \isv \in \GammaIsh,
\]
which recovers the matching condition in~\eqref{eq:missingCondition}.  Our later analysis of a viscous model problem (Section~\ref{sec:viscousAnalysis}) and an inviscid model problem (Section~\ref{sec:inviscidAnalysis}), and subsequent numerical results (Section~\ref{sec:numericalResults}), verify that the boundary conditions used to advance the fluid velocity in the fractional-step scheme lead to stable and accurate results for a wide range of solid densities.

Steps~4 and~5 complete the set of steps belonging to the predictor stage of the algorithm.  The predicted fluid pressure is computed in Step~4 by solving a discrete Poisson problem.  This elliptic problem uses a discrete approximation of the Robin condition in~\eqref{eq:AMPpressureBC}.  Finally, interface values for the solid velocity and traction are obtained in Step~5 using the impedance-weighted projections in~\eqref{eq:AMPinterfaceVelocity} and~\eqref{eq:AMPinterfaceTraction}.  These interface values overwrite the corresponding predicted values of the solid on the boundary.

The set of corrector steps consisting of Steps~6--9 essentially mirror those of the predictor.  The fluid grid is recomputed in Step~6 using an updated grid velocity obtained from solid velocity.  In Step~7, the discrete fluid velocity at $t^{n+1}$ is determined, now using an Adams-Moulton corrector.  The fluid pressure is updated in Step~8 according to the solution of a discrete Poisson problem.  Lastly, the solid velocity and traction are set equal to the corrected fluid values in Step~9, so that the fluid and solid interface values agree at the end of the corrector stage of the AMP time-stepping scheme.

%% file: texFiles/ampAlg.tex
{
    \def\alignspace{\hspace{1.2em}}
\begin{algorithm}\caption{\rm Added-mass partitioned (AMP) scheme \label{alg:ampAlgorithm}}
\small
\[
\begin{array}{l}
\hbox{// \textsl{Predictor steps}}\smallskip\\
1.\text{ Predict solid:}\\ 
\alignspace 
\begin{cases}
\Usv_\isv^{(p)} = \Usv_\isv^n
+ \dt \Vsv_\isv^n
+ \frac{\dt^2}{2\rhos} \Grads \cdot \Sigmasv_\isv^n, 
& \qquad \isv \in \OmegaSh. \\
\Qsv_\isv^{(p)} = \Qsv_\isv^n 
- \dt \sum_{m=1}^3 \frac{1}{\dxs_m} 
\Bigl( \Fsv_{m, \; \isv}^+ - \Fsv_{m, \; \isv}^- \Bigr), & \qquad \isv \in \OmegaSh,
\end{cases} 
\medskip\\
2.\text{ Predict fluid grid: advance fluid grid to $t^{n+1}$ using $\Usv_\isv^{(p)}$ for $\isv \in \GammaIsh$, and compute grid velocity.}
\medskip\\
3.\text{ Predict fluid velocity:}\\ 
\alignspace 
\begin{cases}
\vv_\iv^{(p)} = \vv_\iv^n + \frac{\dt}{2}\bigl(3\Nvh(\vv_\iv\sp n,p_\iv\sp n) - \Nvh(\vv_\iv\sp{n-1},p_\iv\sp{n-1})\bigr)
+ \frac{\dt}{2}\bigl(\Lvh(\vv_\iv\sp{(p)}) + \Lvh(\vv_\iv\sp n)\bigr),
   & \qquad \iv \in \OmegaFh\backslash\GammaIfh, \\
 \tnv_m^T \tauv_{\iv}^{(p)} \nv + \zs \tnv_m^T \vv_\iv^{(p)} =
 \tnv_m^T \Sigmasv_\isv^{(p)} \nv
+ \zs \tnv_m^T \Vsv_{\isv}^{(p)},
 & \qquad\iv \in \GammaIfh, \; \isv \in \GammaIsh,  \\
\Grad\cdot\vv_\iv\sp{(p)}=0,
 & \qquad \iv \in \GammaIfh, \\
 \nv^T\vv_\iv^{(p)}=\f{\zf}{\zf+\zp}\nv^T\Vv^p_h(\vv_\iv\sp{(p)})+ \f{\zp}{\zf+\zp}\nv^T\vsv_\isv^{(p)}, \qquad \tnv_m^T\vv_\iv^{(p)}=\tnv_m^T\Vv\sp{p}_h(\vv_\iv\sp{(p)}), & \qquad \iv \in \GammaIfh, \; \isv \in \GammaIsh, \\
%
\text{Velocity boundary conditions on $\partial \OmegaFh \backslash \GammaIfh$.}
\end{cases}
\medskip\\
4.\text{ Predict fluid pressure:}\\ 
\alignspace 
\begin{cases}
\Delta_h p_\iv^{(p)} = -\rho \Grad \vv_\iv^{(p)} : \bigl(\Grad \vv_\iv^{(p)}  \bigr)^T + \alpha_\iv\Grad \cdot \vv_\iv^{(p)},
&\quad \iv \in \OmegaFh, \\
-p_\iv^{(p)}
-\frac{\zp \dt}{\rho} (\nv\cdot\Grad) p_\iv^{(p)}
= 
\nv^T \bigl(\Sigmasv_\isv^{(p)} \nv - \tauv_\iv^{(p)} \nv\bigr) 
+ \zp \dt
 \nv^T \bigl((\vsv_t)_\isv^{(p)} + \nu \Grad \curl \Grad \curl \vv_\iv^{(p)} \bigr), &\quad \iv \in \GammaIFh, \; \isv \in \GammaIsh, \\
\text{Pressure boundary conditions on $\partial \OmegaFh \backslash \GammaIfh$.}
\end{cases}
\medskip\\
5.\text{ Project solid interface:}\\ 
\alignspace 
\begin{cases}
\nv^T \Vsv_\isv^I = \frac{\zf}{\zf+\zp} \nv^T \vv_\iv^{(p)}
     + \frac{\zp}{\zf+\zp} \nv^T \Vsv_\isv^{(p)} 
     + { \f{1}{\zf+\zp}\bigl( \nv^T\Sigmasv^{(p)}_\iv\nv - \nv^T\Sigmav^{(p)}_\iv\nv \bigr)},
     &\qquad \isv \in \GammaIsh, \iv \in \GammaIfh,\\
\tv_m^T \Vsv_\isv^I = \frac{\zf}{\zf+\zs} \tv_m^T \vv_\iv^{(p)}
     + \frac{\zs}{\zf+\zs} \tv_m^T \Vsv_\isv^{(p)}
     + { \f{1}{\zf+\zs}\bigl( \tv_m^T\Sigmasv^{(p)}_\iv\nv - \tv_m^T\Sigmav^{(p)}_\iv\nv \bigr)},
     &\qquad \isv \in \GammaIsh, \iv \in \GammaIfh,\\
\nv^T \Sigmasv^I_{\isv} \nv = \frac{\zf^{-1}}{\zf^{-1}+\zp^{-1}}\nv^T\Sigmav^{(p)}_\iv \nv
+ \frac{\zp^{-1}}{\zf^{-1}+\zp^{-1}} \nv^T\Sigmasv_\isv^{(p)} \nv
+ { \f{1}{\zf^{-1}+\zp^{-1}}\bigl( \nv^T\vsv^{(p)}_\isv - \nv^T\vv^{(p)}_\iv \bigr)},
     &\qquad \isv \in \GammaIsh, \iv \in \GammaIfh. \\
\tv_m^T \Sigmasv^I_{\isv} \nv = \frac{\zf^{-1}}{\zf^{-1}+\zs^{-1}}\tv_m^T\Sigmav^{(p)}_\iv \nv
    + \frac{\zs^{-1}}{\zf^{-1}+\zs^{-1}} \tv_m^T\Sigmasv_\isv^{(p)} \nv
    + { \f{1}{\zf^{-1}+\zs^{-1}}\bigl( \tv_m^T\vsv^{(p)}_\isv - \tv_m^T\vv^{(p)}_\iv \bigr)},
        &\qquad \isv \in \GammaIsh, \iv \in \GammaIfh,\\
 { \vsv_\isv^{(p)} \leftarrow \vsv_\isv^I, \quad \sigmasv_\isv^{(p)}\nv \leftarrow \sigmasv_\isv^I\nv,}  &\qquad \isv \in \GammaIsh, \iv \in \GammaIfh, \\
\text{Apply solid boundary conditions and set all ghost points.}
\end{cases}
\medskip\\
\hbox{// \textsl{Corrector steps}}\smallskip\\
6.\text{ Correct fluid grid: recompute grid velocity using $\Vsv_\isv^I$ for $\isv \in \GammaIsh$.}
\medskip\\
7.\text{ Correct fluid velocity:} \\
\alignspace 
\begin{cases}
\vv_\iv^{n+1} = \vv_\iv^n + \frac{\dt}{2}\bigl(\Nvh(\vv_\iv^{(p)},p_\iv^{(p)}) + \Nvh(\vv_\iv^{n},p_\iv^{n})\bigr)
+ \frac{\dt}{2}\bigl(\Lvh(\vv_\iv^{n+1}) + \Lvh(\vv_\iv^ n)\bigr),
&\qquad    \iv \in \OmegaFh\backslash\GammaIfh, \\
 \tnv_m^T \tauv_{\iv}^{n+1} \nv + \zs \tnv_m^T \vv_\iv^{n+1} =
 \tnv_m^T \Sigmasv_\isv^{I} \nv
+ \zs \tnv_m^T \Vsv_{\isv}^{I},
 & \qquad\iv \in \GammaIfh, \; \isv \in \GammaIsh,  \\
\Grad\cdot\vv_\iv\sp{n+1}=0,
 & \qquad \iv \in \GammaIfh, \\
 \nv^T\vv_\iv^{n+1}=\f{\zf}{\zf+\zp}\nv^T\Vv_h(\vv_\iv\sp{n+1})+ \f{\zp}{\zf+\zp}\nv^T\vsv_\isv^{I}
, \qquad \tnv_m^T\vv_\iv^{n+1}=\tnv_m^T\vv_\iv^{(e)},  \qquad m = 1,2, & \qquad \iv \in \GammaIfh, \; \isv \in \GammaIsh, \\
\text{Velocity boundary conditions on $\partial \OmegaFh \backslash \GammaIfh$.}
\end{cases}
\medskip\\
8.\text{ Correct fluid pressure.} \\
\alignspace 
\begin{cases}
\Delta_h p_\iv^{n+1} = -\rho \Grad \vv_\iv^{n+1} : \bigl(\Grad \vv_\iv^{n+1}  \bigr)^T + \alpha_\iv\Grad \cdot \vv_\iv^{n+1},
&\quad \iv \in \OmegaFh, \\
-p_\iv^{n+1}-\frac{\zp \dt}{\rho} (\nv\cdot\Grad) p_\iv^{n+1}
= 
\nv^T \bigl(\Sigmasv_\isv^{I} \nv - \tauv_\iv^{n+1} \nv\bigr) 
+ \zp \dt
 \nv^T \bigl((\vsv_t)_\isv^{I} + \nu \Grad \curl \Grad \curl \vv_\iv^{n+1} \bigr), &\quad \iv \in \GammaIFh, \; \isv \in \GammaIsh, \\
\text{Pressure boundary conditions on $\partial \OmegaFh \backslash \GammaIfh$.}
\end{cases}
\medskip\\
9.\text{ Correct solid interface.} \\
\alignspace 
\begin{cases}
 { \Vsv_\isv^{n+1} =  \vv_\iv^{n+1}} , &\qquad \isv \in \GammaIsh, \iv \in \GammaIfh,\\ 
 { \Sigmasv^{n+1}_{\isv} \nv = \Sigmav^{n+1}_\iv \nv, } &\qquad \isv \in \GammaIsh, \iv \in \GammaIfh, \\
%
\text{Reset ghost points corresponding to $\isv \in \GammaIsh$.}
\end{cases}
\\
\end{array}
\]
\label{alg:amp}
\end{algorithm}
}

%% file: texFiles/monolithic-new2.tex
\section{Derivation of the fluid impedance} \label{sec:monolithic}

The focus of this section is an analysis of a FSI problem leading to an expression for the fluid
impedance that guides the choice for $\zf$ introduced in~\eqref{eq:fluidImpedanceAlg}.  A value for
the fluid impedance is required in the formulas for the interface projections
in~\eqref{eq:AMPinterfaceVelocity} and~\eqref{eq:AMPinterfaceTraction}, and also for the
implementation of the IMEX fractional-step scheme for the fluid velocity.  The interface projections
were also used in our earlier paper~\cite{fib2014}, as these provide formulas for the interface
velocity and traction that ensure the fluid and solid velocities and tractions match at the
  interface and that smoothly accommodate the limiting cases of light and heavy solids.  The
earlier paper introduced a fluid impedance given by $\zf=\rho H/\dt$, where $H$ was a measure of the depth of the fluid layer,
and it was found that this choice led to a stable AMP algorithm when using an explicit
fractional-step scheme for the fluid. It was also noted that the scheme was quite insensitive to the choice of $H$.
 For the present AMP scheme, we use an IMEX fractional-step
scheme.  Since the viscous terms in the fluid equations are now treated implicitly, the viscous CFL
number, $\Lambda=\nu \dt / h^2$, can be large in which case the choice for $\zf$ used
in~\cite{fib2014} is no longer sufficient for stability.  We have found that the difficulty can be
resolved by considering an analysis of a more general FSI model problem for which the viscous terms in the
fluid equations contribute.

\newcommand{\bsp}{b_p}
\newcommand{\bss}{b_s}

Consider a FSI model problem in which the fluid occupies the two-dimensional domain, $\OmegaF$, given by $0<x<L$, $y>0$, while the solid exists on the domain, $\OmegaS$, for $0<x<L$, $y<0$, see Figure~\ref{fig:rectangularModelProblemCartoonInfinite}.  The fluid-solid interface, $\Gamma$, of length $L$ is linearized about a flat surface, $y=0$.  The equations governing the model problem are
\bse
\begin{align}
\text{Fluid: }& \begin{cases}
\rho \partial_t\vv + \grad p = \mu \Delta \vv,\quad & \xv\in\OmegaF, \\
\Delta p = 0, & \xv \in \OmegaF, \\
\grad \cdot \vv = 0, & \xv \in \Gamma,
\end{cases}\label{eq:modelFluidEquations}\\
\text{Solid: }& \begin{cases}
\rhos \partial_t \vsv = \grad \cdot \sigmasv, & \xv\in\OmegaS, \\
\partial_t \sigmasv = \lambdas (\grad \cdot \vv) \Iv + \mus (\grad \vv + (\grad \vv)^T), \quad & \xv\in\OmegaS, \\
\end{cases}  \label{eq:modelSolidEquations}\\
\text{Interface: }& \begin{cases}
\vv = \vsv, & \xv \in \Gamma,\\
\sigmav \nv = \sigmasv \nv, \quad & \xv \in \Gamma.
\label{eq:modelInterface}
\end{cases}
\end{align}
\ese
Solutions of the model problem are assumed to be periodic in $x$ with period equal to $L$, and bounded as $y\rightarrow\pm\infty$.

\input texFiles/rectangularModelProblemCartoonInfinite

The equations governing the fluid and solid are discretized in the $x$-direction on a 
uniform grid, $x_\ell = \ell \Delta x$ for $\ell = 0, 1,\ldots N_x$, with grid spacing 
$\dx = L/N_x$.  Since the problem is periodic, each variable can be represented as a 
discrete Fourier series
\begin{align}
q(x,y,t) = \sum_{k=-N_x/2}^{N_x/2}  e^{2 \pi i k x/ L} \hat{q}_k (y,t), 
\qquad x\in[0,L],
\label{eq:discreteFourier}
\end{align}
where $\hat{q}_k(y,t)$ are the Fourier coefficient functions and $N_x$ is an integer, assumed to be even for convenience.  Taking a finite Fourier transform of the fluid equations in~\eqref{eq:modelFluidEquations} gives
\bse
\label{eq:fourierFull}
\begin{alignat}{2}
\rho \partial_t v_1 + i \kx p &= \mu (\partial_y^2 - \kx^2) v_1, \qquad && y>0,\label{eq:fourierV1}\\
\rho \partial_t v_2 + \partial_y p &= \mu (\partial_y^2 - \kx^2) v_2, \qquad && y>0, \label{eq:fourierV2}\\
(\partial_y^2 - \kx^2) p &= 0, \qquad && y>0, \label{eq:fourierP}
\end{alignat}
\ese
where $\kx=2\pi k/L$.  The hats and $k$ subscripts on the coefficient functions in~\eqref{eq:fourierFull} have been dropped for notational convenience.  The equations for the Fourier coefficient functions are now discretized in time.  Define the grid functions $\vv^n(y) \approx \vv(y, t^n)$ and $p^n \approx p(y,t^n)$, where $t^n = n \Delta t$ for a (fixed) time step $\Delta t$. An implicit scheme to advance the solution from $t^n$ to $t^{n+1}$, based on backward-Euler time-stepping, is given by
\bse
\label{eq:backwardEuler}
\begin{alignat}{2}
 \rho \f{v_1^{n+1} - v_1^n}{\dt} + i \kx p^{n+1} &= \mu (\partial_y^2 - \kx^2) v_1^{n+1} , \qquad && y>0, \\
 \rho \f{v_2^{n+1} - v_2^n}{\dt} + \partial_y p^{n+1} &= \mu (\partial_y^2 - \kx^2) v_2^{n+1} , \qquad && y>0, \\
(\partial_y^2 - \kx^2)  p^{n+1} &= 0. \qquad && y>0.
\end{alignat}
\ese
Assume that the coefficient functions for the solid variables have been advanced to $t=t^{n+1}$ using an upwind scheme, for example, and that $\bsp^{n+1}$ and $\bss^{n+1}$ are, respectively, the normal and tangential components of the outgoing characteristic variables of the solid at $t\sp{n+1}$.  Using~\eqref{eq:fluidInterfaceConditions}, the boundary conditions for the fluid at $y=0$ take the form
\bse
\label{eq:interfaceNT}
\begin{alignat}{2}
-p^{n+1} + \tau_{22}^{n+1} - \zp v_2^{n+1} &= \bsp^{n+1}, \qquad && y=0,\label{eq:interfaceNormal} \\
\tau_{12}^{n+1} - \zs v_1^{n+1} &= \bss^{n+1}, \qquad && y=0, \label{eq:interfaceTangent}
\end{alignat}
\ese
where the components of the fluid shear stress in~\eqref{eq:interfaceNT} are given by
\begin{align}
\tau_{12}^{n+1} = \mu \left(i \kx v_2^{n+1} + \partial_y v_1^{n+1} \right),
\qquad
\tau_{22}^{n+1} = 2 \mu \partial_y v_2^{n+1}. 
\end{align}
The implicit scheme in~\eqref{eq:backwardEuler} with boundary conditions in~\eqref{eq:interfaceNT} at $y=0$ and boundedness as $y\rightarrow \infty$ determine the grid functions for the fluid at $t\sp{n+1}$ in terms the fluid velocity at $t\sp n$ and the outgoing solid data $(\bsp^{n+1},\bss^{n+1})$.

\newcommand{\Bsp}{B_p}
\newcommand{\Bss}{B_s}

Consider perturbations in the grid functions of the fluid at $t\sp{n+1}$ for $y>0$ subject to perturbations in the interface data $\bsp^{n+1}$ and $\bss^{n+1}$ at $y=0$.  Define
\begin{align*}
&V_1^{n+1} = v_1^{n+1} + \delta V_1,\qquad
V_2^{n+1} = v_2^{n+1} + \delta V_2,\qquad
P^{n+1} = p^{n+1} + \delta P,\qquad y>0,
\end{align*}
and
\begin{align*}
&\Bsp^{n+1} = \bsp^{n+1} + \delta \Bsp,\qquad
\Bss^{n+1} = \bss^{n+1} + \delta \Bss,
\end{align*}
where $(\delta V_1,\delta V_2,\delta P)$ and $(\delta \Bsp,\delta \Bss)$ are small perturbations.  Assuming the fluid velocity at $t\sp n$ is fixed, the variational equations corresponding to~\eqref{eq:backwardEuler} are 
\bse
\begin{alignat}{2}
 \frac{\rho}{\dt} \delta V_1 + i \kx \delta P &= \mu (\partial_y^2 - \kx^2) \delta V_1 , \qquad && y>0, \\
 \frac{\rho}{\dt} \delta V_2 + \partial_y \delta P &= \mu (\partial_y^2 - \kx^2) \delta V_2 , \qquad && y>0, \\
(\partial_y^2 - \kx^2)  \delta P &= 0 , \qquad && y>0.
\end{alignat}
\ese
Solution to these equations that remain bounded as $y\rightarrow \infty$ are
\bse
\begin{align}
\delta V_1(y) &= -\frac{1}{i \kx} \partial_y\delta V_2(y), \\
\delta V_2(y) &= \delta V_0 e^{-\beta y} + \frac{\kx \dt \delta P_0}{\rho}\left(e^{-\kx y}-e\sp{-\beta y}\right), \\
\delta P(y) &= \delta P_0 e^{-\kx y}, 
\end{align}
\ese
where
\[
\delta V_0=\delta V_2(0),\qquad \delta P_0=\delta P(0),\qquad \beta = \left[\kx\sp2+{\rho\over\mu \dt}\right]\sp{1/2}.
\]
Substituting the solution for the perturbations of the fluid variables into the variational equations~\eqref{eq:interfaceNT}
for the interface conditions leads to the linear system
\begin{align}
\left[\begin{array}{cc}
a_{11} & a_{12} \\
a_{21} & a_{22}
\end{array}\right]
\left[\begin{array}{c}
\delta V_0 \\
\delta P_0
\end{array}\right]
=
\left[\begin{array}{c}
\delta \Bsp \\
\delta \Bss
\end{array}\right],
\label{eq:monolithicBC}
\end{align}
where
\[
\begin{array}{ll}
\displaystyle{
a_{11} = -\mu \kx\left(2\gamma+\theta_p\right),
}
& \quad
\displaystyle{
a_{12} = -1+2\Lambda\left(\gamma-1\right),
} \medskip\\
\displaystyle{
a_{21} = i\mu \kx\left(\gamma\sp2+1+\theta_s\gamma\right),
}
& \quad
\displaystyle{
a_{22} = -i\Lambda(\gamma-1)\left(\gamma+1+\theta_s\right).
}
\end{array}
\]
The coefficients $a_{ij}$ in the linear system are defined in terms of the dimensionless parameters
\[
\Lambda=\nu \kx\sp2 \dt,\qquad \gamma={\beta\over \kx}=\sqrt{1+{1\over\Lambda}},\qquad Z_\alpha={\mu \kx \over \bar z_\alpha},\quad\hbox{$\alpha=p$ or $s$}.
\]
The solution of the linear system
\begin{align}
\delta V_0={a_{22}\delta \Bsp-a_{12}\delta \Bss\over a_{11}a_{22}-a_{12}a_{21}},\qquad \delta P_0={a_{11}\delta \Bss-a_{21}\delta \Bsp\over a_{11}a_{22}-a_{12}a_{21}},
\label{eq:variationalSolution}
\end{align}
determines the variation in the interface values of the fluid velocity and pressure in terms of the variations in the outgoing characteristic variables of the solid.

\begin{figure}[h]
\begin{center}
\begin{tikzpicture}
  \useasboundingbox (0,.4) rectangle (8.25,6.5);  
    \draw(0,0)  node[anchor=south west] {\includegraphics[width=.5\textwidth]{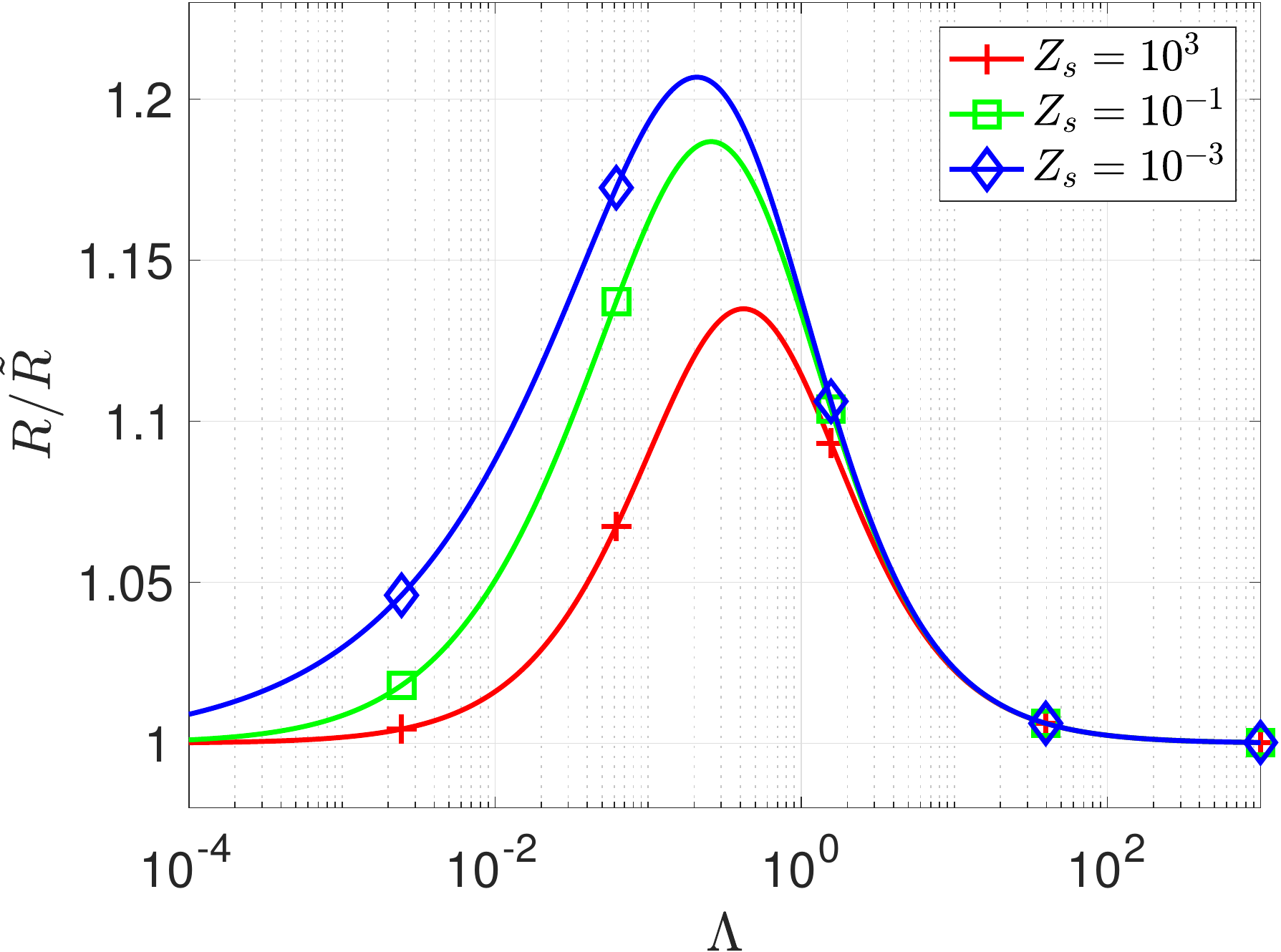}}; 
\end{tikzpicture}
\end{center}
\caption{Behavior of the ratio, $R/\tilde R$, versus the viscous CFL number
        $\Lambda$, for $Z_s=10^{-3}$, $10^{-1}$ and $10^3$ showing that
  $\tilde{R}$ is a reasonable approximation to $R$ over a wide range of $\Lambda$ and $Z_s$.
}
\label{fig:Rplot}
\end{figure}

The AMP algorithm uses impedance-weighted averages to set values for the velocity and pressure at the interface.  For example, the normal component of velocity at the interface is given by~\eqref{eq:AMPinterfaceNormalVelocity}.  In terms of the variational problem, \eqref{eq:AMPinterfaceNormalVelocity} reduces to
\begin{align}
\delta V_0 = -\frac{1}{\zf+\zp} \delta \Bsp,
\end{align}
assuming that the fluid velocity and stress on the right-hand side are held fixed.  In view of the solution in~\eqref{eq:variationalSolution}, we have
\begin{align}
\frac{1}{\zf+\zp} = -{a_{22}\over a_{11}a_{22}-a_{12}a_{21}},
\end{align}
which, after some manipulation, gives
\begin{equation}
z_f=\mu \kx R,\qquad R=2\gamma+{(\gamma+Z_s(\gamma\sp2+1))(1-2\Lambda(\gamma-1))\over\Lambda(\gamma-1)(1+Z_s(\gamma+1))}.
\label{eq:fluidImpedanceFull}
\end{equation}
Of particular interest are the limiting cases when the viscous CFL number, $\Lambda$, is small and large.  A straightforward analysis of the dimensionless parameter $R$ in~\eqref{eq:fluidImpedanceFull} gives
\begin{align}
R \sim
\begin{cases}
\frac{1}{\Lambda}, & \Lambda \rightarrow 0, \\
2, & \Lambda \rightarrow \infty,
\end{cases}
\label{eq:limits}
\end{align}
to leading order.  In view of~\eqref{eq:limits}, we define the relatively simple approximation
\[
\tilde R\eqdef {1\over\Lambda}+2.
\]
The plots in Figure~\ref{fig:Rplot} show that $R/\tilde R\approx1$ over a wide range of values for $\theta_s$, so that the fluid impedance given by
\begin{equation}
\zf=\mu \kx\tilde R=\mu \kx\left({1\over\Lambda}+2\right)={\rho\over \kx\dt}+2\mu \kx,
\label{eq:fluidImpedance}
\end{equation}
is a good approximation of the more complicated form given in~\eqref{eq:fluidImpedanceFull}.
The model problem analysis of Section~\ref{sec:viscousAnalysis} confirms that this choice leads to a stable
scheme.



%
 
  Formula~\eqref{eq:fluidImpedance} provides the generic form of the fluid impedance we use, but it remains
  to make a choice for $\kx$ so that the approximation can be used for a discrete approximation in physical coordinates (as
  opposed to the Fourier transformed space).
Note that in a discrete approximation, the possible wave numbers $\kx$ appearing in~\eqref{eq:fluidImpedance} are
bounded by approximately $1/h$, where $h$ is a measure of the grid spacing in the tangential direction.
For the present model problem with the
the pseudo-spectral approximation~\eqref{eq:discreteFourier}, for example, we have $|\kx|\le \pi/\dx$, while a second-order difference
approximation would roughly imply $|\kx| \le 2/\dx$. 
Experience~\cite{rbinsmp2017} shows that added-damping instabilities
are generally caused by relatively high-frequency modes on the grid, and this suggests taking $\kx=1/h$ which leads to 
a definition of the fluid impedance of the form
\ba
\zf \eqdef  \Cam \Bigl({\rho\, h \over \dt}\Bigr) + \Cad \Bigl(\f{\mu}{h}\Bigr),
\ea
as was done in~\eqref{eq:fluidImpedanceAlg}.  The extensive numerical results
in Section~\ref{sec:numericalResults} and~\cite{fibr2019} confirm that this is an appropriate choice, and furthermore that
the scheme is rather insensitive to the choice of $h$, $\Cam$ and $\Cad$.

%% file: texFiles/rectangularModelProblemCartoonInfinite.tex
{
\newcommand{\lbfont}{\small}
\def\ysb{-.15} 
\def\ysa{-2.5}   
\def\ya{0} 
\def\yb{2.5} 
\def\xL{8}
\begin{figure}[hbt]
	\newcommand{\textFont}{\normalss}
	\begin{center}
		\begin{tikzpicture}[scale=.9]
		\useasboundingbox (0,-2.75) rectangle (\xL,2.5);  

                \draw[fill=red!20,draw=red!20] (0,\ysa) rectangle (\xL,\ysb);
                \draw[-,thick,red,line width=2pt] 
                (0  ,\ysa) -- 
                (0  ,\ysb) -- 
                (\xL,\ysb) --
                (\xL,\ysa);
		\draw[thick,black] (4,-1.2) node {solid: $\OmegaS$};

		\draw[thick,fill=blue!20,draw=blue!20,line width=2pt] (0,\ya) rectangle (\xL,\yb);
                \draw[thick,black] (.5*\xL,.5*\ya+.5*\yb) node {fluid: $\OmegaF$};
		\draw[thick,black] (4,\ya) node[anchor=south] {interface: $\Gamma$};
                \draw[-,thick,blue,line width=2pt] 
                (0  ,\yb) -- 
                (0  ,\ya) -- 
                (\xL,\ya) --
                (\xL,\yb);


		\draw[] (0  ,\ya) node[anchor=east,black,yshift=-4pt] {\lbfont$y=0$};

		\draw (0  ,\ysa) node[anchor=north,black,yshift=-4pt] {\lbfont$x=0$};
		\draw (\xL,\ysa) node[anchor=north,black,yshift=-4pt] {\lbfont$x=L$};
		\draw (0,  \ysa) node[anchor=east,black] {\lbfont$y\rightarrow -\infty$};
		\draw (0,   \yb) node[anchor=east,black] {\lbfont$y\rightarrow +\infty$};
		%
		%
		\end{tikzpicture}
	\end{center}
	\caption{The rectangular geometry for the viscous model problem.} 
        \label{fig:rectangularModelProblemCartoonInfinite}
\end{figure}
}

%% file: texFiles/viscousAnalysis-new.tex
\section{Stability analysis of a viscous model problem} \label{sec:viscousAnalysis}

The stability of the AMP algorithm is explored in the context of two FSI model problems.  The first model problem, discussed in this section, involves a viscous incompressible (Stokes) fluid in contact with a simplified elastic solid.  This analysis extends the work in~\cite{fib2014} to the case of a viscous fluid where both added-mass and added-damping effects are important, and for an IMEX-type scheme in the fluid.  The stability analysis of a second FSI model problem involving an inviscid incompressible fluid is carried out in Section~\ref{sec:inviscidAnalysis}.

We will compare the stability of the AMP scheme to that of the traditional partitioned (TP)
scheme and the anti-traditional partitioned (ATP) scheme. 
In the TP scheme, the solid provides a Dirichlet (no-slip)
boundary condition for the fluid, and then the fluid supplies a Neumann (traction) 
boundary condition for the solid. 
The ATP scheme reverses the role of the solid and fluid. 
In this scheme, the solid provides a Neumann (traction) boundary condition for the fluid
and the fluid supplies a Dirichlet (no-slip) boundary condition for the solid.


\subsection{Viscous model problem}


The viscous model problem analyzed here is similar to the one discussed in Section~\ref{sec:monolithic} and illustrated in Figure~\ref{fig:rectangularModelProblemCartoonInfinite}.  An incompressible Stokes fluid satisfies the system of equations in~\eqref{eq:modelFluidEquations} for $\xv\in\OmegaF$.  The solid satisfies~\eqref{eq:modelSolidEquations} for $\xv\in\OmegaS$, but with $\lambdas$ set equal to $-\mus$.  This choice simplifies the system of equations for the solid somewhat since the compressive and shear wave speeds both equal $\cb=\sqrt{\mus/\rhos}$, i.e.~the shear wave speed which is particularly relevant for the viscous model problem.  It is convenient to consider the hyperbolic equations for the solid in characteristic form.  These equations are
\bse
\label{eq:advFull}
\begin{alignat}{2}
\partial_t \aone - \bar{c} \partial_y \aone &= \bar{c} (\partial_x d - \partial_x \btwo),
\qquad & \xv \in \OmegaS, \label{eq:advA1}\\
\partial_t \bone + \bar{c} \partial_y \bone &= \bar{c} (\partial_x \atwo - \partial_x d),
\qquad & \xv \in \OmegaS, \label{eq:advB1}\\
\partial_t \atwo - \bar{c} \partial_y \atwo &= \bar{c} \partial_x \bone, 
\qquad & \xv \in \OmegaS, \label{eq:advA2}\\
\partial_t \btwo + \bar{c} \partial_y \btwo &= - \bar{c} \partial_x \aone,
\qquad & \xv \in \OmegaS, \label{eq:advB2}\\
\partial_t d &= 0,
\qquad & \xv \in \OmegaS, \label{eq:advD}
\end{alignat}
\ese
where
\begin{align}
\am = \sigmas_{m2} + \bar{z} \vs_m , \qquad
\bm = \sigmas_{m2} - \bar{z} \vs_m, \qquad
\hbox{$m=1$ or $2$,} \qquad \bar{z}=\rhos\cb,
\label{eq:CharDef}
\end{align}
are the variables associated with the incoming $(\aone,\atwo)$ and outgoing $(\bone,\btwo)$ characteristics at the interface, and
\begin{align}
d = \sigmas_{11} + \sigmas_{22}.
\label{eq:CharDefII}
\end{align}
The interface is linearized about a flat surface $\Gamma$ given by $y=0$, and the matching conditions between the fluid and the solid are given in~\eqref{eq:modelInterface}.

\subsection{Discretization}

The discretization of the equations in the $x$-direction follows the approach used previously in Section~\ref{sec:monolithic}.  The equations for the fluid are transformed using the finite Fourier series in~\eqref{eq:discreteFourier}, which results in a system of equations for the corresponding Fourier coefficient functions given in~\eqref{eq:modelFluidEquations}.  These equations are then discretized in time using an IMEX-type scheme given by
\bse
\label{eq:IMEXfluid}
\begin{align}
v_1^{n+1} &= v_1^n 
- \frac{i \kx \dt}{\rho} p^{n}
+ \nu \dt \left( - \kx^2 + \partial_y^2  \right) v_1^{n+1}, \label{eq:viscousV1}\\
v_2^{n+1} &= v_2^n 
- \frac{\dt}{\rho} \partial_y p^{n}
+ \nu \dt \left( - \kx^2 + \partial_y^2\right) v_2^{n+1}, \label{eq:viscousV2}\\
\left(-\kx^2 + \partial_y^2\right) p^{n+1} &= 0 .\label{eq:viscousP}
\end{align}
\ese
Here, $v_1^n(y)$, $v_2^n(y)$ and $p^n(y)$ approximate $v_1(y,t\sp n)$, $v_2(y,t\sp n)$ and $p(y,t\sp n)$, respectively, at $t\sp n=n\dt$ for a fixed time step $\dt$.  Recall that $\kx=2\pi k/L$, and that the hats and $k$ subscripts on the Fourier coefficients have been suppressed.  Note that the components of the fluid velocity are advanced in time using (implicit) backward Euler for the viscous terms and (explicit) forward Euler for the pressure gradient terms.  An elliptic equation is solved at each time step to update the pressure.  It is convenient to keep the discrete equations for the fluid variables continuous in~$y$, and we assume that solutions are bounded as $y\rightarrow\infty$.

The characteristic equations for the solid in~\eqref{eq:advFull} are similarly transformed using the finite Fourier series in~\eqref{eq:discreteFourier}, and then the resulting equations are discretized in time and space using an upwind-type scheme given by
\bse
\label{eq:upwindScheme}
\begin{align}
a_{1,j}^{n+1} &= a_{1,j}^n + \ly (a_{1,j+1}^n-a_{1,j}^n) + i \lx (d_{j}^{n+1} - b_{2,j}^{n+1}), \label{eq:schemeA1}\\
b_{1,j}^{n+1} &= b_{1,j}^n - \ly (b_{1,j}^n-b_{1,j-1}^n) + i \lx (a_{2,j}^{n+1} - d_{j}^{n+1}), \label{eq:schemeB1}\\
a_{2,j}^{n+1} &= a_{2,j}^n + \ly (a_{2,j+1}^n-a_{2,j}^n) + i \lx b_{1,j}^{n+1}, \label{eq:schemeA2}\\
b_{2,j}^{n+1} &= b_{2,j}^n - \ly (b_{2,j}^n-b_{2,j-1}^n) - i \lx a_{1,j}^{n+1}, \label{eq:schemeB2}\\
d_j^{n+1} &= d_j^n, \label{eq:schemeD}
\end{align}
\ese
where, for example, $a_{1,j}^{n+1}\approx a_1(y_j,t_n)$ with $y_j = j \dy$ and $t\sp n=n\dt$, and where $\lx = \bar{c}\kx \dt$ and $\ly = \bar{c} \dt / \dy$.  The grid in the $y$-direction is collocated about the interface at $y=0$.  The terms involving transverse derivatives are treated implicitly to stabilize the pseudo-spectral approximation. For reference, the solid velocity and stress are related to the characteristic variables by
\bse
\begin{align}
\vs_{1,j}^n = \frac{1}{2 \bar{z}} \left( a_{1,j}^n - b_{1,j}^n \right), \qquad
\vs_{2,j}^n = \frac{1}{2 \bar{z}} \left( a_{2,j}^n - b_{2,j}^n \right), \label{eq:solidGridVelocity}\\
\sigmas_{12,j}^n = \frac{1}{2 } \left( a_{1,j}^n + b_{1,j}^n \right), \qquad 
\sigmas_{22,j}^n = \frac{1}{2 } \left( a_{2,j}^n + b_{2,j}^n \right). \label{eq:solidGridStress}
\end{align}
\ese
We assume bounded solutions of~\eqref{eq:upwindScheme} as $y_j\rightarrow-\infty$.

\subsection{Interface coupling}

We explore the stability of partitioned schemes for the model problem that use different interface coupling approaches.  For any of the approahes, corresponding to the AMP, TP, and ATP schemes, the discrete equations require a certain number of boundary conditions at the interface.  For example, the evolution of the fluid equations in~\eqref{eq:IMEXfluid} require three boundary conditions on the interface, $y=0$, to determine the interface velocity and pressure.  Similarly, the evolution of the solid equations in~\eqref{eq:upwindScheme} require two boundary conditions at $y=0$ corresponding to the two incoming characteristic variables.

We first describe the coupling based on the AMP interface conditions given in Section~\ref{sec:interface}. We assume the fluid and solid solutions are known at time $t^n$.  The solid variables are advanced first to $t^{n+1}$ on grid points $j=0,-1,-2,\hdots$ using the evolution equations in~\eqref{eq:upwindScheme}.  The solid interface velocity and stress are computed using 
\begin{align}
\vs_{m,0}^{n+1} = \frac{1}{2 \zb} \left(a_{m_0}^{n+1} - b_{m,0}^{n+1}  \right),
\qquad
\sigmas_{m2,0}^{n+1} = \frac{1}{2} \left(a_{m,0}^{n+1} + b_{m,0}^{n+1}  \right),
\qquad m=1,2.
\end{align}
The fluid velocity is advanced using~\eqref{eq:viscousV1}--\eqref{eq:viscousV2}.  Two boundary conditions are required at $y=0$ to obtain the fluid velocity at $t^{n+1}$.  The condition on the outgoing solid tangential characteristic in~\eqref{eq:AMPTangentialCharacteristic} becomes
\begin{align}
\mu \left( i \kx v_2^{n+1} + \partial_y v_1^{n+1} \right) - \zb v_1^{n+1} 
= \sigmas_{12,0}^{n+1} - \zb \vs_{1,0}^{n+1},
\qquad y=0. \label{eq:fracTangential} 
\end{align}
The normal component of the velocity is projected to obtain the proper limiting behaviors for heavy and light solids.  This condition, taken from~\eqref{eq:normalVelociityFix}, reduces to
\begin{align}
&v_2^{n+1} = \frac{\zf}{\zf+\zb}  V^p(v_2^{n+1}) + \frac{\zb}{\zf+\zb} \vs_{2,0}^{n+1}, \qquad y=0,
\label{eq:fracProject}
\end{align}
where the fluid impedance is given by
\begin{align}
\zf = \frac{\rho}{\kx \dt} + 2 \mu \kx,
\end{align}
according to the derivation in Section~\ref{sec:monolithic}.  The predicted velocity, $V^p(v_2^{n+1})$, in~\eqref{eq:fracProject} is given by
\begin{align}
V^p(v_2^{n+1}) = v_2^{n} 
- \frac{\dt }{\rho} \partial_y p^{n}
- \nu \dt \left(\kx^2 v_2^{n+1} + i \kx \partial_y v_1^{n+1}\right) , \qquad y=0.
\label{eq:predictedVelocity}
\end{align}
This definition is analogous to that in~\eqref{eq:veDef}, but with the substitution $\partial_y v_2^{n+1} = -i k_x v_1^{n+1}$ noting that $\nabla\cdot\vv=0$ on the boundary.  The pressure is updated using~\eqref{eq:viscousP} along with the AMP pressure condition described in~\eqref{eq:AMPpressureBC}.  For the present scheme, this condition reduces to 
\begin{align}
-p^{n+1}  + \frac{\zb \dt}{\rho} \partial_y p^{n+1} 
=\;& \sigmas_{22,0}^{n+1} + 2 i k_x \mu  v_1^{n+1} 
 \nonumber \\
& \qquad
- \zb \dt \left[\dot{\vs}_{2,0}^{n+1} + \nu\left( \kx^2 v_2^{n+1} + i \kx \partial_y v_1^{n+1} \right)\right],
\qquad y=0,
\label{eq:fracPressure} 
\end{align}
again using $\partial_y v_2^{n+1} = -i k_x v_1^{n+1}$.  The acceleration of the solid on the interface, denoted by $\dot{\vs}_{2,0}\sp{n+1}$ in~\eqref{eq:fracPressure}, is taken to be
\begin{align}
\dot{\vs}_{2,0}^{n+1} = \frac{\vs_{2,0}^{n+1}-\vs_{2,0}^{n}}{\dt}.
\end{align}
After solving for the fluid velocity and pressure, interface quantities from the fluid are obtained using
\bse
\begin{align}
v_{m,f}^{n+1} &= v_m^{n+1}, \qquad m=1,2,\\
p_f^{n+1} &= p^{n+1}, \\
\sigma_{12,f}^{n+1} &= \mu \left(\partial_y v_1^{n+1} + i \kx v_2^{n+1} \right), \\
\sigma_{22,f}^{n+1} &= -p^{n+1} + 2 \mu \partial_y v_2^{n+1},
\end{align}
\ese
where all fluid quantities on the right-hand side are evaluated at $y=0$.
The interface velocity and traction are projected from fluid and solid values 
using~\eqref{eq:AMPinterfaceVelocity} and~\eqref{eq:AMPinterfaceTraction}.
These equations reduce to
\bse
\label{eq:projections}
\begin{align}
v_m^I &= \frac{\zf}{\zf+\zb} v_{m,f}^{n+1} 
+ \frac{\zb}{\zf+\zb} \vs_{m,0}^{n+1}
+ \frac{1}{\zf+\zb} ( \sigmas_{m2,0}^{n+1} - \sigma_{m2,f}^{n+1}  ), \label{eq:project1}\\
\sigma_{m2}^I &= \frac{\zf^{-1}}{\zf^{-1}+\zb^{-1}} \sigma_{m2,f}^{n+1}
+ \frac{\zb^{-1}}{\zf^{-1}+\zb^{-1}} \sigmas_{m2,0}^{n+1}
+ \frac{1}{\zf^{-1}+\zb^{-1}} ( \vs_{m,0}^{n+1} - v_{n,f}^{n+1} ), \label{eq:project2}
\end{align}
\ese
where $m=1,2$.  Finally, the ghost points at $j=1$ for the incoming solid characteristics are set using
\begin{align}
a_{m,1}^{n+1} &= \sigma_{m2}\sp I + \zb v_m\sp I, \qquad m=1,2,\label{eq:a1CharBoundary}
\end{align}
which is a first-order accurate approximation (consistent with the order of accuracy of the upwind scheme).

We next consider the coupling conditions for the TP and ATP schemes.  These conditions can be obtained from the coupling conditions for the AMP scheme in the limits of heavy ($\zb\rightarrow\infty$) and light ($\zb\rightarrow0$) solids.  For the AMP algorithm, the fluid velocity and pressure conditions are given in~\eqref{eq:fracTangential}, \eqref{eq:fracProject} and~\eqref{eq:fracPressure}, while the final interface values are defined by the projections in~\eqref{eq:projections}.  For the TP algorithm ($\zb\rightarrow\infty$), the AMP conditions in~\eqref{eq:fracTangential} and~\eqref{eq:fracProject} reduce to Dirichlet conditions on the fluid velocity given by
\begin{align}
v_m^{n+1} 
&=  \vs_{m,0}^{n+1}, \qquad y=0,\qquad m=1,2.\label{eq:TPTangential}
\end{align}
The pressure condition in equation~\eqref{eq:fracPressure} becomes a Neumann condition given by
\begin{align}
\partial_y p^{n+1} 
&= 
- \dot{\vs}_{2,0}^{n+1} 
-\nu \left( \kx^2 v_2^{n+1} + i \kx \partial_y v_1^{n+1} \right), \qquad y=0.
\label{eq:TPPressure} 
\end{align}
For the TP scheme, the interface velocity is taken to be the solid velocity, $v_m^I = \vs_{m,0}^{n+1}$, and the interface traction is taken to be the fluid traction, $\sigma_{m2}^I = \sigma_{m2,f}^{n+1}$, $m=1.2$.

For the ATP scheme, we consider the light-solid limit ($\zb\rightarrow 0$) of the AMP conditions.  In this limit, the condition on the outgoing solid tangential characteristic in~\eqref{eq:fracTangential} reduces to a Neumann condition for the velocity given by
\begin{align}
\mu \left( i k_x v_2^{n+1} + \partial_y v_1^{n+1} \right)
&= \sigmas_{12,0}^{n+1}, \qquad y=0,\label{eq:ATPV1} 
\end{align}
while the condition in~\eqref{eq:fracProject} becomes
\begin{align}
v_2^{n+1} = V^p(v_2^{n+1}), \qquad y=0,
\label{eq:ATPV2}
\end{align}
where $V^p(v_2^{n+1})$ is given by~\eqref{eq:predictedVelocity}.  Using~\eqref{eq:viscousV2}, it can be shown that the condition in~\eqref{eq:ATPV2} can be replaced by
\begin{align}
\partial_y \left( i \kx v_1^{n+1} + \partial_y v_2^{n+1} \right)
&= 0, \qquad y=0, \label{eq:ATPDiv}
\end{align}
which is equivalent to setting the fluid velocity to be divergence-free on the interface.  For the ATP scheme, the pressure condition in~\eqref{eq:fracPressure} reduces to 
\begin{align}
-p^{n+1} + 2 \mu \partial_y v_2^{n+1} 
&= \sigmas_{22,0}^{n+1}, \qquad y=0. \label{eq:ATPP}
\end{align}
For the ATP scheme, the interface velocity is taken to be the fluid velocity, $v_m^I = v_{m,f}^{n+1}$, and the interface traction is taken to be the solid traction, $\sigma_{m2}^I = \sigmas_{m2,0}^{n+1}$, $m=1,2$.

%
%

\subsection{Stability analysis}

In order to assess the stability of the AMP, TP and ATP schemes, we search for normal mode solutions to the discrete evolution equations.  In the fluid, solutions are of the form
\begin{align}
v_m^n(y) = A^n \tilde{v}_m(y), \qquad p^n(y) = A^n \tilde{p}(y),\qquad m=1,2,
\label{eq:fluidSoln}
\end{align}
where $A$ is an amplification factor.  Substituting~\eqref{eq:fluidSoln} into~\eqref{eq:IMEXfluid} and integrating gives
\bse
\label{eq:fluidSolnII}
\begin{align}
\tilde{v}_1(y) 
  &= v_{1,f}^0 e^{-\gamma \kx y} 
    - \frac{i p_f^0}{\mu \kx A (\gamma^2 - 1)} \left( e^{-\kx y} - e^{-\gamma \kx y} \right), \\
\tilde{v}_2(y) 
  &= v_{2,f}^0 e^{-\gamma \kx y}
    + \frac{p_f^0}{\mu \kx A (\gamma^2 - 1)}\left( e^{-\kx y} - e^{-\gamma \kx y} \right), \\
\tilde{p}(y) &= p_f^0 e^{-\kx y},
\end{align}
\ese
where
\begin{align}
\gamma = \sqrt{1+\frac{1}{\Lambda}\left(\frac{A-1}{A}\right)}, \qquad
\Lambda = \nu k_x^2 \dt.
\end{align}
Here, $\Lambda$ represents the viscous CFL number and we have imposed boundedness of the solution in~\eqref{eq:fluidSolnII} as $y\rightarrow\infty$.
The constants, $v_{m,f}^0$ and $p_f^0$, are obtained by imposing the appropriate boundary conditions at $y=0$, namely, \eqref{eq:fracTangential},\eqref{eq:fracProject} and~\eqref{eq:fracPressure} for the case of the AMP scheme.  For the TP scheme, the three constraints are the two boundary conditions for the components of the velocity in~\eqref{eq:ATPV1} and the condition on the pressure in~\eqref{eq:TPPressure}, while the ATP scheme uses the boundary conditions in~\eqref{eq:ATPV1}, \eqref{eq:ATPDiv} and~\eqref{eq:ATPP}.



Having found solutions for the velocity and pressure of the fluid, these solutions can be used (along with the appropriate boundary conditions at $y=0$ for the AMP, TP or ATP coupling) to eliminate the fluid variables on the boundary in favor of the solid variables. 
The issue of stability, then, reduces to examining the behavior of the evolution equations for the solid with the appropriate boundary conditions.  Solutions of these evolution equations are sought in the form
\begin{equation}
\av_j^n =  \phi^j A^n \tilde{\rv},\qquad \av_j^n = [a_{1,j}^n\,,\,b_{1,j}^n\,,\,a_{2,j}^n\,,\,b_{2,j}^n\,,\,d_j^n]^T,
\label{eq:normalModes}
\end{equation}
where $\phi$ is a spatial eigenvalue and $\tilde{\rv}$ is a constant eigenvector.
The scheme is said to be weakly stable if there are no non-trivial solutions with $|A| > 1$.  
Our strategy for determining regions of stability will be to search for unstable modes with 
$|A| > 1$, and then identify regions of the parameter space where no non-trivial solutions exist.  
To do this, we begin by finding the general solution for the spatial grid functions satisfying 
the discrete equations and regularity condition as $j\rightarrow-\infty$, assuming $|A| > 1$.  
We then apply the conditions at the interface to determine whether non-trivial solutions exist.

When the normal mode ansatz in~\eqref{eq:normalModes} is substituted into the evolution equations for the solid in~\eqref{eq:upwindScheme}, a homogeneous system arises of the form
\newcommand{\SolidMatrixB}{\mathcal{F}}
\begin{align}
\SolidMatrixB(\phi) \tilde{\rv} = 0.
\end{align}
The matrix $\SolidMatrixB$ is given by
\begin{align*}
\SolidMatrixB(\phi) = \left[\begin{array}{ccccc}
\eta(\phi) & 0 &  0& -i A \lx & i A \lx  \\
0 & \eta(1/\phi) & i A \lx & 0 & -i A \lx \\
0 & i A \lx & \eta(\phi) & 0 & 0 \\
-i A \lx & 0 & 0 & \eta(1/\phi) & 0 \\
0 & 0 & 0 & 0 & 1-A
\end{array}\right],
\end{align*}
where
\begin{align}
\eta(\phi) = 1-A+\ly (\phi-1), \qquad \lx = \bar{c}\kx \dt, \qquad \ly = \bar{c} \dt / \dy.
\end{align}
The determinant of $\SolidMatrixB$ is given by
\begin{align}
f(\phi) \defeq \det (\SolidMatrixB(\phi)) = (1-A) \left(\eta(\phi) \eta(1/\phi) + (A \lx)^2\right)^2.
\label{eq:SolidDetCond}
\end{align}
The system is singular when $f(\phi) = 0$, which occurs when $A=1$ or $\eta(\phi) \eta(1/\phi) + (A \lx)^2 = 0$.
Since we are searching for unstable modes with $\vert A\vert>1$, we are interested in the latter 
case which leads to roots given by
\newcommand{\bParam}{\xi}
\begin{align}
\phi_\pm = \bParam \pm \sqrt{\bParam^2 - 1}, \qquad 
\bParam = 1 - \frac{(A \lx)^2 + (1-A)^2}{2 \ly (1-A-\ly)}.
\end{align}
The product of the roots is equal to one ($\phi_+ \phi_- = 1$).
Since we are searching for solutions that are bounded as $j\rightarrow -\infty,$ 
we are only interested in the root with modulus greater than one.

\medskip
\noindent\textbf{Lemma}: \textit{If $|A|>1$ and if $\lx$ and $\ly$ are chosen to satisfy a CFL 
condition, then there is precisely one root, either $\phi_+$ or $\phi_-$, denoted by 
$\phi_*$ that has magnitude strictly greater than one, i.e.~$|\phi_*| > 1.$}
\medskip

\noindent This result follows from a similar argument to that given 
in~\cite{DifferenceApproximationsForTheInitialBoundary}.  
We first consider the scheme applied to the pure initial-value problem (Cauchy problem).  
Setting $\phi=e\sp{i\vartheta}$ in~\eqref{eq:SolidDetCond}, we determine a region of 
the $(\lx,\ly)$ plane for which $\vert A\vert\le1$ for all $\vartheta\in[0,2\pi]$.  
This region is found numerically as discussed in Section~\ref{sec:CFLRegionCauchy}
and shown in Figure~\ref{fig:viscousCFL}, and it corresponds 
to a region where the Cauchy problem is (Fourier) stable.  It is also shown that this stability 
region includes a region of $(\lx,\ly)$ satisfying a reasonable CFL restriction, 
namely $\lx\sp2+\ly\sp2\le1$.
Next, since $|A| \le 1$ when $|\phi| = 1$, we have that $|\phi| \ne 1$ when $|A| > 1$.  
Thus, if $|A| > 1$ and if $(\lx,\ly)$ remains within the CFL restriction, then $\phi$ cannot 
cross the unit circle, $|\phi|=1$, as $(\lx,\ly)$ vary.  
It is therefore only necessary to prove that the lemma holds for one set of parameters.  
For $\lx=0$, the discretization reduces to four uncoupled upwind schemes for linear advection. 
In this case, equation~\eqref{eq:SolidDetCond} is equivalent to 
$\eta(\phi) \eta(1/\phi) = 0,$
which has solutions
$\phi_+ = (A-1+\ly)/\ly$ and $\phi_- = 1/\phi_+$.
When $|A|>1$ and $\ly \in (0,1]$, $|\phi_+|>1$ and therefore $\phi_* = \phi_+.$ 
Thus, the condition holds for all $(\lx,\ly)$ provided the CFL condition is satisfied.

When $\phi=\phi_*,$ there are two eigenvectors, namely
\begin{align}
\tilde{\rv}_1 = 
\left[
\frac{\eta(1/\phi_*)}{i A\lx},\, 0,\, 0,\, 1,\, 0
\right]^T, \qquad 
\tilde{\rv}_2 = 
\left[
0,\, -\frac{\eta(\phi_*)}{i A\lx},\, 1,\, 0,\, 0
\right]^T
.
\end{align}
The solution which remains bounded as $j \rightarrow -\infty$ is given by
\begin{align*}
a_{1,j}^{n} = k_1 \frac{ \eta(1/\phi_*) }{i A \lx} \phi_*^{j} A^n, \qquad
b_{1,j}^{n} = -k_2 \frac{ \eta(\phi_*) }{i A \lx} \phi_*^{j} A^n, \qquad
a_{2,j}^{n} = k_2 \phi_*^{j} A^n, \qquad
b_{2,j}^{n} = k_1 \phi_*^{j} A^n,
\end{align*}
where $k_1$ and $k_2$ are constants to be determined by the two interface conditions in~\eqref{eq:a1CharBoundary}.  The application of these interface conditions leads to another homogeneous system of equations given by
\begin{align}
\dMat(A) \kv = \left[\begin{array}{cc}
\dMat_{11} & \dMat_{12} \\ 
\dMat_{21} & \dMat_{22}
\end{array}\right]
\left[\begin{array}{c}
k_1 \\ k_2
\end{array}\right]
= 0.
\label{eq:StabilitySystemViscous}
\end{align}
The coefficients of the matrix $\dMat$
are provided in Section~\ref{sec:ComponentsOfDmatViscous} for the AMP, TP, and
ATP schemes.
Solutions for the amplification factor $A$ are roots of the transcendental 
equation given by
\begin{align}
g(A) \defeq \det(\dMat(A)) = \dMat_{11} \dMat_{22} - \dMat_{12} \dMat_{21} = 0.
\label{eq:ViscousDet}
\end{align}
These roots dependent on the choice of interface coupling and four dimensionless parameters ($\Lambda, Z, \lx, \ly$), where
\begin{align}
Z = \frac{\mu k_x}{\zb}.
\end{align}

Proving stability of the partitioned scheme for a choice of the interface coupling and dimensional parameters is equivalent to 
showing that no roots of~\eqref{eq:ViscousDet} exist such that $|A|>1$.  The number of roots with $|A|>1$ can be assessed using the argument principle.  Define
\newcommand{\gi}{{G}}
\begin{align}
\mathcal{P} \defeq \frac{1}{2 \pi i}\oint_{\vert\zeta\vert=1} \frac{\gi\sp\prime(\zeta)}{\gi(\zeta)} \, d\zeta,
\qquad
\gi(\zeta)=g(1/\zeta). 
\end{align}
There are branch points of $\gi(\zeta)$ in the region $\vert\zeta\vert>1$, and a single-valued branch of $\gi(\zeta)$ can be defined so that its branch cuts lie outside the unit disk.  The only singularity of $\gi(\zeta)$ in the region $|\zeta|\le 1$ is a pole of order~$2$ at the origin, and thus $\mathcal{P} = N-2$, where $N$ corresponds to the number of roots of $g(A)$ with $\vert A\vert>1$.

\begin{figure}[h]
\centering
\includegraphics[width=.48\textwidth]{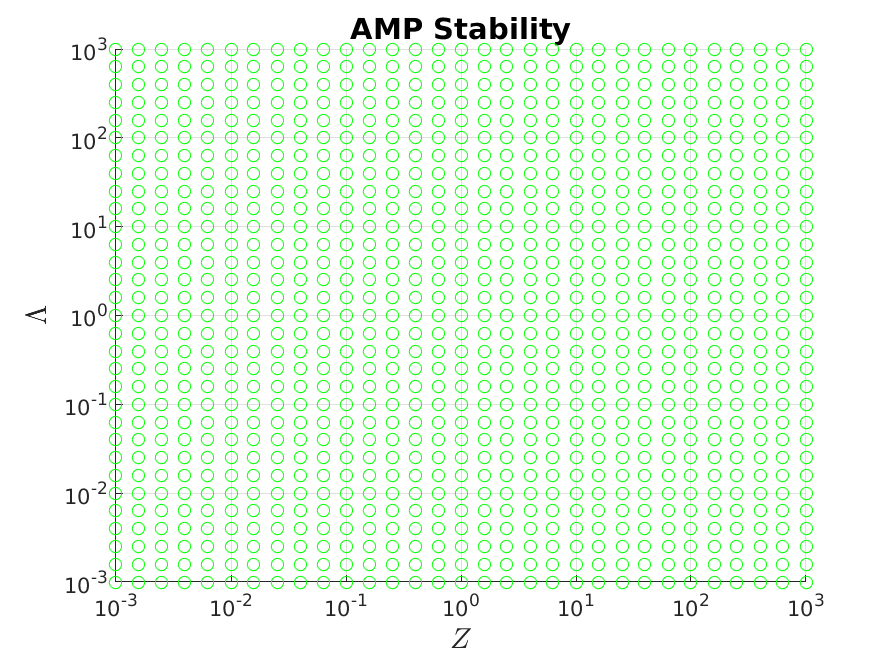}
\includegraphics[width=.48\textwidth]{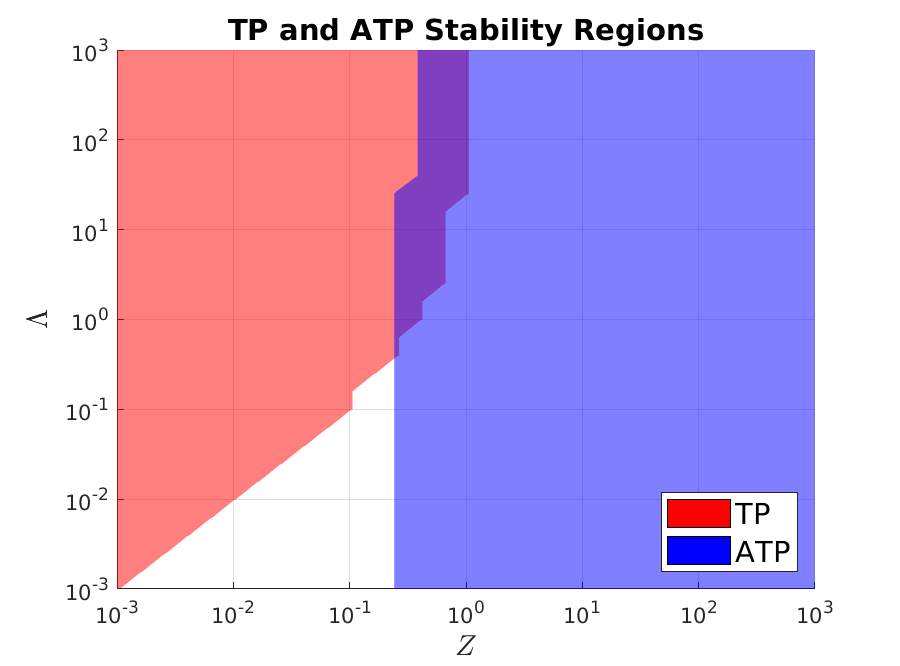}
\caption{
  Left: Green circles represent points for which 
  the AMP algorithm is stable in the CFL region $\lx^2+\ly^2\le 1$.
  Right: stability regions for the TP (red) and ATP (blue) schemes.
  \label{fig:argPrinciple}
}
\end{figure}

An analytic evaluation of the integral for $\mathcal{P}$ is unavailable, and so we consider a numerical evaluation.  The four-dimensional parameter space ($Z$,$\Lambda$,$\lx$,$\ly$) is discretized on a 31$\times$31$\times$20$\times$20 array.  The parameters $Z$ and $\Lambda$ are equally spaced on a logarithmic scale on the interval $[10^{-3}, 10^{3}]$, while $\lx$ and $\ly$ are equally spaced on the interval $[0.05, 0.95]$.  At each grid point, $\mathcal{P}$ is computed numerically with $\vert\mathcal{P}+2\vert\le\delta$ corresponding to stability, where $\delta$ is a small parameter taken to be $10\sp{-5}$.  The results of this computation are shown in Figure~\ref{fig:argPrinciple} for the AMP, TP and ATP schemes.  A grid point in the $(\Lambda,Z)$ plane is marked as stable if the computations of $\mathcal{P}$ for all values of $\lx$ and $\ly$ in the search region are stable.  The point is marked as unstable otherwise.  The results shown in the left plot indicate that the AMP scheme applied to the viscous model problem is stable for all points in the $(\Lambda,Z)$ plane, whereas the results shown in the right plot indicate that the TP and ATP schemes have large regions of instability.  For example, the region in red shows the stable region for the TP scheme, which occurs for heavy solids ($Z$ small) and coarser meshes ($\Lambda$ large).  The stability region for the ATP scheme shown in blue corresponds to light solids ($Z$ large).  The following theorem summarizes the results for the AMP scheme.

\medskip
\noindent\textbf{Theorem}: \textit{The AMP scheme applied to viscous model problem is weakly stable $|A| \le 1$ provided $\lx^2 + \ly^2 \le 1$, which gives the usual CFL-type time-step restriction}
\begin{align}
\dt \le \frac{1}{\cb} \left[ \frac{1}{\dy^2} + k_x^2\right]^{1/2}.
\end{align}
\textit{This is a sufficient but not a necessary condition.  The proof follows from the argument principle and a numerical evaluation of $\mathcal{P}$.}

%% file: texFiles/inviscidAnalysis-new.tex
\section{Stability analysis of an inviscid model problem} \label{sec:inviscidAnalysis}

In this section, the stability of the AMP, TP and ATP algorithms are explored for a two-dimensional FSI model problem involving an inviscid incompressible fluid.  A similar stability analysis was performed in~\cite{fib2014} for a first-order accurate scheme on a staggered grid.  The algorithms considered here are second-order accurate, as in the full AMP algorithm discussed in Section~\ref{sec:algorithm}, and thus the present analysis extends the previous results for less dissipative partitioned schemes.

\subsection{Inviscid model problem}


The fluid is assumed to be inviscid and incompressible, and it occupies the region $\OmegaF$ given by $0<x<L$ and $0<y<H$.  The solid lies below the fluid in the semi-infinite region $\OmegaS$ given by $0<x<L$ and $y<0$, see Figure~\ref{fig:rectangularModelProblemCartoon}.  The solid is taken to be an ``acoustic'' solid that supports a motion in the vertical direction only.  The fluid-solid interface $\Gamma$ is linearized about the flat surface, $y=0$, and it is assumed that the vertical components of the velocity and stress between the fluid and solid match along $\Gamma$.  The fluid pressure is taken to be zero along the top boundary, $y=H$, and the solutions for the fluid and solid are assumed to be periodic in the $x$-direction with period equal to~$L$.  The governing equations for this FSI model problem are
\bse
\begin{align}
\text{Fluid: } &\begin{cases}
\rho \partial_t \vfy + \partial_y p = 0, \quad & \xv\in\OmegaF, \\
(\partial_x^2 + \partial_y^2) p =0, \quad & \xv\in\OmegaF, \\
p = 0, & x\in(0,L), \quad y=H,
\end{cases} \label{eq:fluidModelInviscid} \\
\text{Solid: } &\begin{cases}
\rhos \partial_t \vsy = \partial_x \ssx + \partial_y \ssy, \quad & \xv\in\OmegaS, \\
\partial_t \ssx = \rhos\cp^2 \partial_x \vsy, & \xv\in\OmegaS, \\
\partial_t \ssy = \rhos\cp^2 \partial_y \vsy, & \xv\in\OmegaS, 
\end{cases} \label{eq:solidModelInviscid} \\
\text{Interface: } &\begin{cases}
\vfy = \vsy, & \xv\in\Gamma, \\
p = -\ssy, \quad & \xv\in\Gamma.
\end{cases}
\end{align}
\ese
We note that the horizontal component of the fluid velocity is omitted from the equations for the fluid in~\eqref{eq:fluidModelInviscid} since it decouples in the simplified FSI problem and can be determined from the remaining fluid variables once the solution is known.  The solid variables are assumed to be bounded as $y\rightarrow-\infty$.  The model problem described here corresponds to the MP-IA model problem considered in~\cite{fib2014}.

\input texFiles/rectangularModelProblemCartoon

\subsection{Discretization}

The equations for the fluid in~\eqref{eq:fluidModelInviscid} are transformed to Fourier space using~\eqref{eq:discreteFourier}.  The resulting equations for the Fourier coefficient functions are
\bse
\label{eq:fluidModelInv}
\begin{alignat}{2}
\rho \partial_t \vfyh + \partial_y \ph &= 0, \qquad &  0<y<H, \label{eq:fluidModelVelInv}\\
\partial_y^2 \ph - k_x^2 \ph &=0, &  0<y<H. \label{eq:fluidModelPresInv}
\end{alignat}
\ese
The solution of the elliptic equation for the pressure in~\eqref{eq:fluidModelPresInv}, with $p=0$ at $y=H$, is 
\begin{align}
p(y,t) = p_I(t) \frac{\sinh(k_x(H-y))}{\sinh(k_x H)},
\label{eq:pressureSubProblemSolution}
\end{align}
where $p_I(t)$ is an interface pressure which is specified later by the chosen coupling at the interface.  For the purpose of the analysis, we only require the fluid velocity on the interface, which we define to be $v_I(t)$.  The evolution of the interface velocity involves the gradient of pressure at $y=0$, which is given by
\begin{align}
\partial_y p(0,t) = - \frac{p_I(t)}{\Heff},\qquad \Heff = \frac{\tanh(k_x H)}{k_x}.
\label{eq:interfacePressureGradient}
\end{align}
The effective fluid height, $\Heff$, in the formula for the pressure gradient depends on the normalized wave number, and it takes the limiting values of $\Heff=H$ as $\kx\rightarrow0$ and $\Heff=0$ as $\kx\rightarrow\infty$.

The equations for the solid in~\eqref{eq:solidModelInviscid} are also transformed to Fourier space, and then expressed in the characteristic form
\bse
\begin{alignat}{2}
\partial_t \bc - \cp \partial_y \bc &= i k_x \cp \dc, \qquad && y<0,
\label{eq:advSrcB} \\
\partial_t \ac + \cp \partial_y \ac &= -i k_x \cp \dc, \qquad && y<0, 
\label{eq:advSrcA} \\
\partial_t \dc &= \frac{i k_x \cp}{2} (\bc - \ac), \qquad && y<0,
\label{eq:SrcD}
\end{alignat}
\ese
where $a(y,t)$, $b(y,t)$ and $d(y,t)$ are Fourier coefficients of the characteristic variables
\begin{align}
\bc = \ssyh + \zp\vsyh, \qquad
\ac = \ssyh - \zp\vsyh, \qquad
\dc = \ssxh.
\end{align}
Equations~\eqref{eq:advSrcB} and~\eqref{eq:advSrcA} are advection equations for the incoming and outgoing characteristic variables $\bc$ and $\ac$, respectively, while~\eqref{eq:SrcD} describes the evolution of the tangential component of the solid stress.

The solid characteristic variables are discretized in the $y$-direction using a uniform collocated grid, $y_j = j\dy$ for $j = 0, -1, -2 \ldots$, with grid spacing $\dy$.  Although a staggered grid is used for the analysis in~\cite{fib2014}, and for the stability analyses of other FSI algorithms in~\cite{BanksSjogreen2011,lrb2013}, it is cleaner to use a collocated grid for the second-order accurate scheme examined here.  Define the grid functions $\ac^n_j \approx \ac(y_j, t^n)$, $\bc^n_j \approx \bc(y_j, t^n)$ and $\dc^n_j \approx \dc(y_j,t^n)$, where $t^n = n \Delta t$ for a (fixed) time step $\Delta t$.  The advection equations in~\eqref{eq:advSrcB} and~\eqref{eq:advSrcA} are approximated using a second-order accurate Lax-Wendroff-type scheme having the form
\bse
\label{eq:LaxWend}
\begin{align}
\bc^{n+1}_{j} &= \bc^{n}_j + 
\frac{\lambda_y}{2} \delta_0 \bc^{n}_j + 
\frac{\lambda_y^2}{2} \delta_+ \delta_- \bc^{n}_j 
+ i \lambda_x \dc^n_j
+ \frac{i \lambda_x \lambda_y}{4} \delta_0 \dc^{n}_j
- \frac{\lambda_x^2}{4} (\bc^n_j - \ac^n_j), \label{eq:bLax} \\
\ac^{n+1}_{j} &= \ac^{n}_j - 
\frac{\lambda_y}{2} \delta_0 \ac^{n}_j + 
\frac{\lambda_y^2}{2} \delta_+ \delta_- \ac^{n}_j 
- i \lambda_x \dc^n_j
+ \frac{i \lambda_x \lambda_y}{4} \delta_0 \dc^{n}_j
+ \frac{\lambda_x^2}{4} (\bc^n_j - \ac^n_j),  \label{eq:aLax}
\end{align}
\ese
where $\lambda_x = k_x \cp \dt$ and $\lambda_y = \cp \dt / \dy$.  The centered approximations of the spatial derivatives in~\eqref{eq:LaxWend} are defined in terms of the undivided difference operators $\delta_+ u^{n}_j = u^{n}_{j+1} - u^{n}_j$, $\delta_- u^{n}_j = u^{n}_j - u^{n}_{j-1}$ and $\delta_0 u^{n}_j = u^{n}_{j+1} - u^{n}_{j-1}$ for an arbitrary grid function $u^{n}_j$.  The evolution equation in~\eqref{eq:SrcD} is approximated using a second-order accurate BDF scheme of the form
\begin{align}
\dc^{n+1}_j &= \frac{4}{3}\dc^{n}_j - \frac{1}{3} \dc^{n-1}+ \frac{i \lambda_x}{3}
\left(\bc^{n+1}_j - \ac^{n+1}_j \right). \label{eq:dTrap}
\end{align}
For reference, the Fourier coefficients of the solid velocity and normal stress are related to the incoming and outgoing characteristic variables by
\begin{align}
{\bar v}_{2,j}\sp n = \frac{1}{2\zp} (\bc^n_j-\ac^n_j), \qquad
{\bar\sigma}_{22,j}\sp n = \frac{1}{2} (\bc^n_j+\ac^n_j),
\end{align}
and the boundedness condition at infinity implies that
\begin{align}
\vert \ac^n_j \vert^2 +
\vert \bc^n_j \vert^2 +
\vert \dc^n_j \vert^2 < \infty, \qquad \text{as } j \rightarrow -\infty.
\label{eq:bcInf}
\end{align}

\subsection{Interface coupling}

The evolution of the solid requires a boundary condition on the interface, $y=0$, corresponding to the incoming characteristic variable.  Similarly, the evolution of the fluid requires a boundary condition at $y=0$, which can be interpreted as an interface pressure needed to complete the solution of the elliptic problem for the fluid pressure in~\eqref{eq:pressureSubProblemSolution}.  These two boundary conditions come from the matching conditions at the interface involoving the vertical velocity and normal stress.  In terms of the discretization of the Fourier coefficents for the model problem, the coupling at the interface takes on different forms depending on the choice of the partitioned scheme.  Since the behavior of the fluid for the model problem is determined by the solution of an elliptic problem for the pressure and an evolution equation for the vertical velocity on the interface, the coupling of the fluid to the behavior of the solid can be reduced to a modified boundary condition for the solid depending on the partitioned scheme employed as we discuss below.

Let us first consider the AMP algorithm and its interface coupling scheme.  The general algorithm is described in Section~\ref{sec:algorithm} and we follow this basic description but with suitable modifications for the choices of the discretizations made for the present model problem.  Let us assume that the discrete characteristic variables for the solid given by $\bc^{n}_j$, $\ac^{n}_j$ and $\dc^{n}_j$ are known at time levels $t\sp n$ and $t\sp{n-1}$ for all grid points, $y_j=j\dy$, $j=0,-1,\ldots$, including a ghost point at $j=1$.  We also assume that the interface velocity, $v_I(t)$, is known at $t=t\sp n$ and $t\sp{n-1}$.  The evolution equations for the solid variables in~\eqref{eq:LaxWend} and~\eqref{eq:dTrap} are used to advance the variables to $t\sp{n+1}$, and from these we compute the solid velocity and normal stress on the interface using
\begin{align}
{\bar v}_{2,0}\sp{n+1} = \frac{1}{2\zp} (\bc^{n+1}_0-\ac^{n+1}_0), \qquad
{\bar\sigma}_{22,0}\sp{n+1} = \frac{1}{2} (\bc^{n+1}_0+\ac^{n+1}_0).
\label{eq:interfaceSolidVelocityAndStress}
\end{align}
We next compute an extrapolated interface velocity, $v_I\sp{(e)}$, using a second-order accurate BDF-type integration of the evolution equation for the fluid velocity.  This gives the formula
\begin{align}
v_I^{(e)} = \frac{4}{3} v_I^n - \frac{1}{3} v_I^{n-1} - \frac{2\dt}{3 \rho}\partial_y p(0,t^{n+1}),
\label{eq:extrapInterfaceVelocity}
\end{align}
where the pressure gradient on the interface is given in~\eqref{eq:interfacePressureGradient}.  This extrapolated value for the fluid is averaged with the interface velocity for the solid from~\eqref{eq:interfaceSolidVelocityAndStress} to give
\begin{align}
v_I^{n+1} = \frac{\zf}{\zf + \zp}\,v_I^{(e)} + \frac{\zp}{\zf + \zp}\,{\bar v}_{2,0}\sp{n+1},
\label{eq:interfaceVelocityModel}
\end{align}
where we take $\zf = \rho H/\dt$.  The choice for $\zf$ here is the same one used in~\cite{fib2014}, and it agrees with~\eqref{eq:fluidImpedanceAlg} for $\Cam h=H$.  We have found that the results are insensitive to the choice of $\zf$ for this inviscid model problem.  For later convenience, define the weights
\begin{align}
\theta_f=\frac{\zf}{\zf + \zp},\qquad {\bar\theta}_p=\frac{\zp}{\zf + \zp},
\label{eq:weights}
\end{align}
and we note that $\theta_f+{\bar\theta}_p=1$.  The AMP condition for the fluid pressure 
given in~\eqref{eq:AMPpressureBC} requires an acceleration of the interface, and for this we use the backward difference formula
\begin{align}
\dot{v}_I^{n+1} = \frac{1}{2\Delta t} \left(3v_{I}^{n+1} - 4v_{I}^{n} + v_{I}^{n-1}  \right).
\label{eq:interfaceAcceleration}
\end{align}
Using~\eqref{eq:extrapInterfaceVelocity} in~\eqref{eq:interfaceVelocityModel} to eliminate $v_I\sp{(e)}$, and then using the computed result in~\eqref{eq:interfaceAcceleration} to eliminate $v_{I}^{n+1}$ gives an alternate expression for the acceleration having the form
\begin{align}
\dot{v}_I^{n+1} = \theta_f\left(\frac{p_I^{n+1}}{\rho \Heff}\right) + {\bar\theta}_p\left(\frac{3{\bar v}_{2,0}\sp{n+1} - 4v_{I}^{n} + v_{I}^{n-1}}{2\Delta t} \right),
\label{eq:interfaceAccelerationNew}
\end{align}
which shows that the acceleration of the interface is also a weighted average of the acceleration from the fluid and that computed from the solid.  
Finally, the AMP condition for the fluid pressure in~\eqref{eq:AMPpressureBC} reduces to
\begin{align}
-p_I^{n+1}-\frac{\zp \dt}{\rho\Heff}\, p_I^{n+1} = {\bar\sigma}_{22,0}^{n+1} - \zp \dt \,\dot{v}_I^{n+1}.
\label{eq:AMPpressureBCmodel}
\end{align}
Using~\eqref{eq:interfaceAccelerationNew} in~\eqref{eq:AMPpressureBCmodel} to eliminate $\dot{v}_I^{n+1}$, and then solving for the interface pressure gives
\begin{align}
-p_I^{n+1} = \frac{M}{M+{\bar\theta}_p}\left[{\bar\sigma}_{22,0}^{n+1} - \zp{\bar\theta}_p\left(\frac{3}{2}{\bar v}_{2,0}\sp{n+1} - 2v_{I}^{n} + \frac{1}{2}v_{I}^{n-1} \right)\right],
\label{eq:interfacePressureModel}
\end{align}
where $M$ is a mass ratio given by
\begin{align}
M = \frac{\rho\Heff}{\zp \dt} = \frac{\rho\Heff}{\rhos \cp \dt}.
\label{eq:AcousticMassRatio}
\end{align}
We may now complete the time-step for the solid variables by specifying values for the incoming and outgoing characteristic variables in the ghost points using
\begin{align}
\bc^{n+1}_1 = -\bc^{n+1}_{-1} + 2 a_I\sp{n+1},\qquad \ac^{n+1}_1 = 2 \ac^{n+1}_0 - \ac^{n+1}_{-1},
\label{eq:interfaceBCsForSolid}
\end{align}
where $a_I\sp{n+1}$ is an interface value for the incoming characteristic variable determined by
\begin{align}
a_I\sp{n+1} = -p_I^{n+1} + \zp v_I^{n+1},
\label{eq:incomingCharacteristicInterfaceModel}
\end{align}
Thus, the boundary condition for the incoming characteristic variable is specified by $a_I\sp{n+1}$ (using a second-order accurate average) and the outgoing characteristic variable at the ghost point is given by a simple extrapolation.

For the AMP algorithm, the interface velocity and pressure used in~\eqref{eq:incomingCharacteristicInterfaceModel} to define the interface value for the incoming characteristic variable are given by the formulas in~\eqref{eq:interfaceVelocityModel} and~\eqref{eq:interfacePressureModel}, respectively, while the corresponding interface values for the TP and ATP algorithms use difference formulas.  For the TP algorithm, the interface velocity is given by the solid, while the interface stress is determined by the fluid.  This specification leads to
\begin{align}
v_I\sp{n+1}={\bar v}_{2,0}\sp{n+1},\qquad p_I\sp{n+1}=\rho\Heff\, \dot{v}_I^{n+1},
\label{eq:interfaceValuesTPscheme}
\end{align}
where the solid velocity is given in~\eqref{eq:interfaceSolidVelocityAndStress} and the interface acceleration is given in~\eqref{eq:interfaceAcceleration}.  For the ATP algorithm, the interface velocity is given by the fluid, while the interface stress is determined by the solid, and this gives
\begin{align}
v_I\sp{n+1}=v_I\sp{(e)},\qquad p_I\sp{n+1}=-{\bar\sigma}_{22,0}^{n+1},
\label{eq:interfaceValuesATPscheme}
\end{align}
where the fluid velocity is given in~\eqref{eq:extrapInterfaceVelocity} and the solid stress is given in~\eqref{eq:interfaceSolidVelocityAndStress}.  It is worth noting that values specified for the TP and ATP algorithms in~\eqref{eq:interfaceValuesTPscheme} and~\eqref{eq:interfaceValuesATPscheme}, respectively, agree with values given by the AMP algorithm in the limiting cases of a heavy solid ($\zp\rightarrow\infty$) and a light solid ($\zp\rightarrow0$).

\subsection{Stability analysis}

\newcommand{\nn}{\alpha}

In order to assess the stability of the AMP, TP and ATP schemes, we consider normal mode solutions of the discrete equations for the solid with the appropriate interface conditions.  For this analysis, we let
\begin{align}
\bc^{n}_j = \amp^n \tbc_j, \qquad 
\ac^{n}_j = \amp^n \tac_j, \qquad 
\dc^{n}_j = \amp^n \tdc_j,
\end{align}
where $\amp$ is an amplification factor and $(\tbc_j,\tac_j,\tdc_j)$ are spatial grid functions. 
As in stability analysis of Section~\ref{sec:viscousAnalysis}, showing the scheme is weakly stable 
  is equivalent to showing there are no non-trivial solutions where $|\amp| > 1$
  for a given set of parameters.
  We first make the assumption $|A| > 1$ and solve the spatial problem by finding 
  general solutions of the discretization and boundedness condition as $j\rightarrow-\infty.$
  The interface conditions are then used to determine if the amplification factor satisfies the
  assumption $|A|>1.$

The discrete equations for the solid in~\eqref{eq:LaxWend} and~\eqref{eq:dTrap}, with $\tacv_j = [\tbc_j,\tac_j,\tdc_j]^T$, can be written as a system of linear, second-order difference equations of the form
\begin{align}
[H_0 + H_1 \delta_0 + H_2 \delta_+ \delta_-] \tacv_j = 0,
\label{eq:DifferenceEquation}
\end{align}
where
\[
\begin{array}{c}
\displaystyle{
H_{0} = 
-A^2 
\left[\begin{array}{ccc}
1&0&0 \medskip\\
0&1&0 \medskip\\
-\frac{i \lx}{3}&\frac{i \lx}{3}&1
\end{array}\right] 
+A 
\left[\begin{array}{ccc}
1- \frac{\lx^2}{4}&\frac{\lx^2}{4}&i \lx \smallskip\\
\frac{\lx^2}{4}&1-\frac{\lx^2}{4}&-i \lx \smallskip\\
0&0&\frac{4}{3}
\end{array}\right] 
-
\left[\begin{array}{ccc}
0&0&0 \medskip\\
0&0&0 \medskip\\
0&0&\frac{1}{3}
\end{array}\right],
} \bigskip\\
\displaystyle{
H_{1} = 
A \left[\begin{array}{ccc}
\frac{\ly}{2} &0&\frac{i \lx \ly}{4} \\
0&-\frac{\ly}{2}&\frac{i \lx \ly}{4} \\
0 & 0&0 
\end{array}\right], \qquad
H_{2} = 
A \left[\begin{array}{ccc}
\frac{\ly^2}{2} &0&0 \\
0&\frac{\ly^2}{2}&0 \\
0 & 0&0 
\end{array}\right].
}
\end{array}
\]
Solutions of the difference equations in~\eqref{eq:DifferenceEquation} are sought in the form $\tacv_j = \tilde\rev \p^j$, where $\tilde\rev$ is a constant vector and $\p$ is a scalar.  This implies the homogeneous linear system
\begin{align}
\left[ H_0 + H_1 \left(\p - \p\sp{-1}\right) + H_2 \left(\p - 2 + \p\sp{-1} \right) \right] \tilde\rev = 0.
\label{eq:evProblem}
\end{align}
Non-trivial vectors $\tilde\rev$ exist for values of $\p$ satisfying
\begin{align}
\det\left[ H_0 + H_1 \left(\p - \p\sp{-1}\right) + H_2 \left(\p - 2 + \p\sp{-1} \right) \right] = 0.
\label{eq:quarticPhi}
\end{align}
It can be shown that the algebraic condition in~\eqref{eq:quarticPhi} is equivalent to a quartic polynomial in $\p$ for which there are four roots, $\p_\nn$, $\nn=1,2,3,4$.  Since we seek spatial grid functions that are bounded as $j\rightarrow-\infty$, we are only interested in the roots satisfying $\vert\p_\nn\vert\ge1$.

\medskip\noindent
\textbf{Lemma}: \textit{If $|A| > 1$ and if $\lx$ and $\ly$ are chosen to satisfy a CFL condition, then the roots of the quartic polynomial implied by~\eqref{eq:quarticPhi} satisfy}
\begin{align}
|\phi_1| > 1, \qquad |\phi_2| > 1, \qquad |\phi_3| < 1, \qquad |\phi_4| < 1.
\label{eq:SizeOfEigenvalues}
\end{align}

\medskip\noindent
The analogous lemma for the viscous stability analysis 
  was proven in Section~\ref{sec:viscousAnalysis}.
  The following proof is structured similarly. 
  First, the stability for the
  scheme applied to the pure initial-value problem (Cauchy problem) is considered.
  This corresponds to setting $\phi=e^{i\vartheta}$ in~\eqref{eq:quarticPhi} and 
  finding the region of stability $|A|\le 1$ in the $\lx$-$\ly$ plane for all 
  $\vartheta \in [0,2\pi].$
  This region is found in Section~\ref{sec:CFLRegionCauchy} and shown in Figure~\ref{fig:inviscidCFL}.
  This region contains the quarter circle defined by $\lx^2+\ly^2\le 1$ and $\lx,\ly\ge 0.$
%
When $|\phi|=1,$ we must have $|A| \le 1.$ Therefore $|A|>1$ implies $|\phi| \ne 1.$
  Due to the continuity of $\phi,$ if $|A| > 1$ and if $(\lx,\ly)$ remains within 
  the CFL restriction, $\phi$ cannot cross the unit circle, $|\phi|=1$, as $(\lx,\ly)$ vary. 
Therefore, it is only necessary to prove that the conditions on the roots 
in~\eqref{eq:SizeOfEigenvalues} hold for one set of parameters.
For $\lx=0$, the difference scheme reduces to a pair of decoupled Lax-Wendroff schemes (and a BDF scheme applied to $\partial_td=0$), where it is well known (see~\cite{GKSII} for example) that the conditions in~\eqref{eq:SizeOfEigenvalues} hold provided $\vert\ly\vert\le1$.  Thus, the condition holds for all $(\lx,\ly)$ provided the CFL condition is satisfied.

In view of the above Lemma, we can now write the general solution for $\acv^n_j=[\bc^{n}_j,\ac^{n}_j,\dc^{n}_j]\sp T$, satisfying the regularity condition and assuming $\vert A\vert>1$, in the form
\begin{align}
\acv^n_j = \amp^n \left(k_1 \tilde\rev_1 \p_1^j + 
k_2 \tilde\rev_2 \p_2^j \right), \label{eq:EigenvectorExpansion}
\end{align}
where $k_1$ and $k_2$ are scalar constants.  Let $\tilde\rev_\nn = [\qn{\nn},\rn{\nn},\sn{\nn}]^T$ for $\nn = 1,2$.  We note that the third equation in the matrix system~\eqref{eq:evProblem} implies
\begin{align}
\sn{\nn} = \frac{A^2 i \lx}{3 A^2- 4 A + 1} (\qn{\nn} - \rn{\nn}),\qquad \nn = 1,2,
\end{align}
so that the third component of $\acv^n_j$ can be eliminated in terms of the first two components.  The final two constraints are given by the interface conditions in~\eqref{eq:interfaceBCsForSolid}, and this leads to a homogeneous system for $\kv=[k_1,k_2]\sp T$ of the form
\begin{align}
\dMatI(A) \kv= \left[\begin{array}{cc}
\dComp_{11} & \dComp_{12} \\
\dComp_{21} & \dComp_{22}
\end{array}\right] 
\left[\begin{array}{c}
k_1 \\ k_2
\end{array}\right] = 0.
\label{eq:StabilitySystem}
\end{align}
The coefficients of the matrix $\dMatI$ involve the amplitude factor $\amp$, and these are given in Section~\ref{sec:ComponentsOfDmat} for the AMP, TP and ATP schemes. For non-trivial solutions, we seek solutions for $\amp$ such that
\begin{align}
\det(\dMatI) = \dComp_{11} \dComp_{22} - \dComp_{12} \dComp_{21} = 0.
\label{eq:determinantCondition}
\end{align}
We first proceed with an analysis of the 1D case assuming $\lx = 0$.  After finding solutions for~$\amp$ when $\lx = 0$, we use a continuation procedure to determine solutions for $\amp$ numerically as $\lx$ increases from zero.


\subsubsection{1D schemes ($\lx=0$)} \label{sec:inviscidAnalysis1D}

As noted earlier, the discrete evolution equations in~\eqref{eq:LaxWend} for the solid reduce to a pair of standard Lax-Wendroff schemes (without source terms) for the variables associated with incoming and outgoing characteristics when $\lx=0$.  Accordingly, the first two components of~\eqref{eq:DifferenceEquation} become
\begin{align}
A \tbc_j = \left( 1 + \frac{\lambda_y}{2} \delta_0 
+ \frac{\lambda_y^2}{2} \delta_+ \delta_- \right) \tbc_j,\qquad
A \tac_j = \left( 1 - \frac{\lambda_y}{2} \delta_0 
+ \frac{\lambda_y^2}{2} \delta_+ \delta_- \right) \tac_j.
\end{align}
Setting $\tbc_j=\tbc_0\p_{\bc}\sp j$ and $\tac_j=\tac_0\p_{\ac}\sp j$ yields the eigenvalues
\begin{align}
\p_{\bc,\pm} =\frac{\amp-(1-\ly\sp2) \pm \sqrt {\amp\sp2-(2\amp-1)(1-\ly\sp2)}}{\ly \left( \ly+1 \right) }, \qquad 
\p_{\ac,\pm} =  \frac{1}{\p_{\bc,\mp}}.
\label{eq:eigenvalueOneD}
\end{align}
If $\vert\amp\vert>1$ and if $\ly$ satisfies the CFL constraint $\vert\ly\vert\le1$, then one eigenvalue from each plus-minus pair in~\eqref{eq:eigenvalueOneD} has magnitude greater than one, while the other has magnitude less than one in agreement with~\eqref{eq:SizeOfEigenvalues}.  Let $\p_1$ to be the eigenvalue from $\p_{\bc,\pm}$ with magnitude greater than one and let $\p_2$ to be the corresponding eigenvalue $\p_{\ac,\pm}$.  Assuming the principal branch of the square-root function in~\eqref{eq:eigenvalueOneD}, we set
\begin{align}
(\p_1,\p_2)=
\begin{cases}
(\p_{\bc,+},\p_{\ac,+}), & \hbox{if $\Re(A) \ge 1 - \ly^2$}, \\
(\p_{\bc,-},\p_{\ac,-}), & \hbox{if $\Re(A) < 1 - \ly^2$}.
\end{cases}
\end{align}
%
%
Since the evolution equations for $\bc_j\sp n$ and $\ac_j\sp n$ are decoupled, the corresponding eigenvectors in~\eqref{eq:EigenvectorExpansion} are unit vectors, $\tilde\rev_{1} = [1,0,0]^{T}$ and $\tilde\rev_{2} = [0,1,0]^{T}$, i.e.~$\qn{2} = \rn{1} = \sn{1} = \sn{2} = 0$, and the matrix $\dMatI$ in~\eqref{eq:StabilitySystem} becomes diagonal. Thus, the determinant condition in~\eqref{eq:determinantCondition} reduces to
\begin{align}
\det(\dMatI)= \dComp_{11} \dComp_{22} = 0.
\end{align}
The component, $\dComp_{22}$, is associated with the extrapolation of the outgoing characteristic variable from~\eqref{eq:interfaceBCsForSolid}, and this factor is zero only if $\p_2=1$.  However, $\vert\p_2\vert>1$, so that $\dComp_{22}\ne0$ and thus $k_2=0$ in~\eqref{eq:EigenvectorExpansion}.  The remaining component, $\dComp_{11}$, is associated with the interface condition on the incoming characteristic variable, and this term takes different forms depending on the choice of the interface coupling (see Section~\ref{sec:ComponentsOfDmat}).  Manipulation of the constraint, $\dComp_{11}=0$, for the AMP, TP and ATP schemes leads to polynomials in $\amp$ of degrees~6, 6 and~4, respectively, whose coefficients depend on $\ly$ and $M_0=\rho H/(\zp \dt)$.
  Due to the algebraic complexity, the software package Maple is used to generate the polynomials and compute the roots for given values of $(\ly,M_0)$.  We note that some roots of the polynomials may not be solutions of the determinant constraint and so these spurious roots are discarded.  If no valid roots are found satisfying $\vert\amp\vert>1$ for a given value of $(\ly,M_0)$, then we conclude that $k_1$ in~\eqref{eq:EigenvectorExpansion} must also be zero so that only the trivial solution exists and thus the scheme is (weakly) stable.

We find it convenient to explore the roots of the three polynomials in terms of the parameters $\ly$ and $\Mr$, where
\begin{align}
\Mr = M_0\ly= \frac{\rho H}{\rhos \dy}.
\label{eq:massRatio}
\end{align}
The latter parameter is the ratio of the mass contained in a fluid column of height $H$ to that of the solid with the same cross section and of height $\dy$.  Starting with the AMP scheme, we compute valid roots of its polynomial for an $800\times800$ array of parameter values on a logarithmic scale for the range $10^{-6} \le \lambda_y \le 1$ and $10^{-6} \le \Mr \le 10^{7}$.  For these values, no roots with $\vert\amp\vert>1$ are found, which provides a strong indication of stability of the 1D AMP scheme for $0<\ly\le1$ and $\Mr>0$.  Shaded contours of the maximum values of $\vert\amp\vert$ are displayed in the left plot of Figure~\ref{fig:acousticSolidStokesStabilityTP} for a range of values of $\Mr$ and $\ly$.  The roots of the polynomials associated with the 1D TP and ATP schemes are also computed, and parameter values for which $\vert\amp\vert>1$ are found indicating that regions of instability exist for these two schemes.  The right plot of Figure~\ref{fig:acousticSolidStokesStabilityTP} shows regions of stability of the TP and ATP schemes where no roots with $\vert\amp\vert>1$ are found.  We note that for a fixed value of $\ly=\cp \dt / \dy\in(0,1]$, decreasing the mesh spacings in both time and space corresponds to increasing the mass ratio $\Mr$ defined in~\eqref{eq:massRatio}.  The right plot in the figure indicates that the TP scheme becomes unstable as the mesh is refined, while the ATP scheme ultimately becomes stable.  In contrast, the AMP scheme is always stable, and these observations agree with the analysis of the first-order accurate difference schemes performed in~\cite{fib2014}.

{
\newcommand{\figWidth}{6cm}
\begin{figure}[h]
\begin{center}
  \includegraphics[width=\figWidth]{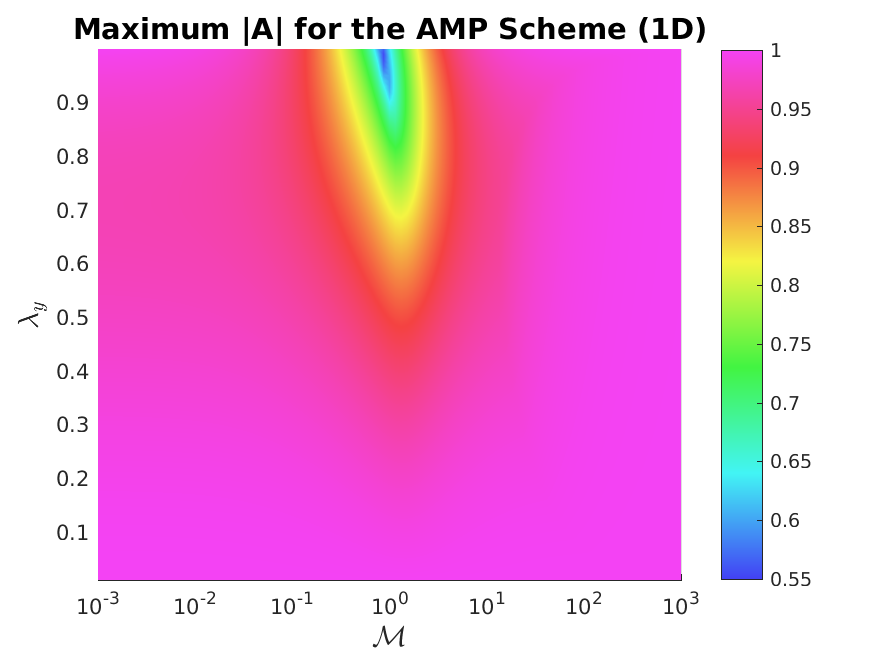}
  \includegraphics[width=\figWidth]{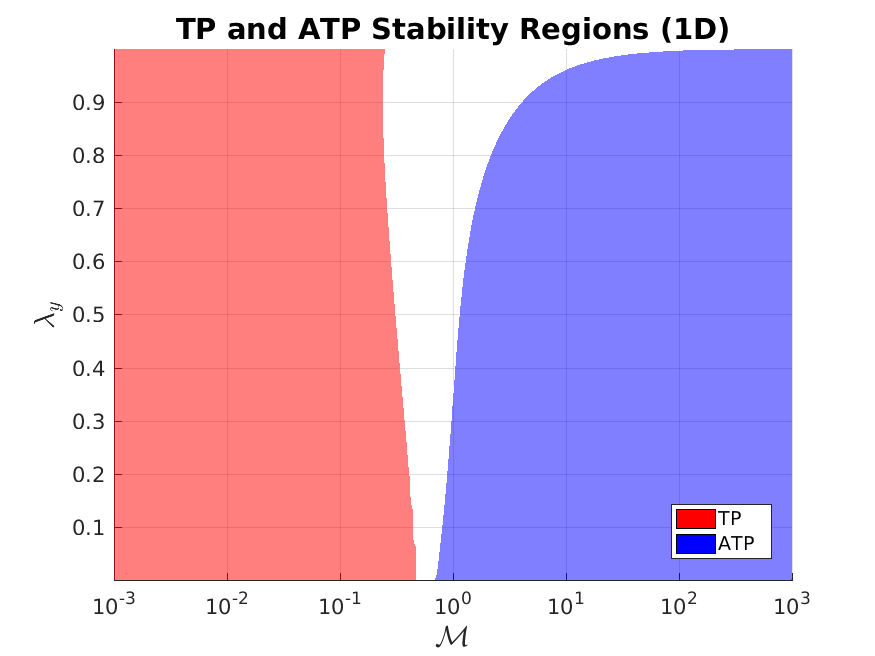}
  \caption{
    Left: Shaded contours of the maximum amplification factor $\vert\amp\vert$ 
    for the AMP scheme for the parameters $\lambda_y = \cp \dt / \dy$ and $\Mr = \rho H / (\rhos \dy)$.
    Right: Stability regions for the 1D TP and ATP schemes. 
    Grid refinement corresponds to increasing $\Mr$, and thus the
    TP scheme becomes unstable as the grid is refined while the ATP scheme becomes stable.
  }  \label{fig:acousticSolidStokesStabilityTP}
\end{center}
\end{figure}
}


The AMP, TP and ATP schemes can be implemented for the one-dimensional FSI model problem, and then run numerically for different values of the parameters to check the results of the stability analysis.  
Initial conditions are specified so that the exact solution decays 
exponentially in time, see Section~\ref{sec:SolutionToOneDimensional} for details.
Since the interface velocity is expected to decay in time, we can assess the stability of the scheme by checking the 
magnitude of the interface velocity.  If this velocity grows in time, becoming larger than the initial interface velocity, then we consider the scheme to be unstable for the choice of parameter values.
For the AMP scheme, no parameter values are found for which the numerical results are unstable.  For the TP and ATP schemes, the colored marks in Figure~\ref{fig:acousticSolidStokesStabilityTPNumerical} indicate stability or instability of the numerical scheme for selected parameter values.  The black curves in the figure correspond to the stability boundaries obtained from the stability analysis.  We note that there is good agreement between the numerical results of the schemes and that given by the analysis.


{
\newcommand{\figWidth}{6cm}
\begin{figure}[h]
\begin{center}
  \includegraphics[width=\figWidth]{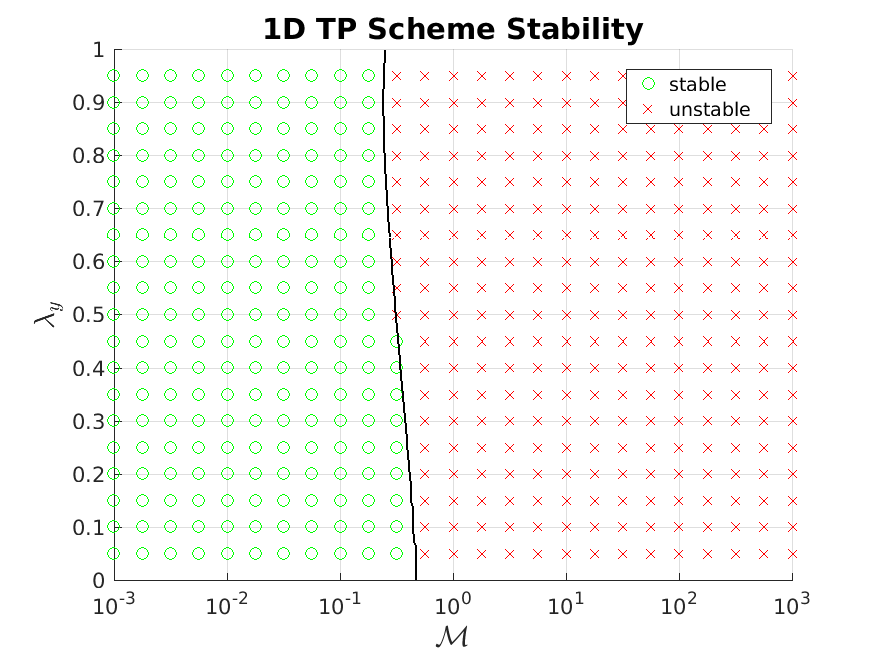}
  \includegraphics[width=\figWidth]{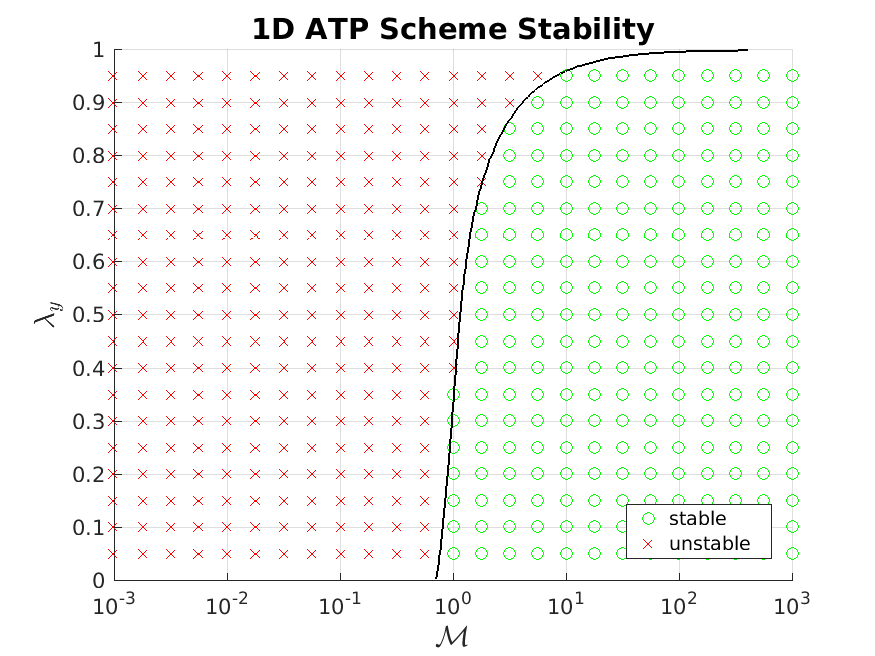}
  \caption{Verification of the stability regions for the TP scheme (left)
    and the ATP scheme (right) using numerical solutions of the 1D FSI model problem. Green marks 
    correspond to a stable numerical solution, while red markers correspond to numerical instability.  The black curves give the stability boundaries obtained from the stability analysis.
  }  \label{fig:acousticSolidStokesStabilityTPNumerical}
\end{center}
\end{figure}
}


%

\subsubsection{2D schemes ($\vert\lx\vert\ne0$)}

\newcommand{\zz}{{z}}  
\newcommand{\NN}{{\mathcal{Z}}}  

The strategy for finding solutions of the determinant condition in~\eqref{eq:determinantCondition} for $\vert\lx\vert \ne 0$ is to use the solutions for $\amp$ of the 1D determinant condition as a starting point for a numerical root solver that performs a continuation as $\lx$ varies.  In addition to the principal continuation parameter $\lx$, the roots also depend on the parameters $\ly$ and $M=\rho\Heff/(\zp\dt)$, where $\Heff$ was defined previously in~\eqref{eq:interfacePressureGradient}.  As before, we find it convenient to examine the roots in terms of the mass ratio $\Mr$ introduced in~\eqref{eq:massRatio} and the height ratio, $\eta=\Heff/H$, which varies on $(0,1]$, and we note that $M=\Mr\eta$.  Using our continuation procedure, we are able to compute the roots of the determinant condition as a function of the parameters $\lx$, $\ly$ and $M=\Mr\eta$.  Though there are multiple roots, $\amp_\zz$, $\zz=1,2,\ldots,\NN$, for a given set of parameters, we are only interested in the amplification factor with the maximum modulus.  Define
\begin{align}
A_{\max}(\lx,\,\ly,\,\Mr \eta)=\max_{1\le\zz\le\NN}\vert\amp_\zz(\lx,\,\ly,\,\Mr \eta)\vert,
\end{align}
to be the maximum amplification factor in magnitude for a given set of parameters.  The pseudo-spectral approximation of the model problem given in~\eqref{eq:solidModelInviscid} depends on the normalized wave number $k_x$.  The corresponding discrete evolution equations for each Fourier mode of the approximation must all be stable for the chosen partitioned scheme to be stable (either the AMP, TP or ATP schemes).  Since $\lx$ and $\eta$ both depend on $k_x$, we are interested in the values of $A_{\max}$ for all possible values for $\lx$ and $\eta$.  Define
\begin{align}
\mathcal{A}_{\max}(\Mr,\ly) = 
\max_{\lx \in \mathcal{D}(\ly),\, \eta \in (0,1]} A_{\max}(\lx,\,\ly,\,\Mr \eta),
\label{eq:AmaxDefinition2D}
\end{align}
where $\mathcal{D}(\ly)$ gives the interval of $\lx$ satisfying the Cauchy stability bound for a given value of $\ly$ (see Section~\ref{sec:CFLRegionCauchy}).  We take $\mathcal{D}(\ly) = [0, (1 - \ly^2)^{1/2}]$.
Since the maximum amplification factor for the 1D case (with $\lx=0$) must be less than or equal to $\mathcal{A}_{\max}$, we are effectively examining whether the maximum amplification can increase when $\vert\lx\vert>0$. Figure~\ref{fig:lysmAMP2D} shows shaded contours of $\mathcal{A}_{\max}(\Mr,\ly)$ on the left, and the stability regions of the TP and ATP schemes on the right.  These plots are similar to the ones for the 1D analysis in Figure~\ref{fig:acousticSolidStokesStabilityTP}, indicating that the stability results provided by the 1D analysis are essentially unchanged when variations in the transfer direction are considered.  In fact, Figure~\ref{fig:lysmAMP2Doverlay} shows the stability regions for the TP and ATP schemes given by the 1D and 2D analyses overlayed.  Here, we observe that the stability regions given by the 2D analysis are only slightly smaller than those given by the 1D analysis.

We can also search for a maximum amplification factor over all mass ratios,
$\mathcal{A}_{\text{CFL}},$ to see the CFL region for the AMP scheme in $(\lx,\ly)$ plane. 
Define
\begin{align}
\mathcal{A}_{\text{CFL}}(\lx,\ly) = 
\max_{0 < \Mr \eta < \infty} A_{\max}(\lx,\,\ly,\,\Mr \eta) .
\end{align}
This definition allows us to find a time step restriction regardless of the
mass ratio of the problem.
Figure~\ref{fig:lylxAMP2D} shows a surface plot of $\mathcal{A}_{\text{CFL}}$
and the stability region, $|A| \le 1$ in the $\lx$-$\ly$ plane.

\medskip\noindent
\textbf{Theorem}: \textit{By numerical evaluation of $\ACFL(\lx,\ly)$,
we have found that the AMP algorithm is weakly stable, i.e.~$|\ACFL| \le 1$,
provided that $\lx^2 + \ly^2 \le 1$.  This implies the time-step restriction}
\begin{align}
\dt \le \frac{1}{\cp} \left[\frac{1}{\dy^2} + k_x^2 \right]^{-1/2}.
\end{align}
\textit{This is a sufficient but not a necessary condition.}

%
{
\newcommand{\figWidth}{6cm}
\begin{figure}[h]
\begin{center}
  \includegraphics[width=\figWidth]{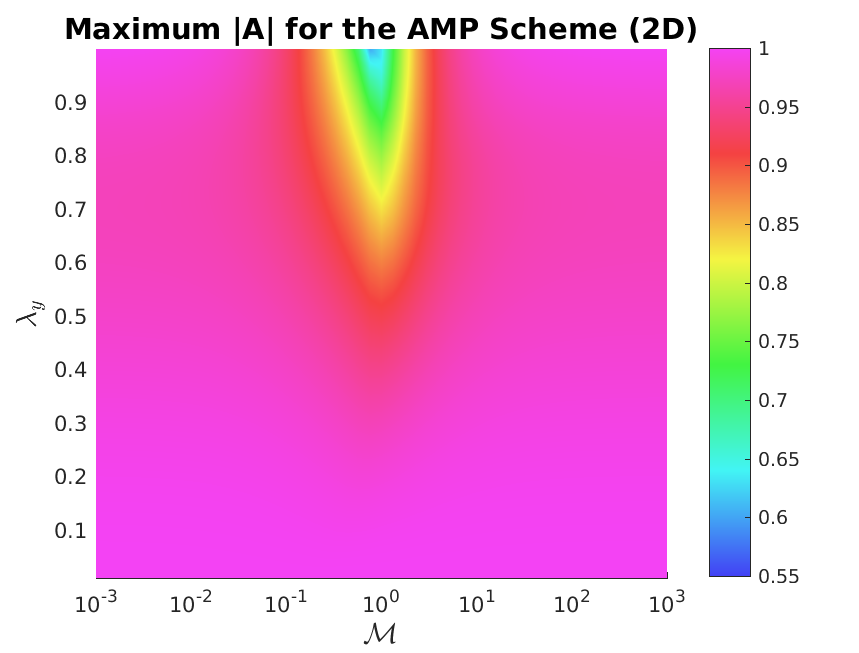}
  \includegraphics[width=\figWidth]{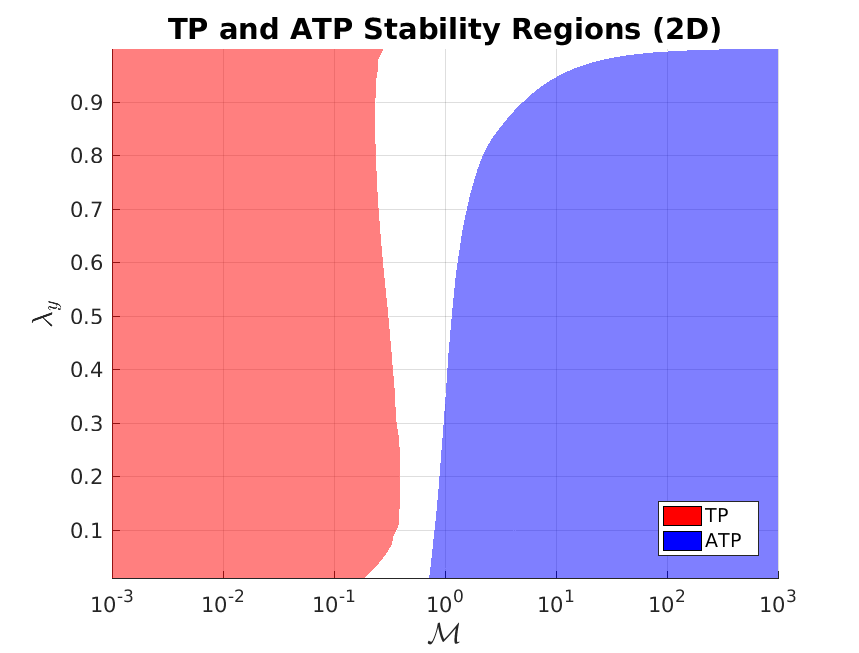}
  \caption{Left: Shaded contours of the maximum amplification factor 
    $\mathcal{A}_{\max}(\Mr,\ly)$ 
    (as defined in equation \eqref{eq:AmaxDefinition2D})
    for the AMP scheme for the parameters $\ly=\cp \dt/\dy$ and $\Mr=\rho H/ (\rhos \dy).$
    Right: Stability regions (ie. where $\mathcal{A}_{\max}(\Mr,\ly) \le 1$)
    for the 2D TP and ATP schemes. 
    Grid refinement corresponds 
    to increasing $\Mr$, and thus the
    TP scheme becomes unstable as the grid is refined while the ATP scheme becomes stable.
  } \label{fig:lysmAMP2D}
\end{center}
\end{figure}
}
{
\newcommand{\figWidth}{6cm}
\begin{figure}[h]
\begin{center}
  \includegraphics[width=\figWidth]{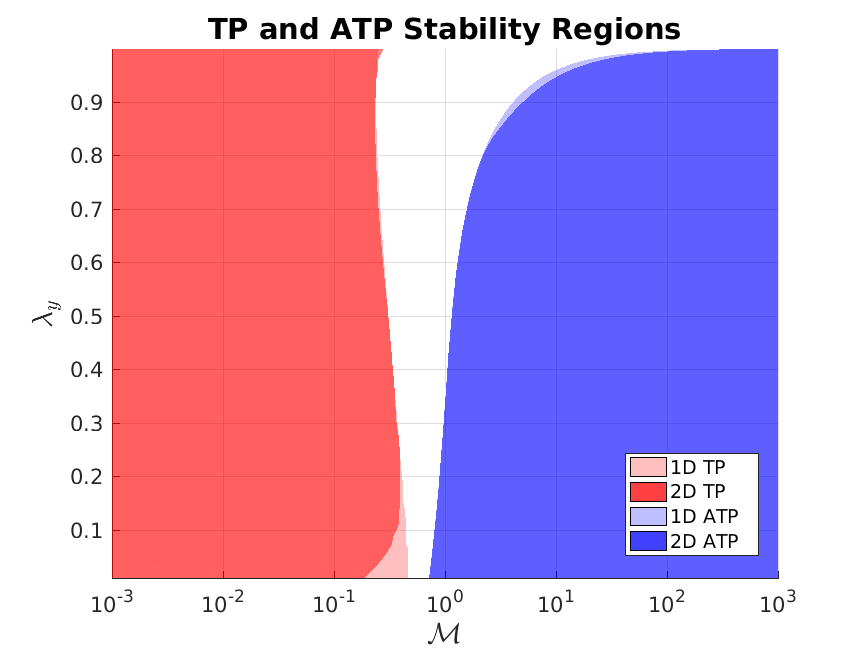}
  \caption{Stability regions for the TP and ATP schemes in 1D and 2D for the 
    parameters $\ly=\cp \dt/\dy$ and $\Mr=\rho H/ (\rhos \dy).$
    The lighter shaded colors represent the 1D regions and the darker shaded colors represent
    the 2D regions.
    Grid refinement corresponds to increasing $\Mr$, and thus the
    TP scheme becomes unstable as the grid is refined while the ATP scheme becomes stable.
  } \label{fig:lysmAMP2Doverlay}
\end{center}
\end{figure}
}
{
\newcommand{\figWidth}{6cm}
\begin{figure}[h]
\begin{center}
  \includegraphics[width=\figWidth]{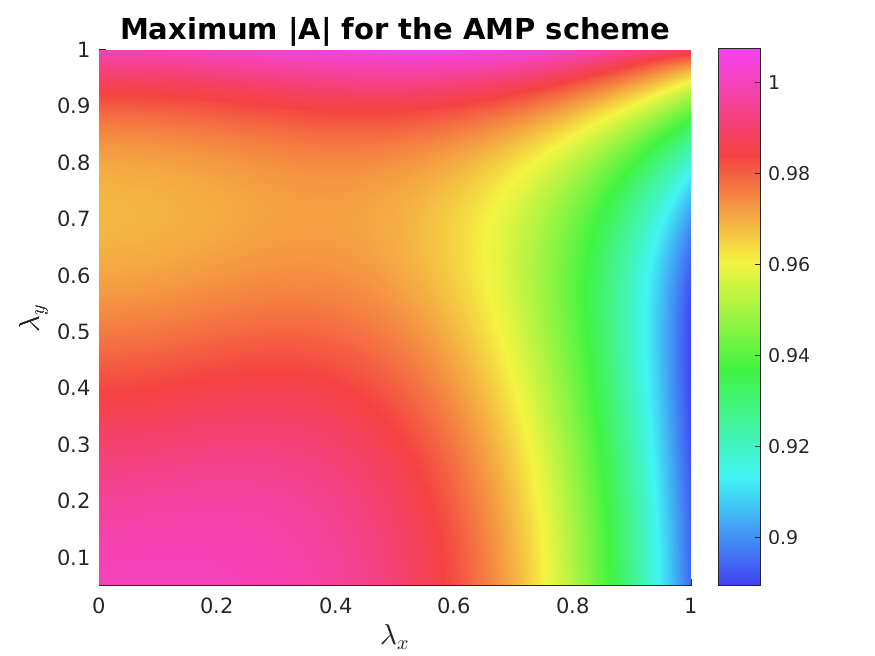}
  \includegraphics[width=\figWidth]{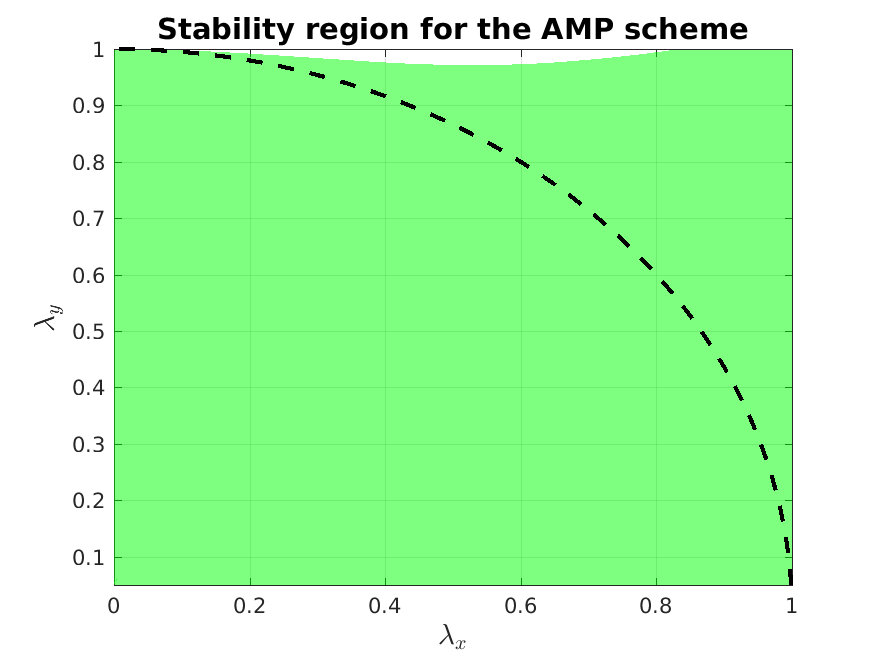}
  \caption{
    Left: surface plot of maximum $|A|$ over $\Mr \eta \in [10^{-6},10^{6}]$
    for the AMP scheme.
    Right: the green fill indicates the stability region $|A| \le 1.$
    The dotted line is a plot of
    $\lx^2 + \ly^2 = 1$. 
    For reference, 
    $\ly = \cp \dt / \dy,\, \Mr = \rho H / (\rhos \dy),\, \lx = \cp k_x \dt,\, 
    \eta = \tanh(k_x H)/ (k_x H).$
  } \label{fig:lylxAMP2D}
\end{center}
\end{figure}
}

%% file: texFiles/rectangularModelProblemCartoon.tex
{
\newcommand{\lbfont}{\small}
\def\ysb{-.15} 
\def\ysa{-2.5}   
\def\ya{0} 
\def\yb{2.5} 
\def\xL{8}
\begin{figure}[hbt]
	\newcommand{\textFont}{\normalss}
	\begin{center}
		\begin{tikzpicture}[scale=.9]
		\useasboundingbox (0,\ysa) rectangle (\xL,3.2);  
		%
                \draw[fill=red!20,draw=red!20] (0,\ysa) rectangle (\xL,\ysb);
                \draw[-,thick,red,line width=2pt] 
                (0  ,\ysa) -- 
                (0  ,\ysb) -- 
                (\xL,\ysb) --
                (\xL,\ysa);
		\draw[thick,black] (4,-1.2) node {solid: $\OmegaS$};
		\draw[thick,fill=blue!20,draw=blue,line width=2pt] (0,\ya) rectangle (\xL,\yb);
		\draw[thick,blue] (0,\ya) rectangle (\xL,\yb);
                \draw[thick,black] (.5*\xL,.5*\ya+.5*\yb) node {fluid: $\OmegaF$};
		\draw[thick,black] (4,\ya) node[anchor=south] {interface: $\Gamma$};

		\draw[-,thick,red,line width=2pt] (0  ,\ysb) -- (\xL,\ysb);
                \draw[-,thick,red,line width=2pt] (0  ,\ysa) -- (0,\ysb);
		\draw[-,thick,red,line width=2pt] (0  ,\ya) -- (0  ,\ya) node[anchor=east,black,yshift=-4pt] {\lbfont$y=0$};
		\draw[-,thick,red,line width=2pt] (\xL,\ysa) -- (\xL,\ysb);
		\draw[-,thick,blue] (0,\yb) -- (\xL,\yb);
		\draw[-,thick,blue] (\xL,0) -- (\xL,\yb);
		\draw (0  ,\ysa) node[anchor=north,black,yshift=-4pt] {\lbfont$x=0$};
		\draw (\xL,\ysa) node[anchor=north,black,yshift=-4pt] {\lbfont$x=L$};
		\draw (0,  \ysa) node[anchor=east,black] {\lbfont$y\rightarrow-\infty$};
		\draw (0,   \yb) node[anchor=east,black] {\lbfont$y=H$};
		%
		%
		\end{tikzpicture}
	\end{center}
	\caption{The rectangular geometry for the FSI model problems.} \label{fig:rectangularModelProblemCartoon}
\end{figure}
}

%% file: texFiles/numericalResults.tex
\section{Numerical results for an elastic piston} \label{sec:numericalResults}

We now present numerical results for two FSI problems to verify the accuracy and stability of the AMP scheme.  The two FSI problems considered involve the interaction of a fluid column with an elastic piston.  In the first problem, we examine longitudinal motion of the piston, while transverse motion of the piston is considered in the second problem.  Exact solutions are found for both FSI problems, and these are used to verify the accuracy and stability of the AMP algorithm for a range of the problem parameters.

\input texFiles/rectangularGeometryNumericalResults

\subsection{Longitudinal motion of an elastic piston}\label{sec:elasticPistonSolution}

The geometry of the elastic piston problem is shown in Figure~\ref{fig:rectangularModelProblemNumericalResults}.  The plot on the left shows the configuration at $t=0$.  The fluid occupies the physical domain between $y=0$ and $y=H$ initially, while the solid lies in its reference domain between $\ys=-\Hs$ and $\ys=0$.  It is assumed that there is no dependence in the $x$-direction so that the fluid-solid interface remains flat at a position $y=y_I(t)$ as shown in the plot on the right.  In the fluid domain, $\Omega(t)$, it is assumed that the horizontal component of velocity $v_1$ is zero, and thus the vertical component $v_2$ is a function of $t$ alone according to the continuity equation.  The fluid pressure is a linear function $y$, and is given by
\begin{equation}
p(y,t)={(H-y)p_I(t)+(y-y_I(t))p_H(t)\over H-y_I(t)},
\label{eq:oneDimFluidPressure}
\end{equation}
where $p_I(t)$ is the pressure on the interface and $p_H(t)$ is a specified fluid pressure at $y=H$.  The momentum equation for the fluid in the vertical direction reduces to
\begin{equation}
\rho \dot{v}_2=-p_y={p_I(t)-p_H(t)\over H-y_I(t)},
\label{eq:oneDimFluidMom}
\end{equation}
where $\rho$ is the (constant) fluid density.

With the displacement of the solid in the horizontal direction assumed to be zero, the equation for the vertical component of the displacement in the reference domain, $\bar\Omega(0)$, becomes
\begin{equation}
\us_{2,tt} = \cp^2 \us_{2,\ys\ys},\qquad -\Hs<\ys<0,
\label{eq:oneDimWaveEquation}
\end{equation}
where $\cp$ is the longitudinal wave speed.  The general solution of~\eqref{eq:oneDimWaveEquation}, assuming a zero-displacement condition at $\ys=-\Hs$, is
\begin{equation}
\us_2(\ys,t) = \bar F(t - (\ys+\Hs)/\cp ) -  \bar F(t + (\ys+\Hs)/\cp),
\label{eq:oneDimSolidDisplacement}
\end{equation}
where $\bar F(\bar\tau)$ is any smooth function.  The vertical position of the fluid-solid interface is determined by the solid displacement evaluated at $\ys=0$.  Thus,
\begin{equation}
y_I(t)=\us_2(0,t),
\label{eq:oneDimInterface}
\end{equation}
and the matching conditions on velocity and stress imply
\begin{equation}
v_2(t) = \us_{2,t}(0,t),\qquad -p_I(t)=(\lambdas + 2 \mus) \us_{2,\bar{y}}(0,t).
\label{eq:oneDimMatching}
\end{equation}
The interface conditions in~\eqref{eq:oneDimInterface} and~\eqref{eq:oneDimMatching} can be used in~\eqref{eq:oneDimFluidMom} to eliminate the fluid variables, which leads to a boundary condition involving the solid displacement and its derivatives.  This boundary condition can be written in the form
\begin{equation}
\rho\bigl(H-\us_2(0,t)\bigr)\us_{2,tt}(0,t)+(\lambdas + 2 \mus) \us_{2,\bar{y}}(0,t)=-p_H(t).
\label{eq:oneDimSolidBC}
\end{equation}

We are interested in finding a function $\bar F(\bar\tau)$ in~\eqref{eq:oneDimSolidDisplacement} so that $\us_2(\ys,t)$ satisfies initial conditions on $\us_2$ and $\us_{2,t}$, and also satisfies the nonlinear boundary condition in~\eqref{eq:oneDimSolidBC} for a specified fluid pressure $p_H(t)$.  While this could be done in principle, we choose instead to construct an exact solution by making a choice for $\bar F(\bar\tau)$, and then backing out the corresponding initial conditions on $\us_2(\ys,t)$ and the fluid pressure $p_H(t)$.  This is done in our verification tests using the choice
\[
\bar F(\bar\tau)=\alpha\cos(\omega\bar\tau),
\]
where $\alpha$ and $\omega$ are parameters.  With $\us_2(\ys,t)$ and $p_H(t)$ known from~\eqref{eq:oneDimSolidDisplacement} and~\eqref{eq:oneDimSolidBC}, respectively, the remaining fluid variables can be obtained from~\eqref{eq:oneDimFluidPressure}, \eqref{eq:oneDimInterface} and~\eqref{eq:oneDimMatching}.  In particular, we note that for this choice, the position of the fluid-solid interface oscillates with frequency $\omega$ and an amplitude $a$ given by
\[
y_I(t)=a\sin(\omega t),\qquad a=2\alpha\sin(\omega\Hs/\cp).
\]

\input tables/ElasticPistonConvergence.tex

Numerical results are obtained for the case $H=1, \rho=1$ and $\mu=0.01$ for the fluid, 
and using $\Hs=0.5$ and $\mus=\lambdas=\rhos=\delta$ for the solid.  The interface position is specified by
$a=0.1$ and $\omega=2\pi$. 
The density ratio, $\rhos/\rho = \delta$, is taken to be $10\sp{-3}$, 1 and $10\sp3$, 
representing FSI problems with light, moderate and heavy solids, respectively.  
Numerical solutions are computed using the AMP algorithm on a two-dimensional rectangular 
configuration (as shown in Figure~\ref{fig:rectangularModelProblemNumericalResults}) with 
periodic boundary conditions taken at $x=0$ and $x=L$ consistent with a one-dimensional solution.  
Table~\ref{fig:ElasticPistonConvergence} gives the maximum-norm errors for solutions computed 
using the AMP algorithm at $t_{{\rm final}}=0.6$ with grid resolutions $h = 1/(20j)$ for
$j=1,2,4,8.$  
The errors in the table indicate that the solution is 
converging at second-order accuracy.
The solution is plotted in Figure \ref{fig:elasticPistonPlot} for the case $\delta=1$.

{
\begin{figure}
\includegraphics[width=.49\textwidth]{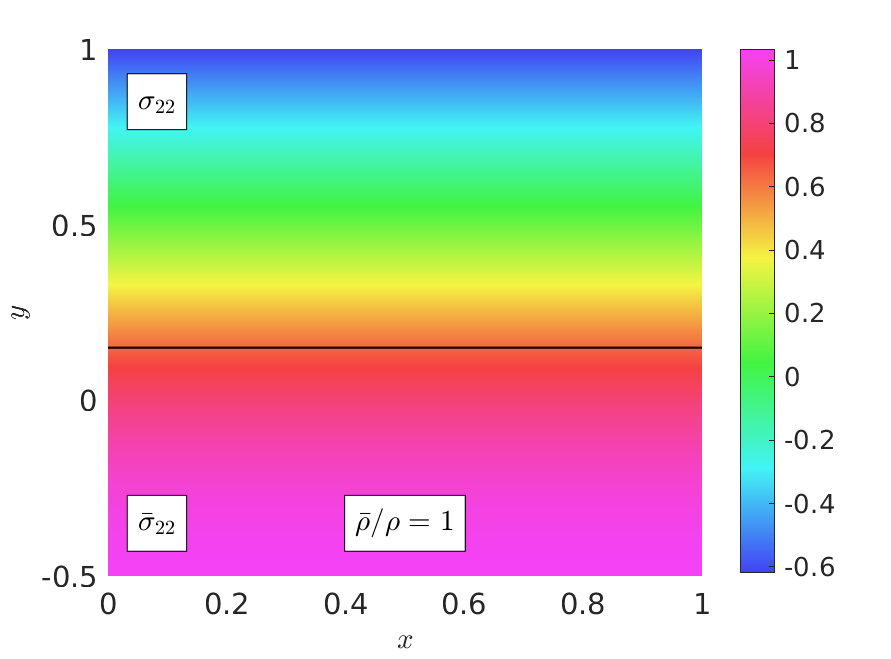}
\includegraphics[width=.49\textwidth]{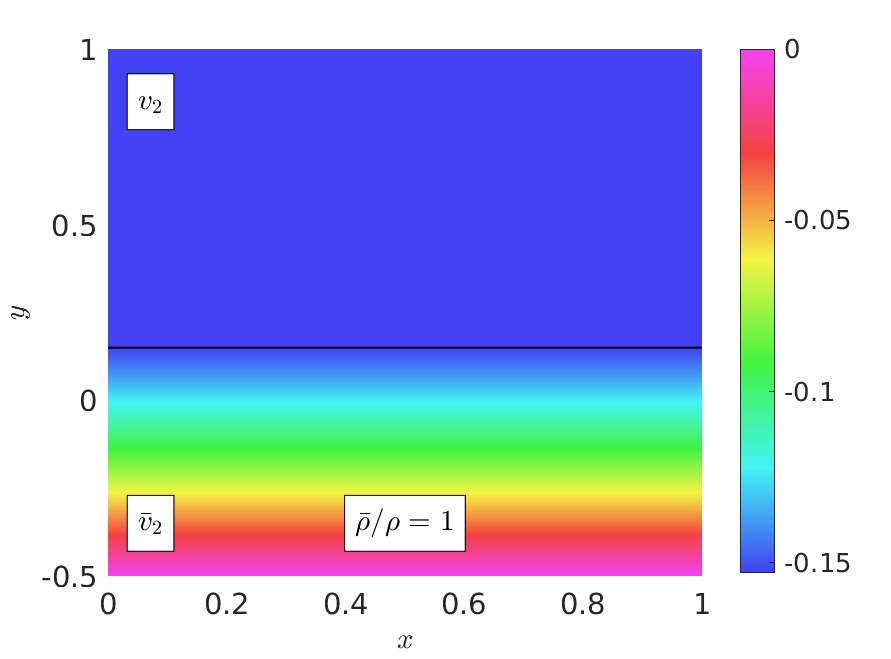}
\caption{Longitudinal motion of an elastic piston: Shaded contours of the numerical solution at  $t=0.6$ for
  $(\sigma_{22},\sigmas_{22})$ on the left and $(v_2,\vs_2)$ on the right computed using the AMP algorithm for $\rhos / \rho = 1.$ \label{fig:elasticPistonPlot}}
\end{figure}
}

%


\subsection{Transverse motion of an elastic piston}\label{sec:elasticPistonShearSolution}

Exact solutions can also be constructed for an FSI problem involving transverse motion of an elastic piston, see Figure~\ref{fig:rectangularModelProblemShearSolution}.  For this case, the vertical components of the fluid velocity and solid displacement are taken to be zero, while the corresponding horizontal components are assumed to be functions of $y$ and $t$ alone.  As a result, the interface remains in its initial position, $y_I(t) = y_I(0) = 0$, and so the solid reference coordinate~$\ys$ is equivalent to the physical coordinate~$y$.

\input texFiles/rectangularGeometryShearSolution

For this problem, the equations governing the horizontal components of the fluid velocity and solid displacement reduce to
\begin{equation}
v_{1,t} = \nu v_{1,yy}, \qquad 0<y<H, 
\label{eq:oneDimTransFluid}
\end{equation}
and
\begin{equation}
\us_{1,tt} = \cs^2 \us_{1,yy}, \qquad -\Hs<y<0, 
\label{eq:oneDimTransSolid}
\end{equation}
respectively, where $\nu=\mu/\rho$ is the kinematic viscosity of the fluid and $\cs$ is shear wave speed of the solid.  Solutions of~\eqref{eq:oneDimTransFluid} and~\eqref{eq:oneDimTransSolid} are sought in the form
\begin{equation}
v_1(y,t) = \hat{v}_1(y) e^{i \omega t}, \qquad \us_1(y,t) = \hat{\us}_1(y) e^{i \omega t},
\label{eq:oneDimTransSolutions}
\end{equation}
where $\omega$ is a parameter (ultimately an eigenvalue).  The coefficient functions, $\hat{v}_1(y)$ and $\hat{\us}_1(y)$, in~\eqref{eq:oneDimTransSolutions} satisfy
\begin{equation}
\hat{v}_1'' + \lambda^2 \hat{v}_1 = 0,\qquad \hat{\us}_1'' + k_s^2 \hat{\us}_1 = 0,
\label{eq:oneDimTransHelmholtz}
\end{equation}
where
\begin{equation}
\lambda^2 = - \frac{i \omega}{\nu}, \qquad k_s^2 = \frac{\omega^2}{\cs^2}.
\label{eq:oneDimTransParams}
\end{equation}
Solutions of the second-order ODEs in~\eqref{eq:oneDimTransHelmholtz} satisfying a no-slip condition in the fluid at $y=H$ and a zero-displacement condition in the solid at $y=-\Hs$ are
\begin{equation}
\hat{v}_1(y) =  b \sin(\lambda (H-y)),\qquad \hat{\us}_1(y) = \bar{b}\sin(k_s (\Hs+y)),
\label{eq:oneDimTransHelmholtzSolutions}
\end{equation}
where $(b,\bar{b})$ are constants.  The interface matching conditions on velocity and shear stress imply
\begin{equation}
\hat{v}_1(0) = i\omega\hat{\us}_1(0),  \qquad \mu \hat{v}_1'(0) = \mus \hat{\us}_1'(0).
\label{eq:oneDimTransMatching}
\end{equation}
Using~\eqref{eq:oneDimTransHelmholtzSolutions} in~\eqref{eq:oneDimTransMatching} leads to the system
\begin{equation}
\left[\begin{array}{cc}
\sin(\lambda H) & -i \omega \sin(k_s \Hs) \smallskip\\
\mu \lambda  \cos(\lambda H) & \mus k_s \cos(k_s \Hs)
\end{array}\right]
\left[\begin{array}{c}
b \smallskip\\ \bar{b}
\end{array}\right]
= 0 .
\label{eq:oneDimTransSystem}
\end{equation}
Nontrivial solutions for $(b,\bar{b})$ exist if the determinant of the coefficient matrix in~\eqref{eq:oneDimTransSystem} is zero, which implies
\begin{equation}
\mathcal{D}(\omega)=\frac{\mus}{\mu \omega} \tan(\lambda H) +  \frac{i \lambda}{k_s} \tan(k_s \Hs) = 0.
\label{eq:oneDimTransConstraint}
\end{equation}
Assuming a value for $\omega$ is found satisfying~\eqref{eq:oneDimTransConstraint}, the constants $b$ and $\bar{b}$ are given by
\begin{align}
\bar{b} = \frac{\us_0}{\sin(k_s H)}, \qquad b = \frac{i \omega \us_0}{\sin(\lambda H)},
\end{align}
where $\us_0$ is the interface displacement at $t=0$.
The fluid pressure is constant and equal to zero due to the matching of the
normal component of stress at the interface ($\sigma_{22} = \sigmas_{22} = 0$).
The real part of~\eqref{eq:oneDimTransSolutions} is used as the solution.

{
\begin{table}[h]
\caption{\label{fig:ElasticShearFrequencies}
Solutions to the dispersion relation in~\eqref{eq:oneDimTransConstraint} for $H=1$, $\Hs = 0.5$, $\rho = 1$, $\mu=0.1$ and $\delta=\rhos=\mus=\lambdas$.}
\begin{center}
\begin{tabular}{cc}
\hline
$\delta$ & $\omega$ \\
\hline
$10^{3}$ & $3.141+7.930 \cdot 10^{-4}$\\
$10^{0}$ & $2.351+5.433 \cdot 10^{-1}$\\
$10^{-3}$ & $6.285+1.784 \cdot 10^{-3}$\\
\hline
\end{tabular}
\end{center}
\end{table}
}

\input tables/ShearSolutionConvergence.tex

Values of $\omega$ with $\Re(\omega) > 0$ satisfying $\mathcal{D}(\omega)=0$ correspond to solutions of the transverse elastic piston problem that decay in time.  Selected values of $\omega$ are listed in Table~\ref{fig:ElasticShearFrequencies} for $H=1$, $\Hs=0.5$, $\rho=1$ and $\mu=0.1$, and for different values of $\delta=\rhos=\mus=\lambdas$. 
Table~\ref{fig:ShearSolutionConvergence} gives the maximum-norm errors for solutions computed 
using the AMP algorithm for $\us_0=0.1$.  The results are presents for solutions at $t_{{\rm final}}=0.3$ using grid resolutions $h = 1/(20j)$, for
$j=1,2,4,8.$  
The errors in the table indicate that the solution is 
converging at second-order accuracy.
The solution is plotted in Figure \ref{fig:elasticShearPlot} for the case $\delta=1$.

{
\begin{figure}
\includegraphics[width=.49\textwidth]{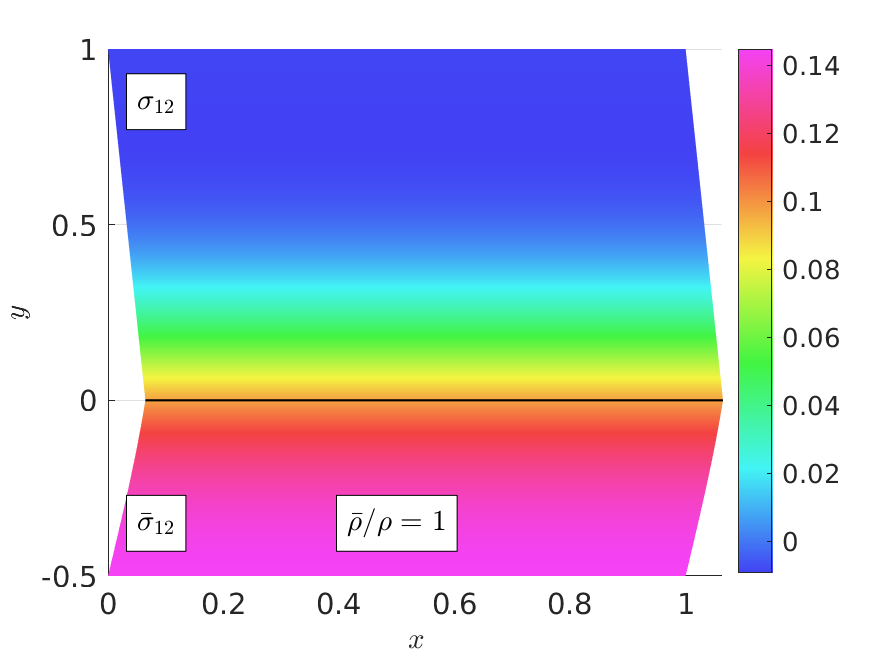}
\includegraphics[width=.49\textwidth]{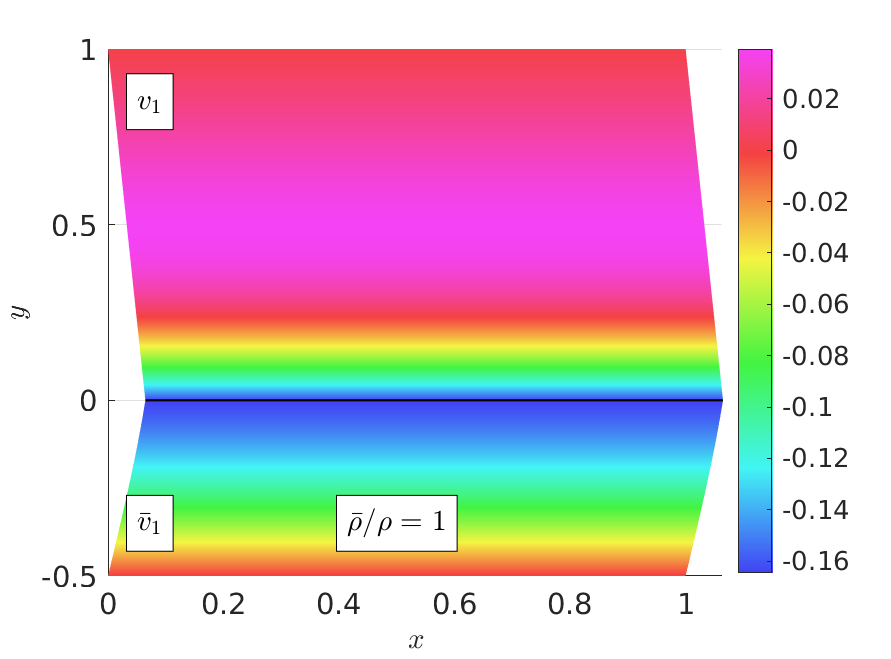}
\caption{Shaded contours of the numerical solution for
  $(\sigma_{12},\sigmas_{12})$ on the left and $(v_1,\vs_1)$ on the right at 
  $t=0.3$ for $\delta = 1.$ \label{fig:elasticShearPlot}}
\end{figure}
}



%% file: texFiles/rectangularGeometryNumericalResults.tex
{
\newcommand{\lbfont}{\small}

\def\ysa{-2.5}   

\def\ya{-.3} 
\def\ysb{-.45} 

\def\yaDef{.3}
\def\ysbDef{.15} 

\def\yb{2.5} 
\def\xL{4}
\def\xR{8}
\begin{figure}[hbt]
	\newcommand{\textFont}{\normalss}
	\begin{center}
		\begin{tikzpicture}[scale=.9]
		\useasboundingbox (0-.5,\ysa) rectangle (\xR+\xL+.5,\yb+.5);  
                \draw[fill=red!20,draw=red,line width=2pt] (0,\ysa) rectangle (\xL,\ysb);
		\draw[thick,fill=blue!20,draw=blue,line width=2pt] (0,\ya) rectangle (\xL,\yb);

		\draw[thick,black] (\xL/2,.5*\ysa+.5*\ysb) node {solid: $\bar{\Omega}(0)$};
                \draw[thick,black] (.5*\xL,.5*\ya+.5*\yb) node {fluid: $\Omega(0)$};
		\draw[thick,black] (\xL/2,\ya) node[anchor=south] {interface: $\Gamma(0)$};
		\draw[-,thick,red,line width=2pt] (0  ,\ya) -- (0  ,\ya) node[anchor=east,black,yshift=-4pt] {\lbfont$\ys=0$};

		\draw (0  ,\ysa) node[anchor=north,black,yshift=-4pt] {\lbfont$x=0$};
		\draw (\xL,\ysa) node[anchor=north,black,yshift=-4pt] {\lbfont$x=L$};
		\draw (0,  \ysa) node[anchor=east,black] {\lbfont$\ys=-\Hs$};
		\draw (0,   \yb) node[anchor=east,black] {\lbfont$\ys=H$};

                \draw[fill=red!20,draw=red,line width=2pt] (\xR+0,\ysa) rectangle (\xR+\xL,\ysbDef);
		\draw[thick,fill=blue!20,draw=blue,line width=2pt] (\xR+0,\yaDef) rectangle (\xR+\xL,\yb);

		\draw[thick,black] (\xR+\xL/2,.5*\ysa+.5*\ysbDef) node {solid: $\bar{\Omega}(t)$};
                \draw[thick,black] (\xR+.5*\xL,.5*\yaDef+.5*\yb) node {fluid: $\Omega(t)$};
		\draw[thick,black] (\xR+\xL/2,\yaDef) node[anchor=south] {interface: $\Gamma(t)$};
		\draw[-,thick,red,line width=2pt] (\xR+0  ,\yaDef) -- (\xR+0  ,\yaDef) node[anchor=east,black,yshift=-4pt] {\lbfont$y=y_I(t)$};

		\draw (\xR+0  ,\ysa) node[anchor=north,black,yshift=-4pt] {\lbfont$x=0$};
		\draw (\xR+\xL,\ysa) node[anchor=north,black,yshift=-4pt] {\lbfont$x=L$};
		\draw (\xR+0,  \ysa) node[anchor=east,black] {\lbfont$y=-\Hs$};
		\draw (\xR+0,   \yb) node[anchor=east,black] {\lbfont$y=H$};

		\end{tikzpicture}
	\end{center}
	\caption{Longitudinal motion of an elastic piston: initial solid and fluid domains (left), and deformed piston and fluid domain (right) for $t>0$.} \label{fig:rectangularModelProblemNumericalResults}
\end{figure}
}

%% file: tables/ElasticPistonConvergence.tex
\begin{table}[H]
\caption{Longitudinal motion of an elastic piston: Maximum-norm errors and convergence rates of the numerical solution at $t_{\rm final} = 0.6$ computed using the AMP algorithm for $\rhos/\rho = \delta= 10^3$, $1$ and $10^{-3}.$ \label{fig:ElasticPistonConvergence}}
\begin{center}
{\footnotesize
\begin{tabular}{c}
Heavy solid ($\delta=10^3$):
\end{tabular}
\begin{tabular}{ccccccccccc}
\hline
$h$ & $E^{(p)}$ & $r$ & $E^{(\vv)}$ & $r$ & $E^{(\usv)}$ & $r$ & $E^{(\vsv)}$ & $r$ & $E^{(\sigmasv)}$ & $r$\\
\hline
1/ 20 &  6.0e-04 &      &  5.9e-05 &      &  4.0e-05 &      &  5.9e-05 &      &  1.9e-01 &     \\
1/ 40 &  1.4e-04 &  4.2 &  1.7e-05 &  3.5 &  9.6e-06 &  4.2 &  1.7e-05 &  3.5 &  4.5e-02 &  4.2 \\
1/ 80 &  3.4e-05 &  4.1 &  4.4e-06 &  3.8 &  2.3e-06 &  4.1 &  4.4e-06 &  3.8 &  1.1e-02 &  4.1 \\
1/160 &  8.5e-06 &  4.1 &  1.1e-06 &  3.9 &  5.8e-07 &  4.1 &  1.1e-06 &  3.9 &  2.7e-03 &  4.1 \\
\hline
\\
\end{tabular}
\begin{tabular}{c}
Medium solid ($\delta=1$):
\end{tabular}
\begin{tabular}{cccccccccccc}
\hline
$h$ & $E^{(p)}$ & $r$ & $E^{(\vv)}$ & $r$ & $E^{(\usv)}$ & $r$ & $E^{(\vsv)}$ & $r$ & $E^{(\sigmasv)}$ & $r$\\
\hline
1/ 20 &  1.8e-05 &      &  4.9e-05 &      &  1.2e-05 &      &  4.9e-05 &      &  5.0e-05 &     \\
1/ 40 &  7.5e-06 &  2.4 &  1.2e-05 &  4.0 &  3.0e-06 &  4.2 &  1.2e-05 &  4.0 &  1.3e-05 &  3.7 \\
1/ 80 &  2.3e-06 &  3.3 &  3.0e-06 &  4.0 &  7.1e-07 &  4.1 &  3.0e-06 &  4.0 &  3.6e-06 &  3.8 \\
1/160 &  6.3e-07 &  3.6 &  7.4e-07 &  4.0 &  1.8e-07 &  4.1 &  7.4e-07 &  4.0 &  9.2e-07 &  3.9 \\
\hline
\\
\end{tabular}

\begin{tabular}{c}
Light solid ($\delta=10^{-3}$):
\end{tabular}
\begin{tabular}{cccccccccccc}
\hline
$h$ & $E^{(p)}$ & $r$ & $E^{(\vv)}$ & $r$ & $E^{(\usv)}$ & $r$ & $E^{(\vsv)}$ & $r$ & $E^{(\sigmasv)}$ & $r$\\
\hline
1/ 20 &  8.0e-07 &      &  6.5e-07 &      &  3.3e-06 &      &  2.4e-05 &      &  1.3e-07 &     \\
1/ 40 &  2.4e-07 &  3.3 &  1.6e-07 &  4.0 &  5.3e-07 &  6.3 &  4.2e-06 &  5.7 &  3.4e-08 &  3.9 \\
1/ 80 &  6.6e-08 &  3.7 &  4.1e-08 &  4.0 &  8.9e-08 &  5.9 &  8.3e-07 &  5.0 &  8.8e-09 &  3.9 \\
1/160 &  1.7e-08 &  3.8 &  1.0e-08 &  4.0 &  2.3e-08 &  3.8 &  1.8e-07 &  4.5 &  2.2e-09 &  4.0 \\
\hline
\end{tabular}
}
\end{center}
\end{table}

%% file: texFiles/rectangularGeometryShearSolution.tex
{
\newcommand{\lbfont}{\small}

\def\ys{\bar{y}}
\def\ysa{-2.5}   

\def\ya{-.3} 
\def\ysb{-.45} 

\def\yaDef{-.3}
\def\ysbDef{-.45} 

\def\yb{2.5} 
\def\xL{4}
\def\xR{6}

\def\xDef{.3}

\begin{figure}[hbt]
	\newcommand{\textFont}{\normalss}
	\begin{center}
		\begin{tikzpicture}[scale=.9,every node/.style={transform shape}]
                  \useasboundingbox (0-.5,\ysa) rectangle (\xR+\xL+.5,\yb+.5);  
                \draw[fill=red!20,draw=red,line width=2pt] (0,\ysa) rectangle (\xL,\ysb);
		\draw[thick,fill=blue!20,draw=blue,line width=2pt] (0,\ya) rectangle (\xL,\yb);

		\draw[thick,black] (\xL/2,.5*\ysa+.5*\ysb) node {solid: $\bar{\Omega}_0$};
                \draw[thick,black] (.5*\xL,.5*\ya+.5*\yb) node {fluid: $\Omega_0$};
		\draw[thick,black] (\xL/2,\ya) node[anchor=south] {interface: $\Gamma_0$};
		\draw[-,thick,red,line width=2pt] (0  ,\ya) -- (0  ,\ya) node[anchor=east,black,yshift=-4pt] {\lbfont$\ys=0$};

		\draw (0  ,\ysa) node[anchor=north,black,yshift=-4pt] {\lbfont$x=0$};
		\draw (\xL,\ysa) node[anchor=north,black,yshift=-4pt] {\lbfont$x=L$};
		\draw (0,  \ysa) node[anchor=east,black] {\lbfont$\ys=-\Hs$};
		\draw (0,   \yb) node[anchor=east,black] {\lbfont$\ys=H$};

                \draw[fill=red!20,draw=red,line width=2pt] (\xR+0,\ysa) to[out=90,in=-100] 
                (\xR+\xDef+0,\ysbDef) to[out=0,in=180] 
                (\xR+\xDef+\xL,\ysbDef) to[out=-100,in=-270]
                (\xR+\xL,\ysa) to[out=180,in=0]
                (\xR+0,\ysa);

                \draw[fill=blue!20,draw=blue,line width=2pt] (\xR+\xDef,\ya) to[out=100,in=-85] 
                (\xR+0,\yb) to[out=0,in=180] 
                (\xR+\xL,\yb) to[out=-85,in=100]
                (\xR+\xDef+\xL,\ya) to[out=180,in=0]
                (\xR+\xDef+0,\ya);

		\draw[thick,black] (\xR+\xL/2,.5*\ysa+.5*\ysbDef) node {solid: $\bar{\Omega}(t)$};
                \draw[thick,black] (\xR+.5*\xL,.5*\yaDef+.5*\yb) node {fluid: $\Omega(t)$};
		\draw[thick,black] (\xR+\xL/2,\yaDef) node[anchor=south] {interface: $\Gamma(t)$};
		\draw[-,thick,red,line width=2pt] (\xR+0  ,\yaDef) -- (\xR+0  ,\yaDef) node[anchor=east,black,yshift=-4pt] {\lbfont$y=0$};

		\draw (\xR+0  ,\ysa) node[anchor=north,black,yshift=-4pt] {\lbfont$x=0$};
		\draw (\xR+\xL,\ysa) node[anchor=north,black,yshift=-4pt] {\lbfont$x=L$};
		\draw (\xR+0,  \ysa) node[anchor=east,black] {\lbfont$y=-\Hs$};
		\draw (\xR+0,   \yb) node[anchor=east,black] {\lbfont$y=H$};

		\end{tikzpicture}
	\end{center}
	\caption{Transverse motion of an elastic piston: initial solid and fluid domains (left), and deformed piston and fluid domain (right) for $t>0$.} \label{fig:rectangularModelProblemShearSolution}

\end{figure}

}

%% file: tables/ShearSolutionConvergence.tex
\begin{table}[H]
\caption{Transverse motion of an elastic piston: Maximum-norm errors and convergence rates of the numerical solution computed using the AMP algorithm for $\rhos/\rho = \delta= 10^3$, $1$ and $10^{-3}.$  The exact solutions use values of $\omega$ from Table~\ref{fig:ElasticShearFrequencies} and $t_{\rm final} = 0.3$.\label{fig:ShearSolutionConvergence}}
\begin{center}
{\footnotesize
\begin{tabular}{c}
Heavy solid ($\delta=10^3$): 
\end{tabular}
\begin{tabular}{cccccccccccc}
\hline
$h$ & $E^{(p)}$ & $r$ & $E^{(\vv)}$ & $r$ & $E^{(\usv)}$ & $r$ & $E^{(\vsv)}$ & $r$ & $E^{(\sigmasv)}$ & $r$\\
\hline
1/ 20 &  6.0e-04 &      &  5.9e-05 &      &  4.0e-05 &      &  5.9e-05 &      &  1.9e-01 &     \\
1/ 40 &  1.4e-04 &  4.2 &  1.7e-05 &  3.5 &  9.6e-06 &  4.2 &  1.7e-05 &  3.5 &  4.5e-02 &  4.2 \\
1/ 80 &  3.4e-05 &  4.1 &  4.4e-06 &  3.8 &  2.3e-06 &  4.1 &  4.4e-06 &  3.8 &  1.1e-02 &  4.1 \\
1/160 &  8.5e-06 &  4.1 &  1.1e-06 &  3.9 &  5.8e-07 &  4.1 &  1.1e-06 &  3.9 &  2.7e-03 &  4.1 \\
\hline
\\
\end{tabular}

\begin{tabular}{c}
Medium solid ($\delta=1$): 
\end{tabular}
\begin{tabular}{cccccccccccc}
\hline
$h$ & $E^{(p)}$ & $r$ & $E^{(\vv)}$ & $r$ & $E^{(\usv)}$ & $r$ & $E^{(\vsv)}$ & $r$ & $E^{(\sigmasv)}$ & $r$\\
\hline
1/ 20 &  1.8e-05 &      &  4.9e-05 &      &  1.2e-05 &      &  4.9e-05 &      &  5.0e-05 &     \\
1/ 40 &  7.5e-06 &  2.4 &  1.2e-05 &  4.0 &  3.0e-06 &  4.2 &  1.2e-05 &  4.0 &  1.3e-05 &  3.7 \\
1/ 80 &  2.3e-06 &  3.3 &  3.0e-06 &  4.0 &  7.1e-07 &  4.1 &  3.0e-06 &  4.0 &  3.6e-06 &  3.8 \\
1/160 &  6.3e-07 &  3.6 &  7.4e-07 &  4.0 &  1.8e-07 &  4.1 &  7.4e-07 &  4.0 &  9.2e-07 &  3.9 \\
\hline
\\
\end{tabular}

\begin{tabular}{c}
Light solid ($\delta=10^{-3}$): 
\end{tabular}
\begin{tabular}{cccccccccccc}
\hline
$h$ & $E^{(p)}$ & $r$ & $E^{(\vv)}$ & $r$ & $E^{(\usv)}$ & $r$ & $E^{(\vsv)}$ & $r$ & $E^{(\sigmasv)}$ & $r$\\
\hline
1/ 20 &  8.0e-07 &      &  6.5e-07 &      &  3.3e-06 &      &  2.4e-05 &      &  1.3e-07 &     \\
1/ 40 &  2.4e-07 &  3.3 &  1.6e-07 &  4.0 &  5.3e-07 &  6.3 &  4.2e-06 &  5.7 &  3.4e-08 &  3.9 \\
1/ 80 &  6.6e-08 &  3.7 &  4.1e-08 &  4.0 &  8.9e-08 &  5.9 &  8.3e-07 &  5.0 &  8.8e-09 &  3.9 \\
1/160 &  1.7e-08 &  3.8 &  1.0e-08 &  4.0 &  2.3e-08 &  3.8 &  1.8e-07 &  4.5 &  2.2e-09 &  4.0 \\
\hline
\end{tabular}
}
\end{center}
\end{table}

%% file: texFiles/conclusions.tex
\section{Conclusions} \label{sec:conclusions}


A stable added-mass partitioned (AMP) algorithm was developed for fluid-structure interaction
problems involving viscous incompressible flow and compressible elastic solids.  The new algorithm
is stable, without sub-time-step iterations, for both heavy and very light solids and effectively
suppresses both added-mass and added-damping effects.
The fluid is advanced using a fractional-step IMEX scheme with the viscous terms treated 
implicitly.
Key elements of the new AMP scheme are a Robin interface condition for
the pressure and an impedance-weighted interface projection based on a new form for the fluid
impedance. The fluid impedance is derived from an analysis of a solution to an FSI problem 
near the interface. Stability of the AMP scheme is analyzed for two model problems, 
one involving a first-order accurate discretization of a viscous model problem 
and a second-order accurate discretization of an inviscid model problem.
The \textit{elastic-piston}, a set of new benchmark problems, 
was developed to verify stability and accuracy of the AMP scheme. 
These solutions are exact and include finite interface deformations
either normal or tangent to the surface.

%% file: texFiles/analysisAppendix-new.tex
\section{Results supporting the stability analysis of the model problems} \label{sec:StabilityAnalysisDetails}
The subsections below give additional results used for the stability analysis of the viscous and inviscid model problems of Sections~\ref{sec:viscousAnalysis} and~\ref{sec:inviscidAnalysis}, respectively.

\subsection{CFL region for Cauchy problem} \label{sec:CFLRegionCauchy}

The CFL stability regions are describedfor the discretizations of the model problems discussed in 
Sections~\ref{sec:viscousAnalysis} and~\ref{sec:inviscidAnalysis} applied to the 
pure initial-value problem (Cauchy problem). 
For both discretization, we assume that the spatial eigenfunction is given by $\phi=e\sp{i\vartheta}$
for all $\vartheta\in[0,2\pi]$, i.e.~$\vert\phi\vert=1$.  For the viscous model problem, the amplification factor $A$ can be found
from~\eqref{eq:SolidDetCond} for a given value of $\phi$, while~\eqref{eq:quarticPhi} gives $A$ in terms
of $\phi$ for the inviscid model problem.
The CFL stability region is defined as the region in the $(\lx,\ly)$ plane for which the
amplification factor satisfies $|A|\le 1$.
This region can be found numerically by computing the $A_{{\rm max}}=\max_{\vert\phi\vert=1}\vert A\vert$
for a range of values in the
$(\lx,\ly)$ plane.   The left plot of 
Figure~\ref{fig:viscousCFL} shows shaded countours of $A_{{\rm max}}$, while the plot on
the right indicates the region for which $A_{{\rm max}}\le1$.  Note that the constraint $\lx\sp2+\ly\sp2\le1$ is sufficient for CFL stability.  Similar results are shown in
Figure~\ref{fig:inviscidCFL} for the analysis of the inviscid model problem.

{
\newcommand{\figWidth}{6cm}
\begin{figure}[h]
\begin{center}
\includegraphics[width=\figWidth]{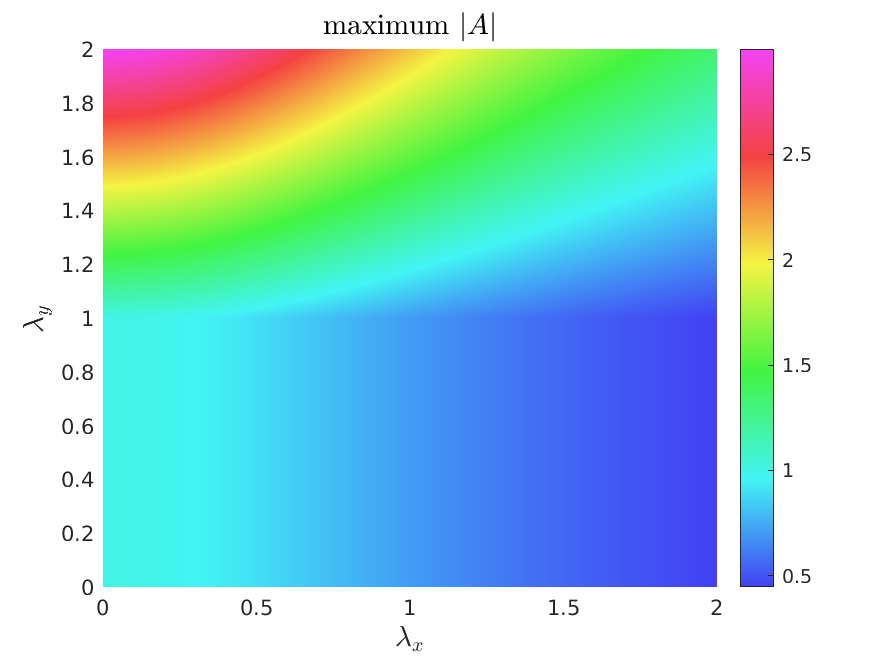}
\includegraphics[width=\figWidth]{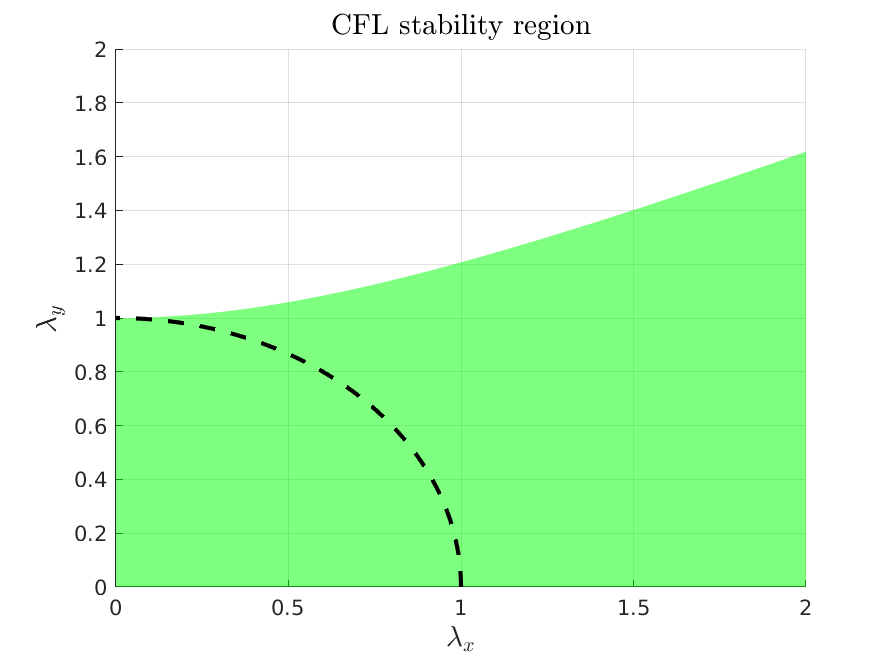}
\caption{CFL stability for the discretization of the viscous model problem.
Left: shaded contours of the maximum amplification factor, $A_{{\rm max}}$
for $\phi=e\sp{i\vartheta}$ with $\vartheta \in [0,2 \pi]$. 
Right: green shaded region indicates the stability region for which $A_{{\rm max}}\le 1.$ 
The dotted line shows the curve $\lx^2 + \ly^2 = 1$ which lies in 
the stability region. \label{fig:viscousCFL}}
\end{center}
\end{figure}

\begin{figure}[h]
\begin{center}
\includegraphics[width=\figWidth]{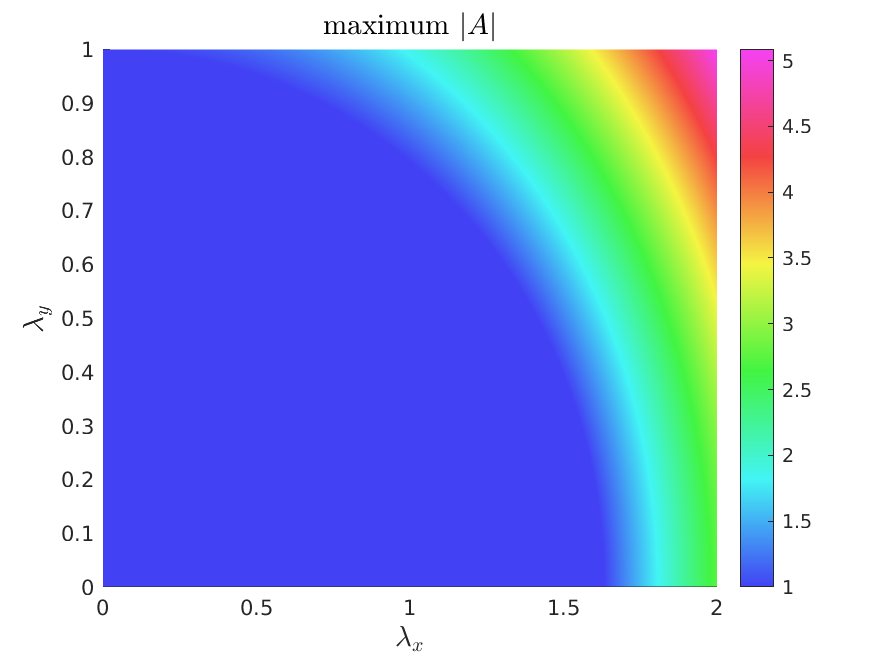}
\includegraphics[width=\figWidth]{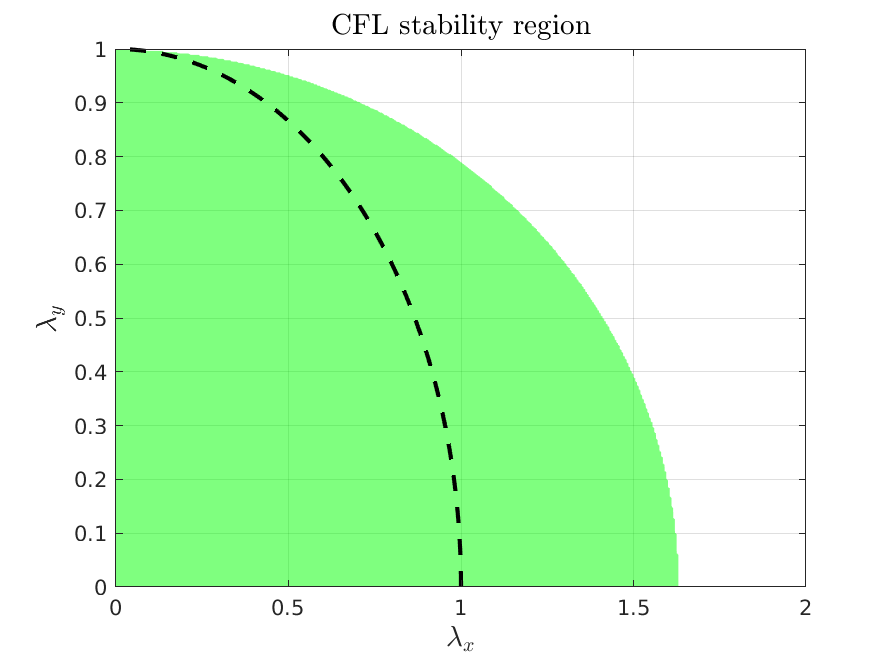}
\caption{CFL stability for the discretization of the inviscid model problem.
Left: shaded contours of the maximum amplification factor, $A_{{\rm max}}$
for $\phi=e\sp{i\vartheta}$ with $\vartheta \in [0,2 \pi]$. 
Right: green shaded region indicates the stability region for which $A_{{\rm max}}\le 1.$ 
The dotted line shows the curve $\lx^2 + \ly^2 = 1$ which lies in 
the stability region. \label{fig:inviscidCFL}}
\end{center}
\end{figure}
}
\subsection{Matrix system in the stability analysis of the viscous model problem} \label{sec:ComponentsOfDmatViscous}

The stability analysis of the viscous model problem in Section~\ref{sec:viscousAnalysis} led to the matrix system
\[
\dMat(A) \kv = \left[\begin{array}{cc}
\dMat_{11} & \dMat_{12} \\ 
\dMat_{21} & \dMat_{22}
\end{array}\right]
\left[\begin{array}{c}
k_1 \\ k_2
\end{array}\right]
= 0,
\]
where the components $\dMat_{ij}$ depend on the choice of coupling at the interface for the partitioned scheme.  For the AMP scheme, the components can be written in the form
\begin{align*}
\dMat_{11} &= i (\phi-1) \eta(1/\phi) \left[\beta_0+\beta_1Z+\beta_2Z\sp2+\beta_3Z\sp3\right],
\\
\dMat_{12} &= \dMat_{21} = 0,
\\
\dMat_{22} &= -\lx (\phi-1) A \left[\beta_0+\beta_1Z+\beta_2Z\sp2+\beta_3Z\sp3\right],
\end{align*}
where
\begin{align*}
\beta_0=& \left[ (\gamma+1)A\sp2-A \right] {\Lambda}^{2}, \\
\beta_1=& \left[  \left( \gamma+1 \right)  \left( 2{\Lambda}^{2}+
 \left( \gamma+3 \right) \Lambda +2 \right) {A}^{2}- \left( 2\Lambda+1 \right)  \left( \Lambda \left( \gamma+1 \right) +\gamma+2 \right) A
+2\Lambda+1 \right] \Lambda, \\
\beta_2=& \left( \gamma+1 \right)  \left( 4\gamma{\Lambda}^{3}+
 2\left( 2\gamma+1 \right) {\Lambda}^{2}+ \left( 2\gamma+3 \right) 
\Lambda +1\right) A\sp2 \\
& \quad  - \left( 4\left( {\gamma}^{2}+\gamma \right) {\Lambda}^{3}
 + 4\left( {\gamma}^{2}+\gamma \right) {\Lambda}^{2} + \left( {
\gamma}^{2}+3\gamma \right) \Lambda +\gamma+1 \right)A,  \\
\beta_3=& \left( 2\Lambda+1 \right)  \left[  \left( 2\Lambda+1 \right) 
 \left( \gamma\sp2+\gamma \right) {A}^{2}-\left( \gamma-1
 \right)  \left( \gamma+2\Lambda+2 \right) A -2\right].
\end{align*}
%
For the TP scheme, the components of the linear system have the form 
\begin{align*}
\dMat_{11} =& i\eta(1/\phi)\Lambda\left[(2\phi-1)\left((\gamma+1)A-1\right)+\gamma A(\gamma+1)Z \right]
              -i\lx Z\left[\Lambda(\gamma+1)A\sp2-(2\Lambda+1)A+1\right], \\
\dMat_{12} =&  -i\eta(\phi)\Lambda\left[(\gamma+1)A-1-\gamma A(\gamma+1)Z \right]
               +i\lx Z\left[\Lambda(\gamma+1)A\sp2-(2\Lambda+1)A+1\right], \\
\dMat_{21} =& \eta(1/\phi)\Lambda Z\left[\gamma A(\gamma+1)-2\gamma\right] \\
            &\quad  -\lx \left[\left((\gamma+1)A\sp2-A\right)\Lambda-\left((\gamma+1)\left((2\gamma+1)\Lambda+1\right)A\sp2
            -\left(2(\gamma+1)\Lambda+\gamma+3\right)A+2\right) Z \right], \\
\dMat_{22} =& \eta(\phi)\Lambda Z\left[\gamma A(\gamma+1)-2\gamma\right] \\
            &\quad +\lx \left[(1-2\phi)\left((\gamma+1)A\sp2-A\right)\Lambda-\left((\gamma+1)\left((2\gamma+1)\Lambda+1\right)A\sp2 - \left(2(\gamma+1)\Lambda+\gamma+3\right)A+2\right)Z\right].
\end{align*}
The components of the linear system for the ATP scheme are 
\begin{align*}
\dMat_{11} =& i\eta(1/\phi)\left[\gamma(\gamma+1)A+(2\phi-1)\left((\gamma+1)(\gamma\sp2+1)A+2(\gamma-1)\right)Z\right] + i\lx\left[(\gamma-1)A\right], \\
\dMat_{12} =& -i\eta(\phi)\left[\gamma(\gamma+1)A-\left((\gamma+1)(\gamma\sp2+1)A+2(\gamma-1)\right)Z\right] + i\lx\left[(\gamma-1)A\right], \\
\dMat_{21} =& -\eta(1/\phi)\left[(\gamma+1)A-2\right] -\lx\left[(\gamma+1)A-\left((\gamma+1)(\gamma\sp2+1)A^2+2(\gamma-1)A\right)Z\right],\\
\dMat_{22} =& \eta(\phi)\left[(\gamma+1)A-2\right] -\lx\left[(\gamma+1)A+(2\phi-1)\left((\gamma+1)(\gamma\sp2+1)A^2+2(\gamma-1)A\right)Z\right].
\end{align*}

\subsection{Matrix system in the stability analysis of the inviscid model problem} \label{sec:ComponentsOfDmat}

The stability analysis of the inviscid model problem in Section~\ref{sec:inviscidAnalysis} led to the matrix system
\[
\dMatI(A) \kv= \left[\begin{array}{cc}
\dComp_{11} & \dComp_{12} \\
\dComp_{21} & \dComp_{22}
\end{array}\right] 
\left[\begin{array}{c}
k_1 \\ k_2
\end{array}\right] = 0,
\]
where the components $\dComp_{ij}$ depend on the choice of coupling at the interface for the partitioned scheme as before.  The first row of $\dMat$ is associated with boundary condition on the incoming characteristic.  For the AMP scheme, the components of this row are given by
\begin{align}
\dComp_{1\nn} &= \frac{\delta_{2,\nn} M^2 + \delta_{1,\nn} M + \delta_{0,\nn} }
{\left(A^2 - \frac{4}{3} A + \frac{1}{3}\right) M^2
+ A^2 M + A^2} ,\qquad \nn=1,2,\qquad\hbox{[AMP scheme]}
\end{align}
where
\begin{align}
\delta_{0,\nn} &= A^2 \left(\qn{\nn} - \rn{\nn} - \qn{\nn} \left( \p_\nn + \frac{1}{\p_\nn} \right) \right),\\
\delta_{1,\nn} &= -A^2 \left( \frac{\qn{\nn}}{6} - \frac{5 \rn{\nn}}{6} 
+ \qn{\nn} \left(\p_\nn + \frac{1}{\p_\nn} \right)\right)
+ 2 A (\qn{\nn} - \rn{\nn})
+ \frac{1}{2} (\rn{\nn} - \qn{\nn}), \\
\delta_{2,\nn} &= \left(A^2 - \frac{4}{3} A + \frac{1}{3}\right)
\left(\qn{\nn} + \rn{\nn} - \qn{\nn} \left( \p_\nn + \frac{1}{\p_\nn}\right) \right) .
\end{align}
For the TP and ATP schemes, the components of the first row are given by
\begin{align}
\dComp_{1\nn} &= 
(\rn{\nn} +\qn{\nn}) \left(\p_\nn + \frac{1}{\p_\nn} \right)
- M \frac{3A^2-4A+1}{A^2} (\rn{\nn}-\qn{\nn}),
\qquad \hbox{[TP scheme]} \\
\dComp_{1\nn} &= 
(\rn{\nn} -\qn{\nn}) \left(\p_\nn + \frac{1}{\p_\nn} \right)
- \frac{1}{4 M} \frac{A^2}{3A^2 - 4A + 1} (\rn{\nn} + \qn{\nn}),
\qquad\hbox{[ATP scheme]}
\end{align}
where $\nn=1,2$.  The second row of $\dMat$ is common to all three schemes as it represents an extrapolation of the outgoing characteristic into the interface ghost point. This extrapolation condition gives
\begin{align}
\dComp_{2\nn} = \rn{\nn} \left(\p_\nn - 2 + \frac{1}{\p_\nn} \right), \qquad \nn = 1,2,\qquad\hbox{[all schemes]}.
\end{align}

\subsection{Exact solution for the 1D inviscid model problem} \label{sec:SolutionToOneDimensional}
The verification of the stability analysis of the inviscid model problem in Section~\ref{sec:inviscidAnalysis1D} requires the exact solution of the model problem.  Assuming that the solution has no variation in the $x$-direction, the ODE for the pressure with the boundary condition at $y=H$ gives
\begin{align}
p(y,t) = \rho \dot{v}_I(t) (H-y),\qquad 0<y<H,\quad t>0,
\label{eq:fluidPressureSoln}
\end{align}
where $v_I(t)$ is the interface velocity (to be determined).  In terms of the characteristic variables of the solid, the solution is
\begin{align}
a(y,t) = 
\begin{cases}
a_0(y+\cp t), & 0 < t < - y/\cp\\
a_I(t+y/\cp), & t > -y/\cp 
\end{cases}, \qquad 
 b(y,t) = b_0(y-\cp t),
\label{eq:solidCharSoln}
\end{align}
where $a_0(y)$ and $b_0(y)$, $y<0$, are given by the initial conditions for $\vs_2$ and $\sigmas_{22}$, and $a_I(t)$ specifies the incoming characteristic variable of the solid at the interface (to be determined).  Evaluating the outgoing characteristic of the solid at the interface gives $b_I(t)=b_0(-\cp t)$, and this can be matched with the fluid to give
\begin{align}
-\rho H \dot{v}_I - \zp v_I = b_I(t).
\label{eq:interfaceVelODE}
\end{align}
The solution of the ODE in~\eqref{eq:interfaceVelODE} for the interface velocity has the form
\begin{align}
v_I(t) = v_I(0) e^{-\lambda t} - \frac{1}{\rho H} \int_0^t e^{\lambda (\tau - t)} b_I(\tau) \, d \tau,
\end{align}
where $\lambda=\zp/(\rho H)$ and $v_I(0)$ is assumed to be given as an initial condition.  The fluid pressure in~\eqref{eq:fluidPressureSoln} is now specified, and the incoming characteristic of the solid at the interface can now be obtained from
\begin{align}
a_I(t) = -\rho H\dot{v}_I  + \zp v_I.
\end{align}
This completes the characteristic description of the exact solution for the solid in~\eqref{eq:solidCharSoln}, which gives
\begin{align}
\vs_2(y,t)={1\over2\zp}\left(a(y,t)-b(y,t)\right),\qquad \sigmas_{22}(y,t)={1\over2}\left(a(y,t)+b(y,t)\right),\qquad y<0,\quad t>0.
\end{align}
The initial conditions used to specify the exact solution in Section~\ref{sec:inviscidAnalysis1D} are
\[
a_0(y)=2\rhos\cp\sp2\exp\left(-(y/H)\sp4\right),\qquad b_0(y)=0,\qquad v_I(0)=2\cp,
\]
and the parameters of the model problem are taken as $\rho=H=\cp=1$, while $\rhos$ is varied as needed to generate the plots in Figure~\ref{fig:acousticSolidStokesStabilityTPNumerical}.

%% file: fibrmp.bbl
\begin{thebibliography}{10}
\expandafter\ifx\csname url\endcsname\relax
  \def\url#1{\texttt{#1}}\fi
\expandafter\ifx\csname urlprefix\endcsname\relax\def\urlprefix{URL }\fi

\bibitem{fib2014}
J.~W. Banks, W.~D. Henshaw, D.~W. Schwendeman, An analysis of a new stable
  partitioned algorithm for {FSI} problems. {Part I}: Incompressible flow and
  elastic solids, J. Comput. Phys. 269 (2014) 108--137\citeCount{39}.

\bibitem{Küttler2008}
U.~K{\"u}ttler, W.~A. Wall, Fixed-point fluid--structure interaction solvers
  with dynamic relaxation, Computational Mechanics 43~(1) (2008) 61--72.
\newline\urlprefix\url{https://doi.org/10.1007/s00466-008-0255-5}

\bibitem{MEHL2016869}
M.~Mehl, B.~Uekermann, H.~Bijl, D.~Blom, B.~Gatzhammer, A.~van Zuijlen,
  Parallel coupling numerics for partitioned fluid–structure interaction
  simulations, Computers \& Mathematics with Applications 71~(4) (2016) 869 --
  891.
\newline\urlprefix\url{http://www.sciencedirect.com/science/article/pii/S0898122115005933}

\bibitem{Wang2018}
Y.~Wang, A.~Quaini, S.~{\v{C}}ani{\'{c}}, A higher-order discontinuous
  galerkin/arbitrary lagrangian eulerian partitioned approach to solving
  fluid--structure interaction problems with incompressible, viscous fluids and
  elastic structures, Journal of Scientific Computing 76~(1) (2018) 481--520.
\newline\urlprefix\url{https://doi.org/10.1007/s10915-017-0629-y}

\bibitem{BASTING2017312}
S.~Basting, A.~Quaini, S.~Čanić, R.~Glowinski, Extended ale method for
  fluid–structure interaction problems with large structural displacements,
  Journal of Computational Physics 331 (2017) 312 -- 336.
\newline\urlprefix\url{http://www.sciencedirect.com/science/article/pii/S0021999116306350}

\bibitem{BadiaNobileVergara2008}
S.~Badia, F.~Nobile, C.~Vergara, Fluid--structure partitioned procedures based
  on {Robin} transmission conditions, J. Comput. Phys. 227~(14) (2008)
  7027--7051.

\bibitem{MokWallRamm2001}
D.~P. Mok, W.~A. Wall, E.~Ramm, Accelerated iterative substructuring schemes
  for instationary fluid structure interaction, in: K.~Bathe (Ed.),
  Computational Fluid and Solid Mechanics, Elsevier, 2001, pp. 1325--1328.

\bibitem{FernandezMullaertVidrascu2014}
M.~A. Fern{\'a}ndez, J.~Mullaert, M.~Vidrascu, {Generalized Robin-Neumann
  explicit coupling schemes for incompressible fluid--structure interaction:
  stability analysis and numerics}, Int. J. Numer. Meth.
  Eng.~\url{http://dx.doi.org/10.1002/nme.4785}.

\bibitem{FernandezMullaertVidrascu2013}
M.~A. Fern{\'a}ndez, J.~Mullaert, M.~Vidrascu, {Explicit Robin-Neumann schemes
  for the coupling of incompressible fluids with thin-walled structures},
  Comput. Method. Appl. Mech. Engrg. 267 (2013) 566--593.

\bibitem{FernandezLandajuela2014}
M.~A. Fern{\'a}ndez, M.~Landajuela, {Fully decoupled time-marching schemes for
  incompressible fluid/thin-walled structure interaction}, Rapport de recherche
  RR-8425, INRIA (Jan. 2014).

\bibitem{Gerardo_GiordaNobileVergara2010}
L.~Gerardo-Giorda, F.~Nobile, C.~Vergara, Analysis and optimization of
  {Robin-Robin} partitioned procedures in fluid--structure interaction
  problems, SIAM J. Numer. Anal. 48~(6) (2010) 2091--2116.

\bibitem{BadiaNobileVergara2009}
S.~Badia, F.~Nobile, C.~Vergara, {Robin-Robin} preconditioned {Krylov} methods
  for fluid--structure interaction problems, CMAME 198~(33--36) (2009)
  2768--2784.

\bibitem{NobilePozzoliVergara2014}
F.~Nobile, M.~Pozzoli, C.~Vergara, Inexact accurate partitioned algorithms for
  fluid--structure interaction problems with finite elasticity in
  haemodynamics, Journal of Computational Physics 273~(0) (2014) 598 -- 617.

\bibitem{fsi2012}
J.~W. Banks, W.~D. Henshaw, D.~W. Schwendeman, Deforming composite grids for
  solving fluid structure problems, J. Comput. Phys. 231~(9) (2012)
  3518--3547\citeCount{34}.

\bibitem{fibr2019}
D.~A. Serino, J.~W. Banks, W.~D. Henshaw, D.~W. Schwendeman, A stable
  added-mass partitioned ({AMP}) algorithm for elastic solids and
  incompressible flow.

\bibitem{splitStep2003}
W.~D. Henshaw, N.~A. Petersson, A split-step scheme for the incompressible
  {Navier-Stokes} equations, in: M.~M. Hafez (Ed.), Numerical Simulation of
  Incompressible Flows, World Scientific, 2003, pp. 108--125\citeCount{50}.

\bibitem{smog2012}
D.~Appel\"o, J.~W. Banks, W.~D. Henshaw, D.~W. Schwendeman, Numerical methods
  for solid mechanics on overlapping grids: Linear elasticity, J. Comput. Phys.
  231~(18) (2012) 6012--6050\citeCount{24}.

\bibitem{flunsi2016}
J.~W. Banks, W.~D. Henshaw, A.~Kapila, D.~W. Schwendeman, An added-mass
  partitioned algorithm for fluid-structure interactions of compressible fluids
  and nonlinear solids, J. Comput. Phys. 305 (2016) 1037--1064\citeCount{8}.

\bibitem{rbinsmp2017}
J.~W. Banks, W.~D. Henshaw, D.~W. Schwendeman, Q.~Tang, A stable partitioned
  {FSI} algorithm for rigid bodies and incompressible flow. {Part I}: Model
  problem analysis., J. Comput. Phys. 343 (2017) 432--468.

\bibitem{DifferenceApproximationsForTheInitialBoundary}
H.~O. Kreiss, Difference approximations for the initial-boundary value problem
  for hyperbolic differential equations, Numerical Solutions of Nonlinear
  Differential Equations (1996) 141–166.

\bibitem{BanksSjogreen2011}
J.~W. Banks, B.~Sj{\"o}green, A normal mode stability analysis of numerical
  interface conditions for fluid/structure interaction, Commun. Comput. Phys.
  10~(2) (2011) 279--304, \url{publications/BanksSjogreenFSI2011.pdf}.

\bibitem{lrb2013}
J.~W. Banks, W.~D. Henshaw, B.~Sj{\"o}green, A stable {FSI} algorithm for light
  rigid bodies in compressible flow, J. Comput. Phys. 245 (2013)
  399--430\citeCount{21}.

\bibitem{GKSII}
B.~Gustafsson, H.-O. Kreiss, A.~Sundstr\"om, Stability theory of difference
  approximations for mixed initial boundary value problems. {II}, Mathematics
  of Computation 26~(119) (1972) 649--686.

\end{thebibliography}
